\documentclass[a4paper,11pt]{article}

\usepackage[pdftex]{hyperref}
\usepackage{stmaryrd}
\usepackage[normalem]{ulem}
\usepackage{graphicx,xcolor}
\usepackage{amsfonts}
\usepackage{amsmath,amssymb,mathrsfs,amsthm}
\usepackage[all]{xy} 
\usepackage{color}
\usepackage{times}
\usepackage{enumerate,vmargin}
\usepackage{verbatim}
\usepackage{paralist}
\usepackage[labelfont=bf,size=small,font=it]{caption}
\usepackage[numbers,sort&compress]{natbib}
\usepackage{bbm}
\usepackage{tikz}
\usetikzlibrary{positioning,decorations.pathreplacing,calligraphy}
\usetikzlibrary{shapes,arrows}

\setpapersize{A4}
\setmarginsrb{2.5cm}{2cm}{2.5cm}{2cm}{.4cm}{4mm}{0cm}{7.5mm}

\hypersetup{
    unicode      = false,     
    pdftoolbar   = true,      
    pdfmenubar   = true,      
    pdffitwindow = true,      
    pdfnewwindow = true,      
    colorlinks   = true,      
    linkcolor    = blue,      
    citecolor    = blue,      
    filecolor    = blue,      
    urlcolor     = blue       
}

\theoremstyle{plain}
\newtheorem{lem}{Lemma}[section]
\newtheorem{thm}[lem]{Theorem}
\newtheorem{prop}[lem]{Proposition}
\newtheorem{cor}[lem]{Corollary}

\theoremstyle{remark}

\newcommand{\RIG}{\mathrm{RIG}}

\newcommand{\bP}{\mathbf P}
\newcommand{\bE}{\mathbf E}

\newcommand{\bQ}{\mathbf Q}

\newcommand{\rH}{\mathrm H}
\newcommand{\Ht}{\mathtt{Hght}}

\newcommand{\rb}{^{(\mathtt{b})}}
\newcommand{\rw}{^{(\mathtt{w})}}
\newcommand{\rbn}{^{(\mathtt{b},n)}}
\newcommand{\rwm}{^{(\mathtt{w},m)}}

\newcommand{\bi}{^{\mathrm{bi}}}

\newcommand{\lf}{\lfloor}
\newcommand{\rf}{\rfloor}

\newcommand{\dhaus}{\operatorname{d^{\mathrm{Haus}}_{E}}}

\newcommand{\dpr}{\operatorname{d_{\mathrm{Pr}}}}

\newcommand{\dghp}{\operatorname{d_{GHP}}}
\newcommand{\dTV}{\operatorname{d_{\mathrm{TV}}}}

\newcommand{\gr}{\operatorname{gr}}
\newcommand{\dgr}{d_{\gr}}

\newcommand{\cT}{\mathcal T}
\newcommand{\cF}{\mathcal F}
\newcommand{\cC}{\mathcal C}

\newcommand{\cS}{\mathcal S}

\newcommand{\cL}{\mathcal L}
\newcommand{\cM}{\mathcal M}
\newcommand{\cA}{\mathcal A}
\newcommand{\cI}{\mathcal I}

\newcommand{\cN}{\mathcal N}
\newcommand{\cD}{\mathcal D}
\newcommand{\cE}{\mathcal E}
\newcommand{\cG}{\mathcal G}

\newcommand{\cZ}{\mathcal Z}
\newcommand{\cH}{\mathcal H}
\newcommand{\cQ}{\mathcal Q}

\newcommand{\bx}{\mathbf x}
\newcommand{\by}{\mathbf y}

\newcommand{\R}{\mathbb R}
\newcommand{\N}{\mathbb N}
\newcommand{\Z}{\mathbb Z}

\newcommand{\eqd}{\stackrel{(d)}=}

\title{\textsc{Random bipartite graphs with i.i.d.~weights\\
and applications to inhomogeneous random intersection graphs}
\thanks{This research is supported by EPSRC grant EP/W033585/1.}
}

\author{Alastair \textsc{Haig}
\thanks{Department of Mathematics, University of Sussex, 
Falmer campus, Brighton, BN1 9QH, England, United Kingdom.  
Email: A.Haig@sussex.ac.uk} \ and\  
Minmin \textsc{Wang}
\thanks{Department of Mathematics, University of Sussex, 
Falmer campus, Brighton, BN1 9QH, England, United Kingdom.  
Email: Minmin.Wang@sussex.ac.uk}}

\date{\today}

\begin{document}

\maketitle
\begin{abstract}

We propose a random bipartite graph with weights assigned to both parts of the vertex sets. Edges are formed independently with probabilities that depend on these weights. This bipartite graph naturally gives rise to a random intersection graph which has nontrivial clustering properties and inhomogeneous vertex degrees. We focus on the situation where the weights are themselves i.i.d.~random variables. In the so-called moderate clustering regime, we identify three types of scaling limit for the large connected components in the graphs at criticality, depending on the tail behaviours of the weight distributions of both parts.   


\smallskip

\end{abstract}

\tableofcontents

\section{Introduction}

Ever since its conception by Erd\H{o}s and R\'enyi in the 1960s, the study of random graphs has continued to yield fascinating discoveries and challenging questions for both combinatorists and probabilists. The current work contributes to a growing body of recent studies that focus on the large-scale structures of random graphs at criticality. Together, these studies point to a compelling picture of universality, wherein models of random graphs belonging to the same universality class are expected to exhibit similar large-scale behaviours within the so-called critical window.
In particular, two notable advances in this direction have been:
\begin{itemize}
\item 
{\bf Identification of the Erd\H{o}s--R\'enyi universality class.} Aldous’ seminal article~\cite{Al97} revealed a surprising connection between the sizes of large connected components in a critical $G(n,p)$ and Brownian excursions. This work also pioneered the use of stochastic analysis in the study of critical random graphs. Addario-Berry, Broutin, and Goldschmidt~\cite{ABBrGo12} deepened this connection, demonstrating that after rescaling edge lengths by a factor of $n^{-1/3}$, one obtains a random metric space ``encoded'' by Brownian excursions for each of these large connected components. The limit graph introduced in~\cite{ABBrGo12}, known as the Erd\H{o}s--R\'enyi continuum graph, has since appeared in the scaling limits of many other models of random graphs; see~\cite{BhvdHvL10, BhSe20, BhDhvdHSe20a, BhBrSeWa14+} to name a few.

\item 
{\bf Discovery of the multiplicative coalescent universality class.} Aldous’ work~\cite{Al97} also highlighted the close relationships between random graphs and the multiplicative coalescent. Together with Aldous and Limic~\cite{AlLi98}, these studies paved the way for discovering more general behaviours beyond those observed in the Erd\H{o}s--R\'enyi graph. To date, the most comprehensive results in this direction concern the rank-1 inhomogeneous random graph model \cite{BoJaRi07, BhSeWa17, BHS18} (also called the Poisson random graph or multiplicative graph). The works~\cite{BrDuWa21} and~\cite{BrDuWa22} address the scaling limits of large components in these graphs under the general asymptotic regime identified in~\cite{AlLi98}. This leads to the introduction of a large family of limit graphs, where the role of Brownian excursions in the construction of the Erd\H{o}s--R\'enyi continuum graph is replaced by excursions of a thinned spectrally positive L\'evy process. Specific members of this family have also been found in the scaling limits of critical configuration models (\cite{BhDhvdHSe20a, CG23}). 
\end{itemize}
In~\cite{W25}, the Binomial bipartite graph model and the associated intersection graphs are investigated. Introduced in~\cite{KaScSC99}, the random intersection graph model provides a simple mechanism for tuning the density of triangles in the graphs, measured by the so-called clustering coefficient. The findings in~\cite{W25} make the case for expanding the Erd\H{o}s--R\'enyi continuum graph into a family of continuum graphs parameterised by a positive number $\theta$, which indicates the clustering level of the discrete graphs. 

While the model of random bipartite graphs in~\cite{W25} can be described as homogeneous, we introduce here a model of {\it inhomogeneous} bipartite graphs by assigning weights to both vertex sets. Our main results, which concern the scaling limit of these graphs at criticality,  give a first indication on how the universality phenomenon extends to random bipartite graphs. 
In particular, these results suggest that, while the weight sequences for the two vertex sets can have different tail behaviours, the large-scale behaviours of the graphs are dominated by the {\it heavier} tail. 
This feature of heavy tail dominance mirrors similar behaviours in the related models, particularly in the multi-type branching trees; see~\cite{Berzunza} and~\cite{HaSt21}.
Key to our approach is a construction of the inhomogeneous bipartite graph via LIFO-queues, based on a similar but simpler version for the rank-1 model that has appeared in~\cite{BrDuWa21, BrDuWa22}. 

\subsection{A model of random bipartite graphs with weights}
\label{sec: model}

Let $n, m\in \N$. Consider a bipartite graph with $n$ black vertices and $m$ white vertices. Suppose that each vertex is further equipped with a {\it weight}, i.e.~a positive number that represents their propensities for forming edges. Denote the respective collections of weights for the black and white vertices as 
\[
\mathbf x = (x_{1}, x_{2}, \dots, x_{n})\in (0, \infty)^n \quad\text{and}\quad \mathbf y=(y_{1}, y_{2}, \dots, y_{m})\in (0, \infty)^m.
\]
The edges in the bipartite graph are formed independently across pairs of vertices of different colours. 
Denote the black vertices as $b_1, b_2, b_3, \dots, b_n$ and the white vertices as $w_1, w_2, \dots, w_m$, and 
write $p_{i, j}$ for the probability of having an edge between $b_i$ and $w_j$. In our model, we assume that 
\begin{equation}
\label{def: edge-prob}
p_{i, j} = 1-\exp\Big(-\frac{x_{i}\cdot y_{j}}{z_{n, m}}\Big), \quad 1\le i\le n, 1\le j\le m,  
\end{equation}
where $z_{n, m}>0$ is some normalising constant that will be specified later on. 
Denote by $B(n,m; \mathbf x, \mathbf y)$ the resulting graph. 

\smallskip
\noindent
{\bf Inhomogeneous random intersection graphs.} The graph $B(n,m; \mathbf x, \mathbf y)$ induces, via the so-called {\it intersection graph} mechanism, a graph on $n$ vertices. Specifically, let $G(n, m; \mathbf x, \mathbf y)$ be the graph on the vertex set $[n]:=\{1, 2, 3, \dots, n\}$ where two vertices $i, j$ share an edge between them if and only if $b_i$ and $b_j$ are at distance $2$ from each other in $B(n,m; \mathbf x, \mathbf y)$. An interesting aspect of the random intersection graphs is their nontrivial {\it clustering} properties, which can be  adjusted by tuning the ratio between $m$ and $n$.  
More precisely, we quantify the clustering level of the graph by looking at the following conditional probability: 
\[
\mathrm{CL} = \mathbb P(V_2 \text{ adjacent to } V_3\,|\, V_1 \text{ adjacent to } V_3 \text{ and }V_2),
\]
where $(V_1, V_2, V_3)$ is a triplet of distinct vertices chosen uniformly. If both $\mathbf x$ and $\mathbf y$ are constant, then it is not difficult to check the following: (i) $\mathrm{CL}\to 0$ if $m/n\to\infty$ and $n\to\infty$; (ii) $\mathrm{CL}\to 1$ if $m/n\to 0$ and $n\to\infty$; (iii) $\mathrm{CL}\to (1+\sqrt{\theta})^{-1}$ if $m=\lfloor \theta n\rfloor$ for some $\theta\in (0, \infty)$ and $n\to\infty$.  Similar phenomena have been observed when $\mathbf x$ is not constant; we refer to~\cite{DeKe09} for further details. In this paper, we will restrict our attention to the {\it moderate clustering regime}; namely, we assume that 
\begin{equation}
\label{hyp: mod-clst}\tag{H-clustering}
\exists\,\theta\in (0, \infty): m=\lfloor \theta n\rfloor.
\end{equation}

\smallskip
\noindent
{\bf Connection with other models.} The expression of $p_{i, j}$ in~\eqref{def: edge-prob} is reminiscent of the rank-1 inhomogeneous random graph model as studied in~\cite{AlLi98, BhvdHvL10, BhvdHvL12, BHS18, BrDuWa21, BrDuWa22}. As in the rank-1 model, vertices with large weights are expected to emerge as ``hubs", i.e.~nodes with high degrees. Compared to the rank-1 model, our random bipartite graphs have two sets of vertex weights, which can exhibit different tail behaviours. How the two weight sets interact with each other and influence the asymptotic behaviours of the graphs is one of the questions that have motivated this work.   

The model of $B(n,m; \mathbf x, \mathbf y)$ and the associated intersection graph $G(n,m; \mathbf x, \mathbf y)$ also encompass models of random graphs that have been introduced previously. In particular:
\begin{itemize}
\item
{\bf The homogeneous case. } This refers to the particular case where both $\mathbf x$ and $\mathbf y$ are constant. Say $x_i\equiv a\in (0, \infty)$ and $y_j\equiv b\in (0, \infty)$. In that case, each of the potential $mn$ edges is present independently with the same probability $p=1-\exp(-ab/z_{n,m})$. If $z_{n,m}=\sqrt{mn}$, then this relates to the Binomial model studied in~\cite{W25}. 
\item 
{\bf The active/passive model. }
If only $\mathbf x$ (resp.~$\mathbf y$) is non constant, this is referred to as the active (i.e.~passive) model in~\cite{BJKR}. Degree distributions and clustering properties for the active model have been studied in~\cite{DeKe09}. See also~\cite{BJKR} and the references therein. 
\end{itemize}
Other models closely related to $G(n, m; \mathbf x, \mathbf y)$ also include the ones studied in~\cite{BaRe11, BaSiTr14, BoJaRi11}. During the preparation of this article, we learned of the works~\cite{Clancy24a, Clancy24b, ClKoLi24} which introduce more general models than $B(n,m; \bx, \by)$. We point out that  the approach in this work is genuinely different and while~\cite{Clancy24a, Clancy24b, ClKoLi24} focus on the component sizes, we prove limit distributions for both the sizes and the graph distances. 

\subsection{The case with i.i.d.~weights}

In this work, we will focus on the case where both sets of weights $\mathbf x, \mathbf y$ are i.i.d.~random variables. The graph $B(n, m; \mathbf x, \mathbf y)$ is then sampled by first conditioning on the weights. Let $F\rb$ and $F\rw$ be the cumulative distribution functions (c.d.f.) of two probability distributions supported on $(0, \infty)$. 
We further assume that both have finite second moments:
\begin{equation}
\label{hyp: 2moment}\tag{H-2nd}
\sigma\rb_{2}:=\int_{(0, \infty)}x^{2}dF\rb(x)<\infty \quad\text{and}\quad \sigma\rw_{2}:=\int_{(0, \infty)}x^{2}dF\rw(x)<\infty,
\end{equation}
which in particular implies the existence of a first moment in both cases. 
Let $(X_{i})_{i\ge 1}$ (resp.~$(Y_{j})_{j\ge 1}$) be a sequence of i.i.d.~variables with common law $dF\rb$ (resp.~$dF\rw$). 
We then set $\mathbf X_{n} = (X_{1}, X_{2}, \dots, X_{n})$ and $\mathbf Y_{m} = (Y_{1}, Y_{2}, \dots, Y_{m})$. We shall use the short-hand notation
\[
B_{n, m} = B(n, m; \mathbf X_{n}, \mathbf Y_{m}), \quad G_{n, m} = G(n, m; \mathbf X_{n}, \mathbf Y_{m})
\]
to denote the random bipartite graph with the vertex weight sequences $\mathbf X_{n}$ and $\mathbf Y_{m}$ and the associated intersection graph. Recall from~\eqref{def: edge-prob} the normalising constant $z_{n, m}$. In the case with i.i.d.~weights, we take 
\[
z_{n, m}=\sqrt{mn},
\]
which amounts to placing $\{G_{n,m}: n, m\ge 1\}$ in a  ``sparse regime''. 
Standard branching process approximation techniques allow us to detect a phase transition in the size of the largest connected component; we refer to~\cite{BoJaRi11} for details. Here, we are interested in the {\it critical} case, i.e.
\begin{equation}
\label{hyp: crit}\tag{H-critical}
\sigma\rb_2\cdot \sigma\rw_2=1.
\end{equation}

\subsection{Main results}

We first gather the various assumptions that have been made on $m, n$ and the weight distributions $F\rb, F\rw$. Recall from~\eqref{hyp: mod-clst} that $m=\lfloor \theta n\rfloor$, which places $G_{n,m}$ in a moderate clustering regime. 
Recall from~\eqref{hyp: 2moment} the assumption of finite second moments and from~\eqref{hyp: crit} the criticality assumption. 
For $r\in (0, \infty)$, denote 
\[
\sigma\rb_r=\int_{(0, \infty)} x^r dF\rb(x)\in [0, \infty] \quad\text{and}\quad \sigma\rw_r=\int_{(0, \infty)} x^r dF\rw(x)\in [0, \infty].
\]
We further assume the following tail behaviours from $F\rb$ and $F\rw$:  either
\begin{equation}
\label{hyp: b-third}\tag{H-b-3rd}
\sigma\rb_3=\int_{(0, \infty)}x^3 dF\rb(x)<\infty,
\end{equation}
or there exists some $\alpha\in (1, 2)$ and $C\rb\in (0,\infty)$ so that 
\begin{equation}
\label{hyp: b-power}\tag{H-b-power}
1-F\rb(x)\sim C\rb x^{-\alpha-1}, \quad x\to\infty; 
\end{equation}
analogously for $F\rw$, either
\begin{equation}
\label{hyp: w-third}\tag{H-w-3rd}
\sigma\rw_3=\int_{(0, \infty)}x^3 dF\rw(x)<\infty,
\end{equation}
or there exists some $\alpha'\in (1, 2)$ and $C\rw\in (0,\infty)$ so that 
\begin{equation}
\label{hyp: w-power}\tag{H-w-power}
1-F\rw(x)\sim C\rw x^{-\alpha'-1}, \quad x\to\infty.
\end{equation}
The symmetry in the model of $B(n,m; \bx, \by)$ means that we can assume without loss of generality that the tail of $F\rw$ is {\it not heavier} than that of $F\rb$. 
The combinations of the two tail behaviours lead to the following three distinctive scenarios: 
\begin{itemize}
\item
{\bf Double finite third moments:} 
this is the case where both~\eqref{hyp: b-third} and~\eqref{hyp: w-third} hold.
\item
{\bf One dominant heavy tail:} 
assume~\eqref{hyp: b-power} holds for $F\rb$ with $\alpha\in(1, 2)$, $C\rb\in (0, \infty)$, and for $F\rw$, either~\eqref{hyp: w-third} is true or~\eqref{hyp: w-power} holds for some $\alpha'\in (\alpha, 2)$ and $C\rw\in (0, \infty)$. 
\item
{\bf Matched heavy tails: }
in this case, ~\eqref{hyp: b-power} and~\eqref{hyp: w-power} both hold with $\alpha'=\alpha\in (1, 2)$, $C\rb\in (0, \infty)$ and $C\rw\in (0, \infty)$.
\end{itemize}

For a connected component $C$ of $B_{n,m}$, we denote by $\bx(C)$ (resp.~$\by(C)$)  the total weights of the black vertices (resp.~white vertices) in $C$, i.e.~$\bx(C)=\sum_{b_i\in C}X_i$ (resp.~$\by(C)=\sum_{w_j\in C}Y_j$). 
Let $\{C^{\bx}_{n, (k)}: 1\le k\le \kappa_{n,m}\}$ be the sequence of the connected components of $B_{n,m}$ that contain at least one edge, ranked in decreasing order of their $\bx$-weights (breaking ties arbitrarily), i.e.~$\bx(C^{\bx}_{n, (1)})\ge \bx(C^{\bx}_{n, (2)})\ge \bx(C^{\bx}_{n, (3)})\ge \cdots$.  We turn each $C^{\bx}_{n, (k)}$ into a metric space by furnishing it with the graph distance $\dgr$ of $B_{n,m}$ and a finite measure $\mu^{\bx}_{n, k}$ which assigns an atom of size $X_i$ to each of the black vertex $b_i$ contained in $C^{\bx}_{n, (k)}$. 
Our main results concern the scaling limit of the measured metric space $(C^{\bx}_{n, (k)},  \dgr, \mu^{\bx}_{n, k})$. 

\begin{thm}[Double finite third moments]
\label{thm: thm1}
Assume~\eqref{hyp: mod-clst} and~\eqref{hyp: crit}. 
Under the assumptions~\eqref{hyp: b-third} and~\eqref{hyp: w-third}, 
there exists a sequence of (random) measured metric spaces $\cG^{(1)}(\theta)=\big\{(\cC^{(1)}_k, d^{(1)}_k, \mu^{(1)}_k): k\in\N\}$ so that 
we have the following convergence in distribution as $n\to\infty$: 
\begin{equation}
\label{eqcv: thm1}
\Big\{\Big(C^{\bx}_{n, (k)}, \tfrac12 n^{-\frac13}\dgr, n^{-\frac23}\mu^{\bx}_{n,k}\Big): 1\le k\le \kappa_{n,m}\Big\} \ \Longrightarrow\  \cG^{(1)}(\theta) 
\end{equation}
with respect to the product topology induced by the Gromov--Hausdorff--Prokhorov topology. 
\end{thm}

\begin{thm}[One dominant heavy tail]
\label{thm: thm2}
Assume~\eqref{hyp: mod-clst} and~\eqref{hyp: crit}.  
Suppose that~\eqref{hyp: b-power} holds for $F\rb$ with $\alpha\in(1, 2)$, $C\rb\in (0, \infty)$, and that either~\eqref{hyp: w-third} is true or~\eqref{hyp: w-power} holds for some $\alpha'\in (\alpha, 2)$ and $C\rw\in (0, \infty)$.  Then there exists a sequence of (random) measured metric spaces $\cG^{(2)}(\theta)=\big\{(\cC^{(2)}_k, d^{(2)}_k, \mu^{(2)}_k): k\in\N\}$ so that 
we have the following convergence in distribution as $n\to\infty$: 
\begin{equation}
\label{eqcv: thm2}
\Big\{\Big(C^{\bx}_{n, (k)}, \tfrac12 n^{-\frac{1}{\alpha+1}}\dgr, n^{-\frac{\alpha}{\alpha+1}}\mu^{\bx}_{n,k}\Big): 1\le k\le \kappa_{n,m}\Big\} \ \Longrightarrow\  \cG^{(2)}(\theta) 
\end{equation}
with respect to the product topology induced by the Gromov--Hausdorff--Prokhorov topology. 

\end{thm}

\begin{thm}[Matched heavy tails]
\label{thm: thm3}
Assume~\eqref{hyp: mod-clst} and~\eqref{hyp: crit}.  
Suppose that~\eqref{hyp: b-power} and~\eqref{hyp: w-power} both hold with $\alpha'=\alpha\in (1, 2)$, $C\rb\in (0, \infty)$ and $C\rw\in (0, \infty)$. Then there exists a sequence of (random) measured metric spaces $\cG^{(3)}(\theta)=\big\{(\cC^{(3)}_k, d^{(3)}_k, \mu^{(3)}_k): k\in\N\}$ so that 
we have the following convergence in distribution as $n\to\infty$: 
\begin{equation}
\label{eqcv: thm1}
\Big\{\Big(C^{\bx}_{n, (k)}, \tfrac12 n^{-\frac{1}{\alpha+1}}\dgr, n^{-\frac{\alpha}{\alpha+1}}\mu^{\bx}_{n,k}\Big): 1\le k\le \kappa_{n,m}\Big\} \ \Longrightarrow\  \cG^{(3)}(\theta) 
\end{equation}
with respect to the product topology induced by the Gromov--Hausdorff--Prokhorov topology. 
\end{thm}

We can also rank the connected components in their $\by$-weights: let $\{C^{\by}_{n, (k)}: 1\le k\le \kappa_{n,m}\}$ be the re-ordering of $\{C^{\bx}_{n, (k)}: 1\le k\le\kappa_{n,m}\}$ in decreasing order of their $\by$-weights. 
Denote by $\mu^{\by}_{n, k}$ the finite measure on $C^{\bx}_{n, (k)}$ which assigns an atom of size $Y_j$ to each of the white vertex $w_j$ in $C^{\bx}_{n, (k)}$. 
 The following result says that with high probability, $C^{\bx}_{n, (k)}$ coincides with $C^{\by}_{n, (k)}$ and $\bx(C^{\bx}_{n, (k)})/\by(C^{\bx}_{n, (k)})\to \sigma\rb_2/\sqrt\theta$ in probability for each $k\ge 1$. 
 
\begin{prop}[Consistency in rankings]
\label{prop: rnk}
Let $\rho= \sigma\rb_2/\sqrt\theta$. 
Assume the conditions from one of the three Theorems~\ref{thm: thm1}-~\ref{thm: thm3}. Then for any $K\ge 1$, 
\[
\lim_{n\to\infty}\mathbb P\Big(\forall\, 1\le k\le K: C^{\by}_{n, (k)}=C^{\bx}_{n, (k)} \Big) = 1.
\]
Moreover, the following weak convergences take place with respect to the product topology induced by the Gromov--Hausdorff--Prokhorov topology: under the conditions of Theorems~\ref{thm: thm1}, we have
\[
\Big\{\Big(C^{\bx}_{n, (k)}, \tfrac12 n^{-\frac13}\dgr, n^{-\frac23}\rho^{-1}\mu^{\by}_{n,k}\Big): 1\le k\le \kappa_{n,m}\Big\} \ \Longrightarrow\  \cG^{(1)}(\theta); 
\]
under the conditions of Theorem~\ref{thm: thm2}, we have
\[
\Big\{\Big(C^{\bx}_{n, (k)}, \tfrac12 n^{-\frac{1}{\alpha+1}}\dgr, n^{-\frac{\alpha}{\alpha+1}}\rho^{-1}\mu^{\by}_{n,k}\Big): 1\le k\le \kappa_{n,m}\Big\} \ \Longrightarrow\  \cG^{(2)}(\theta); 
\]
under the conditions of Theorem~\ref{thm: thm3}, we have
\[
\Big\{\Big(C^{\bx}_{n, (k)}, \tfrac12 n^{-\frac{1}{\alpha+1}}\dgr, n^{-\frac{\alpha}{\alpha+1}}\rho^{-1}\mu^{\by}_{n,k}\Big): 1\le k\le \kappa_{n,m}\Big\} \ \Longrightarrow\  \cG^{(3)}(\theta). 
\]
\end{prop}

As in the cases of Erd\H{o}s--R\'enyi model~\cite{ABBrGo12}, the rank-1 models~\cite{BrDuWa21, BrDuWa22}, the configuration models~\cite{CG23} and a number of other notables random graph models, the limit graphs that appear in the previous theorems can be constructed from certain stochastic processes. We refer to Section~\ref{sec: limit-gr} for the construction of $\cG^{(i)}(\theta), i\in \{1, 2, 3\}$. 

\medskip

Let us now explore the implications on the intersection graph $G_{n,m}$. Denote by $\varrho_n$ the mapping that sends the vertex $b_i$ in the bipartite graph $B_{n,m}$ to the vertex $i$ in $G_{n,m}$, and denote by $\hat C_{n, (k)}$ the image of $C^{\bx}_{n, (k)}$ by $\varrho_n$. Then $\hat C_{n, (k)}$ is a connected component of $G_{n,m}$. 
Denote by $\dgr^{\RIG}$ the graph distance of $G_{n,m}$ and let $\mu_{n,k}^{\RIG}$ be the push-forward of $\mu^{\bx}_{n, k}$ by $\varrho_n$.  
In~\cite{W25}, a simple argument (Proposition 2.4) shows that $\varrho_n$ is in fact an isometry between $(B_{n,m}, \dgr)$ and $(G_{n,m}, 2\cdot \dgr^{\RIG})$, and as a result  the Gromov--Hausdorff--Prokhorov distance between $(\hat C_{n, (k)}, \dgr^{\RIG}, \mu_{n,k}^{\RIG})$ and $(C^{\bx}_{n,(k)}, \tfrac12\dgr, \mu^{\bx}_{n,k})$ is bounded by $1$. For the current model, this bound combined with 
the previous theorems immediately yields the following results on the intersection graph $G_{n,m}$. 

\begin{cor}
Assume~\eqref{hyp: mod-clst} and~\eqref{hyp: crit}.  
\begin{enumerate}
\item
{\sf (Double finite third moments)} Under the assumptions~\eqref{hyp: b-third} and~\eqref{hyp: w-third}, we have the following convergence in distribution as $n\to\infty$: 
\[
\Big\{\Big(\hat C_{n, (k)}, n^{-\frac13}\dgr^{\RIG}, n^{-\frac23}\mu^{\RIG}_{n,k}\Big): 1\le k\le\kappa_{n,m}\Big\} \ \Longrightarrow\  \cG^{(1)}(\theta) 
\]
with respect to the product topology induced by Gromov--Hausdorff--Prokhorov topology. 
\item
{\sf (One dominant heavy tail)} Suppose that~\eqref{hyp: b-power} holds for $F\rb$ with $\alpha\in(1, 2)$, $C\rb\in (0, \infty)$, and that either~\eqref{hyp: w-third} is true or~\eqref{hyp: w-power} holds for some $\alpha'\in (\alpha, 2)$ and $C\rw\in (0, \infty)$.    We have the following convergence in distribution as $n\to\infty$: 
\[
\Big\{\Big(\hat C_{n, (k)}, n^{-\frac{1}{\alpha+1}}\dgr^{\RIG}, n^{-\frac{\alpha}{\alpha+1}}\mu^{\RIG}_{n,k}\Big): 1\le k\le\kappa_{n,m}\Big\} \ \Longrightarrow\  \cG^{(2)}(\theta) 
\]
with respect to the product topology induced by Gromov--Hausdorff--Prokhorov topology. 
\item
{\sf (Matched heavy tails)} Suppose that~\eqref{hyp: b-power} and~\eqref{hyp: w-power} both hold with $\alpha'=\alpha\in (1, 2)$, $C\rb\in (0, \infty)$ and $C\rw\in (0, \infty)$.  We have the following convergence in distribution as $n\to\infty$: 
\[
\Big\{\Big(\hat C_{n, (k)}, n^{-\frac{1}{\alpha+1}}\dgr^{\RIG}, n^{-\frac{\alpha}{\alpha+1}}\mu^{\RIG}_{n,k}\Big): 1\le k\le\kappa_{n,m}\Big\} \ \Longrightarrow\  \cG^{(3)}(\theta) 
\]
with respect to the product topology induced by Gromov--Hausdorff--Prokhorov topology. 
\end{enumerate}
\end{cor}

\smallskip
\noindent
{\bf Discussion.} We have restricted to the moderate clustering regime, namely~\eqref{hyp: mod-clst}. It is possible to consider the cases where $m/n\to\infty$ or $m/n\to 0$. For the homogeneous case, the scaling limit of graph across the full clustering spectrum was investigated in~\cite{W25} and the results there revealed that the diameter of the graph can vary according to the clustering. For the inhomogeneous model, the situation seems to be more complex; we therefore leave it as an open question.

\smallskip
\noindent
{\bf Organisation of the paper.} Our main results are shown by studying the scaling limit of certain stochastic processes that encode $B_{n,m}$. In Section~\ref{sec: tech}, we introduce these graph encoding processes and identify their  limit. We also explain how to obtain the limit graphs from these limit processes. Detailed proofs are found in Section~\ref{sec: proof}. 

\section{Exploration of the bipartite graphs} 
\label{sec: tech}

Graph exploration has become a standard and versatile tool in the study of random graphs. In our current approach, we will rely on a Last-In-First-Out (LIFO) queue that will simultaneously generate the random bipartite graph $B_{n, m}$ and a depth-first traversal of the graph. 
The idea of constructing graphs using LIFO-queues can be traced back to~\cite{BrDuWa21, BrDuWa22}, where a simpler version appears in the context of rank-1 models, and further back to~\cite{LGLJ98} for the genealogies of branching processes.

\subsection{LIFO-queues and depth-first exploration of graphs}
\label{sec: explore}

We begin with the description of a generic LIFO-queue, and explain how it leads to a natural notion of genealogy. The focus of this subsection is an alternative construction of the random bipartite $B(n,m; \bx, \by)$, using two independent sequences of exponential variables.  We then explain how the construction is connected with LIFO-queues.  

\smallskip
\noindent
{\bf A Last-In-First-Out queuing system. } Let us consider a queuing system with a single server, where a total number of $N\in \N\cup\{\infty\}$ clients arrive during the whole process.  Let $i\in \N$ and $i\le N$. Suppose that Client $i$ arrives at time $t_i$ and requests $\Delta_i\in [0, \infty)$ amount of service time. The client is served immediately, even if the server has been occupied by another client prior to its arrival. Once Client $i$ has received in full the $\Delta_i$ service time, noting that its own service can be interrupted by later arrivals, the server returns to the client whose service has been interrupted by the arrival of Client $i$. 

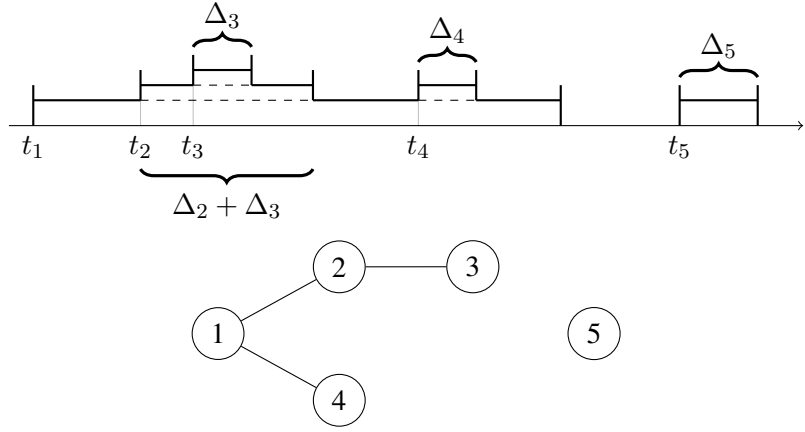
\begin{figure}
\centering
\begin{tikzpicture}
\draw[->] (0,0)--(10.5,0);

\node[below] (t1) at (0.322,0) {$t_1$};
\node[below] (t2) at (1.740,0) {$t_2$};
\node[below] (t3) at (2.436,0) {$t_3$};
\node[below] (d3) at (3.208,0) {};
\node[below] (d2) at (4.020,0) {};
\node[below] (t4) at (5.414,0) {$t_4$};
\node[below] (d4) at (6.178,0) {};
\node[below] (d1) at (7.296,0) {};
\node[below] (t5) at (8.868,0) {$t_5$};
\node[below] (d5) at (9.9,0) {};

\node (t1t) [above = 0.4cm of t1] {};
\node (t1m) [above = 0.2cm of t1] {};
\node (t2a) [above = 0.2cm of t2] {};
\node (t2b) [above = 0.4cm of t2] {};
\node (t2c) [above = 0.6cm of t2] {};
\node (t2d) [above = 0.8cm of t2] {};
\node (t3a) [above = 0.2cm of t3] {};
\node (t3b) [above = 0.4cm of t3] {};
\node (t3c) [above = 0.6cm of t3] {};
\node (t3d) [above = 0.8cm of t3] {};
\node (d3a) [above = 0.2cm of d3] {};
\node (d3b) [above = 0.4cm of d3] {};
\node (d3c) [above = 0.6cm of d3] {};
\node (d3d) [above = 0.8cm of d3] {};
\node (d2a) [above = 0.2cm of d2] {};
\node (d2b) [above = 0.4cm of d2] {};
\node (d2c) [above = 0.6cm of d2] {};
\node (d2d) [above = 0.8cm of d2] {};
\node (t4a) [above = 0.2cm of t4] {};
\node (t4b) [above = 0.4cm of t4] {};
\node (t4c) [above = 0.6cm of t4] {};
\node (d4a) [above = 0.2cm of d4] {};
\node (d4b) [above = 0.4cm of d4] {};
\node (d4c) [above = 0.6cm of d4] {};
\node (d1a) [above = 0.2cm of d1] {};
\node (d1b) [above = 0.4cm of d1] {};
\node (t5a) [above = 0.2cm of t5] {};
\node (t5b) [above = 0.4cm of t5] {};
\node (d5a) [above = 0.2cm of d5] {};
\node (d5b) [above = 0.4cm of d5] {};

\draw[thick] (t1)--(t1t.center);
\draw[thick] (t2a.center)--(t2c.center);
\draw[thick] (t3b.center)--(t3d.center);
\draw[thick] (d3b.center)--(d3d.center);
\draw[thick] (d2a.center)--(d2c.center);
\draw[thick] (t4a.center)--(t4c.center);
\draw[thick] (d4a.center)--(d4c.center);
\draw[thick] (d1)--(d1b.center);
\draw[thick] (d5)--(d5b.center);
\draw[thick] (t5)--(t5b.center);

\draw[thick] (t1m.center)--(t2a.center);
\draw[thick] (t2b.center)--(t3b.center);
\draw[thick] (t3c.center)--(d3c.center);
\draw[thick] (d3b.center)--(d2b.center);
\draw[thick] (d2a.center)--(t4a.center);
\draw[thick] (t4b.center)--(d4b.center);
\draw[thick] (d4a.center)--(d1a.center);
\draw[thick] (t5a.center)--(d5a.center);

\draw[dashed] (t2a.center) -- (d2a.center);
\draw[dashed] (t3b.center) -- (d3b.center);
\draw[dashed] (t4a.center) -- (d4a.center);

\draw[opacity = 0.25] (t2a.center)--(t2);
\draw[opacity = 0.25] (t3b.center)--(t3);
\draw[opacity = 0.25] (t4a.center)--(t4);

\draw[decorate, very thick,
    decoration = {brace,mirror,
        raise=-5pt,
        amplitude=5pt}] (1.740,-0.75)--(4.020,-0.75) node[pos=0.5,below,black]{$\Delta_2+\Delta_3$};
\draw[decorate, very thick,
    decoration = {brace,
        raise=3pt,
        amplitude=5pt}] (t3d.center)--(d3d.center) node[pos=0.5,above=6pt,black]{$\Delta_3$};
\draw[decorate, very thick,
    decoration = {brace,
        raise=3pt,
        amplitude=5pt}] (t4c.center)--(d4c.center) node[pos=0.5,above=6pt,black]{$\Delta_4$};
\draw[decorate, very thick,
    decoration = {brace,
        raise=3pt,
        amplitude=5pt}] (t5b.center)--(d5b.center) node[pos=0.5,above=6pt,black]{$\Delta_5$};

\node[] (rootbase) at (5.25,-2.75) {};
\node[draw,circle] (21) [above right = 0.5cm and 0.5cm of rootbase] {3};
\node[draw,circle] (11) [above left = 0.5cm and 0.5cm of rootbase] {2};
\node[draw,circle] (12) [below left = 0.5cm and 0.5cm of rootbase] {4};
\node[draw,circle] (1) [left = 2cm of rootbase] {1};
\node[draw,circle] (5) [right = 2cm of rootbase] {5};

\draw (1) -- (11);
\draw (11) -- (21);
\draw (1) -- (12);
\end{tikzpicture}
\caption{\label{fig: LIFO} An example of a LIFO queue and its genealogy. There are five clients in total. Client 1's service is interrupted by the arrival of Client 2. During the service of the latter, Client 3 arrives. This corresponds to the branch 1-2-3 in the genealogy below. In due time, Client 3 leaves the queue, and the server gets back to Client 2, completing its remaining service. After that, the service of Client 1 is resumed, but interrupted again by the arrival of Client 4. Only after the departure of Client 4, the service of Client 1 is finally completed. The server is then unoccupied for a while until the final client arrives. }
\end{figure}

\smallskip
\noindent
{\bf Genealogy of a LIFO-queue. } Given the previous queue with $N$ clients, we say that Client $i$ is a child of Client $j$ if and only if the arrival of the former interrupts the service of the latter. In particular, if a client arrives when the server is unoccupied, then that client becomes an ancestor for a subset of the $N$ clients that consists of those arriving while that client is still in the queue. 
It is not difficult to see that this notion of genealogy can be represented as a forest of {\it rooted ordered trees} on $N$ vertices where each vertex stands for a client and siblings are ranked according to their arrivals. 
We point to Fig.~\ref{fig: LIFO} for an example. 
We note that the arrival order of the clients in the queue coincides with the depth-first traversal of this forest. 

\smallskip
\noindent
{\bf A LIFO-queue construction of $B(n,m; \bx, \by)$. } We recall the random bipartite graph $B(n,m;\bx, \by)$ with given weight sequences $\bx=(x_i)_{1\le i\le n}$ and $\by=(y_j)_{1\le j\le m}$. We now present an alternative construction of the graph where the connected components of the graph appear in a sequential way. 
Informally, the construction consists in sampling the neighbourhood of each vertex in turn. However, the particular schemes for sampling depend on the colour of the vertex. For black vertices, this is done in a breadth-first manner, i.e.~the entire neighbourhood is sampled in one step; for white vertices, depth-first samplings are used instead:  
each time a member of the neighbourhood has been identified, we immediately move on to the neighbourhood of this member. 
Formally, let 
\[
\big\{E\rb_{i}: 1\le i\le n\big\}  \quad \text{and} \quad \big\{E\rw_j: 1\le j\le m\big\}
\]
be two independent collections of independent exponential variables of respective rates $x_{i}/z_{n,m}$, $1\le i\le n$, and $y_{j}/z_{n,m}$, $1\le j\le m$. 
During the construction, we operate a queue for the white vertices. At step $k\ge 0$, the status of the queue is recorded in an ordered sequence $\cA_k=\{(l_{k,i}, r_k(l_{k, i})): 1\le i\le h_k\}$, where $l_{k, 1}, l_{k, 2}, l_{k, 3}, \dots$ are the labels of the first, second, third, and so on white vertices  in the queue, and $r_k(l_{k, i})$ is the remaining service time for $l_{k, i}$ at that stage.
Imagine also two dials--a black one and a white one--that will help us to record the progress of the construction. 
Initially, set the black dial $\tau\rb_{0}=0$ and the white dial $\tau\rw_{0}=0$. Set the initial queue $\cA_0$ to be empty. At step $k\ge 1$, given $\tau\rb_{k-1}, \tau\rw_{k-1}$ and $\cA_{k-1}$, do as follows. 

\smallskip
\noindent
{\bf If $\cA_{k-1}$ is empty, } this signifies the start of a new connected component. Let $V_k=\mathrm{argmin}\{E\rb_i: E\rb_i>\tau\rb_{k-1}\}$, namely the first black vertex arriving after $\tau\rb_{k-1}$. Set 
\[
\tau\rb_k=E\rb_{V_k}, \quad \text{and}\quad I(V_k)=(\tau\rw_{k-1}, \tau\rw_k]  \quad\text{with}\quad  \tau\rw_k:=\tau\rw_{k-1}+x_{V_k}.
\] 
Declare $V_k$ as the root of a new rooted tree whose first generation, consisting entirely of white vertices, is given by the following set
\[
\cN(V_k):=\{w_j: j \in [m], E\rw_j\in I(V_k)\}.
\]
Update the queue as follows: rank the elements of $\cN(V_k)$ in increasing order of their arrivals, and put them to the top of $\cA_{k-1}$, along with a service request $r_k(w_j) = y_j$ for each $w_j\in \cN(V_k)$. Call the new queue $\cA_k$. 

\smallskip
\noindent
{\bf If $\cA_{k-1}$ is not empty, } let $U_k$ be the first entry in $\cA_{k-1}$ along with its remaining service time $r_{k-1}(U_k)$. This signifies that $U_k$ is the customer being served from $\tau\rb_{k-1}$ until the moment 
\[
\tau\rb_k := \big(\tau\rb_{k-1} + r_{k-1}(U_k)\big) \wedge \min\big\{E\rb_i: E\rb_i>\tau\rb_{k-1}\big\}.
\]
In words, $U_k$ will receive its service in full unless interrupted by a new arrival. 
If $\tau\rb_k - \tau\rb_{k-1}<r_{k-1}(U_k)$, then a new black vertex arrives and breaks the service of $U_k$.  Let $V_k=\mathrm{argmin}\{E\rb_i: E\rb_i>\tau\rb_{k-1}\}$ be this  black vertex and 
declare $V_k$ as a child of $U_k$. Set 
\[
I(V_k) = (\tau\rw_{k-1}, \tau\rw_k] \quad\text{with} \quad \tau\rw_k:=\tau\rw_{k-1}+x_{V_k}.
\]
The offspring of $V_k$ is defined as
\[
\cN(V_k) = \big\{w_j: j\in [m], E\rw_j\in I(V_k)\big\}.
\]
Update the queue as follows: first, set the remaining service time for $U_k$ to be $r_k(U_k) = r_{k-1}(U_k)-(E\rb_{V_k}-\tau\rb_{k-1})$; secondly, rank the elements of $\cN(V_k)$ in increasing order of their arrivals and add them to the top of $\cA_{k-1}$,  along with a service request $r_k(w_j)=y_j$ for each $w_j\in \cN(V_k)$. Any other previous member of $\cA_{k-1}$ has their remaining service time unchanged. Call the new queue $\cA_k$. 

\noindent
If, instead, $\tau\rb_k - \tau\rb_{k-1}= r_{k-1}(U_k)$, then the server completes the service of $U_k$ at $\tau\rb_k$. In that case, remove $U_k$ from $\cA_{k-1}$ to  yield $\cA_k$ and set $\tau\rw_k=\tau\rw_{k-1}$. 

\smallskip
\noindent
{\bf Stop} when $\cA_{k-1}$ is empty and $\tau\rb_{k-1}> \max_{i\in [n]}E\rb_i$. 

\medskip

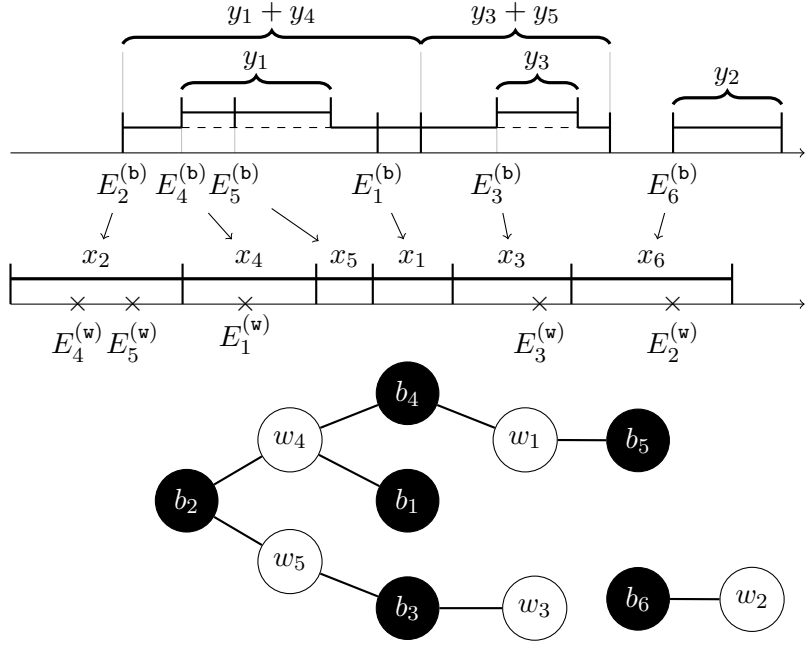
\begin{figure}
    \centering
    \begin{tikzpicture}
        \draw[->] (0,0)--(10.5,0);
        
        \node[below] (b1) at (1.484,0) {$E_2\rb$};
        \node (b1b) [above = 0.2cm of b1] {};
        \node (b1c) [above = 0.4cm of b1] {};
        \node (b1e) [above = 1.2cm of b1] {};
        \node[below] (b2) at (2.262,0) {$E_4\rb$};
        \node (b2b) [above = 0.2cm of b2] {};
        \node (b2c) [above = 0.4cm of b2] {};
        \node (b2d) [above = 0.6cm of b2] {};
        \node[below] (b3) at (2.965,0) {$E_5\rb$};
        \node (b3b) [above = 0.2cm of b3] {};
        \node (b3c) [above = 0.4cm of b3] {};
        \node (b3d) [above = 0.6cm of b3] {};
        \node[below] (w3) at (4.235,0) {};
        \node (w3b) [above = 0.2cm of w3] {};
        \node (w3c) [above = 0.4cm of w3] {};
        \node (w3d) [above = 0.6cm of w3] {};
        \node[below] (b4) at (4.854,0) {$E_1\rb$};
        \node (b4b) [above = 0.2cm of b4] {};
        \node (b4c) [above = 0.4cm of b4] {};
        \node (b4d) [above = 0.6cm of b4] {};
        \node[below] (w1) at (5.426,0) {};
        \node (w1b) [above = 0.2cm of w1] {};
        \node (w1c) [above = 0.4cm of w1] {};
        \node (w1d) [above = 0.6cm of w1] {};
        \node (w1e) [above = 1.2cm of w1] {};
        \node[below] (b5) at (6.43,0) {$E_3\rb$};
        \node (b5b) [above = 0.2cm of b5] {};
        \node (b5c) [above = 0.4cm of b5] {};
        \node (b5d) [above = 0.6cm of b5] {};
        \node[below] (w4) at (7.500,0) {};
        \node (w4b) [above = 0.2cm of w4] {};
        \node (w4c) [above = 0.4cm of w4] {};
        \node (w4d) [above = 0.6cm of w4] {};
        \node[below] (w2) at (7.926,0) {};
        \node (w2b) [above = 0.2cm of w2] {};
        \node (w2c) [above = 0.4cm of w2] {};
        \node (w2e) [above = 1.2cm of w2] {};
        \node[below] (b6) at (8.760,0) {$E_6\rb$};
        \node (b6b) [above = 0.2cm of b6] {};
        \node (b6c) [above = 0.4cm of b6] {};
        \node[below] (w5) at (10.196,0) {};
        \node (w5b) [above = 0.2cm of w5] {};
        \node (w5c) [above = 0.4cm of w5] {};

        \draw[thick] (b1)--(b1c.center);
        \draw[thick] (b2b.center)--(b2d.center);
        \draw[thick] (b3b.center)--(b3d.center);
        \draw[thick] (w3b.center)--(w3d.center);
        \draw[thick] (b4)--(b4c.center);
        \draw[thick] (w1)--(w1c.center);
        \draw[thick] (b5b.center)--(b5d.center);
        \draw[thick] (b6)--(b6c.center);
        \draw[thick] (w4b.center)--(w4d.center);
        \draw[thick] (w2)--(w2c.center);
        \draw[thick] (w5)--(w5c.center);

        \draw[thick] (b1b.center)--(b2b.center);
        \draw[thick] (b2c.center)--(w3c.center);
        \draw[thick] (w3b.center)--(b5b.center);
        \draw[thick] (b5c.center)--(w4c.center);
        \draw[thick] (b6b.center)--(w5b.center);
        \draw[thick] (w4b.center)--(w2b.center);

        \draw[dashed] (b2b.center) -- (w3b.center);
        \draw[dashed] (b5b.center) -- (w4b.center);

        \draw[opacity = 0.25] (b3b.center)--(b3);
        \draw[opacity = 0.25] (b2b.center)--(b2);
        \draw[opacity = 0.25] (b5b.center)--(b5);

        \draw[opacity = 0.25] (b1c.center)--(b1e.center);
        \draw[opacity = 0.25] (w1c.center)--(w1e.center);
        \draw[opacity = 0.25] (w2c.center)--(w2e.center);

        \draw[decorate, very thick,
        decoration = {brace,
        raise=3pt,
        amplitude=5pt}] (b1e.center)--(w1e.center) node[pos=0.5,above=6pt,black]{$y_1+y_4$};
        \draw[decorate, very thick,
        decoration = {brace,
        raise=3pt,
        amplitude=5pt}] (b2d.center)--(w3d.center) node[pos=0.5,above=6pt,black]{$y_1$};
        \draw[decorate, very thick,
        decoration = {brace,
        raise=3pt,
        amplitude=5pt}] (b6c.center)--(w5c.center) node[pos=0.5,above=6pt,black]{$y_2$};
        \draw[decorate, very thick,
        decoration = {brace,
        raise=3pt,
        amplitude=5pt}] (w1e.center)--(w2e.center) node[pos=0.5,above=6pt,black]{$y_3+y_5$};
        \draw[decorate, very thick,
        decoration = {brace,
        raise=3pt,
        amplitude=5pt}] (b5d.center)--(w4d.center) node[pos=0.5,above=6pt,black]{$y_3$};

        \draw[->] (0,-2)--(10.5,-2);

        \node[below = 0.1] (wA1) at (0.889,-2) {$E_4\rw$};
        \node[] at (0.889,-2) {$\times$};
        \node[below = 0.1] (wA2) at (1.616,-2) {$E_5\rw$};
        \node[] at (1.616,-2) {$\times$};
        \node[below] (bA1) at (2.275,-2) {};
        \node (bA1b) [above = 0.4cm of bA1] {};
        \node (bA1m) [above = 0.2cm of bA1] {};
        \node[below] (wA3) at (3.107,-2) {$E_1\rw$};
        \node[] at (3.107,-2) {$\times$};
        \node (wA3b) [above = 0.4cm of wA3] {};
        \node (wA3m) [above = 0.2cm of wA3] {};
        \node[below] (bA2) at (4.042,-2) {};
        \node (bA2m) [above = 0.2cm of bA2] {};
        \node (bA2b) [above = 0.4cm of bA2] {};
        
        \node[below] (bA3) at (4.79,-2) {};
        \node (bA3m) [above = 0.2cm of bA3] {};
        \node (bA3b) [above = 0.4cm of bA3] {};
        
        \node[below] (bA4) at (5.846,-2) {};
        \node (bA4m) [above = 0.2cm of bA4] {};
        \node (bA4b) [above = 0.4cm of bA4] {};
        \node[below = 0.1] (wA4) at (7,-2) {$E_3\rw$};
        \node[] at (7,-2) {$\times$};
        \node[below] (bA5) at (7.416,-2) {};
        \node (bA5b) [above = 0.4cm of bA5] {};
        \node[below = 0.1] (wA5) at (8.758,-2) {$E_2\rw$};
        \node[] at (8.758,-2) {$\times$};
        \node (bA5m) [above = 0.2cm of bA5] {};
        \node[below] (bA6) at (9.541,-2) {};
        \node (bA6b) [above = 0.4cm of bA6] {};
        \node (bA6m) [above = 0.2cm of bA6] {};

        \node[] (root2) at (0,-2) {};
        \node (root2b) [above = 0.4cm of root2.center] {};
        \node[] (root2m) [above = 0.2cm of root2.center] {};

        \draw[thick] (root2.center)--(root2b.center);
        \draw[thick] (bA1)--(bA1b.center);
        \draw[thick] (bA2)--(bA2b.center);
        \draw[thick] (bA3)--(bA3b.center);
        \draw[thick] (bA4)--(bA4b.center);
        \draw[thick] (bA5)--(bA5b.center);
        \draw[thick] (bA6)--(bA6b.center);

        \draw[very thick] (root2m.center) -- (bA1m.center) node[midway,above] (xV1){$x_{2}$};
        \draw[very thick] (bA2m.center) -- (bA1m.center) node[midway,above] (xV2){$x_{4}$};
        \draw[very thick] (bA3m.center) -- (bA2m.center) node[midway,above] (xV3){$x_{5}$};
        \draw[very thick] (bA3m.center) -- (bA4m.center) node[midway,above] (xV4){$x_{1}$};
        \draw[very thick] (bA4m.center) -- (bA5m.center) node[midway,above] (xV5){$x_{3}$};
        \draw[very thick] (bA5m.center) -- (bA6m.center) node[midway,above] (xV6){$x_{6}$};

        \draw[->] (b1)--(xV1); 
        \draw[->] (b2)--(xV2);
        \draw[->] (b3)--(xV3);
        \draw[->] (b4)--(xV4);
        \draw[->] (b5)--(xV5);
        \draw[->] (b6)--(xV6); 

        \node[draw, circle, fill = black, text = white] (Tb4) at (5.25,-4.6) {$b_1$};
        \node[draw, circle, fill = black, text = white] (Tb5) [below = 1cm of Tb4.center] {$b_3$};
        \node[draw, circle, fill = black, text = white] (Tb2) [above = 1cm of Tb4.center] {$b_4$};

        \node[draw, circle, fill = black, text = white] (Tb1) [left = 2.5cm of Tb4.center] {$b_2$};

        \node[draw, circle] (Tw2) [below left = 0.5cm and 1.25cm of Tb4.center] {$w_5$};
        \node[draw, circle] (Tw1) [above left = 0.5cm and 1.25cm of Tb4.center] {$w_4$};
        \node[draw, circle] (Tw3) [above right = 0.5cm and 1.25cm of Tb4.center] {$w_1$};

        \node[draw, circle, fill = black, text = white] (Tb3) [above right = 0.5cm and 2.75cm of Tb4.center] {$b_5$};

        \node[draw, circle] (Tw4) [right = 1.25cm of Tb5.center] {$w_3$};

        \node[draw, circle, fill = black, text = white] (Tb6) [below right = 1cm and 2.75cm of Tb4.center] {$b_6$};

        \node[draw, circle] (Tw5) [below right = 1cm and 4.25cm of Tb4.center] {$w_2$};
        
        \draw[thick] (Tw4)--(Tb5)--(Tw2)--(Tb1)--(Tw1);
        \draw[thick] (Tb4)--(Tw1)--(Tb2)--(Tw3)--(Tb3);
        \draw[thick] (Tb6)--(Tw5);
        
    \end{tikzpicture}
    \caption{\label{fig: bi}A LIFO-queue construction of $B(n,m; \bx, \by)$. The first steps of the construction can be explained as follows: since $E\rb_2$ is the smallest among all the $E\rb_i$'s, we have $V_1=2$ and start the first tree rooted at $b_2$. An inspection of the line below finds the two neighbours of $b_2$: $w_4$ and $w_5$; thus $\cA_1=\{(w_4, y_4), (w_5, y_5)\}$. At step $2$, we start with $U_2=w_4$. The service of $w_4$ gets interrupted by the arrival of $b_4$ at $E\rb_4$, which becomes the first child of $w_4$. We then immediately pivot to identifying the neighbours of $b_4$ and find $w_1$ as a result. At the end of step $2$, the queue status is given by
$\cA_2 = \{(w_1, y_1), (w_4, y_4-(E\rb_4-E\rb_2)), (w_5, y_5)\}$. Note that the arrivals
of $b_1$ and $b_5$ add no new white nodes to the queue, so their ``interruption'' of the currently serviced client is finished instantly. }
\end{figure}

In Fig.~\ref{fig: bi}, we illustrate the previous construction with a simple example. 
Let us point out the following features of the previous construction: 
\begin{itemize}
\item[-]
The construction stops in no more than $m+n$ steps. Indeed, at each step we complete the sampling of the neighbourhood for at least one vertex. If $\cA_{k-1}$ is empty, then we identify $\cN(V_k)$. If $\cA_{k-1}$ is non empty and the service of $U_k$ is interrupted by the arrival of $V_k$, we again identify $\cN(V_k)$. If $\cA_{k-1}$ is non empty but $U_k$ receives its service in full, then we complete the construction of the neighbourhood of $U_k$ at step $k$. We also note that in that case, $V_k$ is undefined. The stopping mechanism also ensures that all the $n$ black vertices will be explored by the end. 

\item[-]
The construction outputs a bipartite forest where each tree component is rooted at a black vertex. Indeed, the construction ensures that at each step all the neighbours of $V_k$ and $U_k$ are vertices that have not appeared so far in the construction. We further rank vertices of the same generation according to their arrival orders. 
The vertex set of this forest contains the $n$ black vertices and a subset of the $m$ white vertices. Add any remaining white vertices as isolated vertices to the forest and call the resulting graph $\cF$. 
We shall see that $\cF$ is a spanning forest of $B(n,m;\bx, \by)$. 
\end{itemize}

\begin{figure}
\centering
\scalebox{0.7}{
\begin{tikzpicture}
    \node[draw, circle, fill = black, minimum size = 0.5cm] (rootL) at (-5,0) {};
    \node[draw, circle, minimum size = 0.5cm] (11L) [above left = 1.2cm and 0.5cm of rootL.center] {};
    \node[draw, circle, minimum size = 0.5cm] (12L) [above right = 1.2cm and 0.5cm of rootL.center] {};

    \node[draw, circle, fill = black, minimum size = 0.5cm] (122L) [above right = 1.2cm and 0.5cm of 12L.center] {};
    \node[draw, circle, fill = black, minimum size = 0.5cm] (121L) [above left = 1.2cm and 0.5cm of 12L.center] {};
    \node[draw, circle, fill = black, minimum size = 0.5cm] (111L) [above left = 1.2cm and 0.5cm of 11L.center] {};

    \node[draw, circle, minimum size = 0.5cm] (1111L) [above left = 1.2cm and 0.5cm of 121L.center] {};
    \node[draw, circle, minimum size = 0.5cm] (1221L) [above right = 1.2cm and 0.5cm of 121L.center] {};

    \node[draw, circle, fill = black, minimum size = 0.5cm] (12211L) [above = 1.2cm of 1111L.center] {};
    \node[draw, circle, fill = black, minimum size = 0.5cm] (12212L) [above = 1.2cm  of 1221L.center] {};
    
    \draw[thick] (1111L)--(111L)--(11L)--(rootL)--(12L)--(121L);
    \draw[thick] (12L)--(122L)--(1221L)--(12211L);
    \draw[thick] (1221L)--(12212L);
    
    \node[draw, circle, fill = black, minimum size = 0.5cm] (rootM) at (0,0) {};
    \node[draw, circle, minimum size = 0.5cm, opacity = 0.5] (11M) [above left = 1.2cm and 0.5cm of rootM.center] {};
    \node[draw, circle, minimum size = 0.5cm, opacity = 0.5] (12M) [above right = 1.2cm and 0.5cm of rootM.center] {};

    \node[draw, circle, fill = black, minimum size = 0.5cm] (122M) [above right = 1.2cm and 0.5cm of 12M.center] {};
    \node[draw, circle, fill = black, minimum size = 0.5cm] (121M) [above left = 1.2cm and 0.5cm of 12M.center] {};
    \node[draw, circle, fill = black, minimum size = 0.5cm] (111M) [above left = 1.2cm and 0.5cm of 11M.center] {};

    \node[draw, circle, minimum size = 0.5cm, opacity = 0.5] (1111M) [above left = 1.2cm and 0.5cm of 121M.center] {};
    \node[draw, circle, minimum size = 0.5cm, opacity = 0.5] (1221M) [above right = 1.2cm and 0.5cm of 121M.center] {};

    \node[draw, circle, fill = black, minimum size = 0.5cm] (12211M) [above = 1.2cm of 1111M.center] {};
    \node[draw, circle, fill = black, minimum size = 0.5cm] (12212M) [above = 1.2cm  of 1221M.center] {};

    \draw[thick] (1111M.center)--(111M)--(11M.center)--(rootM)--([xshift = 0.125cm]12M.west)--(121M);
    \draw[thick] (rootM)--([xshift = -0.125cm]12M.east)--(122M)--([xshift = 0.125cm]1221M.west)--(12211M);
    \draw[thick] (122M)--([xshift = -0.125cm]1221M.east)--(12212M);


    \node[draw, circle, fill = black, minimum size = 0.5cm] (rootR) at (5,0) {};
    \node[circle, minimum size = 0.5cm] (11R) [above left = 1.2cm and 0.5cm of rootR.center] {};
    \node[circle, minimum size = 0.5cm] (12R) [above right = 1.2cm and 0.5cm of rootR.center] {};

    \node[draw, circle, fill = black, minimum size = 0.5cm] (122R) [above right = 1.2cm and 0.5cm of 12R.center] {};
    \node[draw, circle, fill = black, minimum size = 0.5cm] (121R) [above left = 1.2cm and 0.5cm of 12R.center] {};
    \node[draw, circle, fill = black, minimum size = 0.5cm] (111R) [above left = 1.2cm and 0.5cm of 11R.center] {};

    \node[circle, minimum size = 0.5cm] (1111R) [above left = 1.2cm and 0.5cm of 121R.center] {};
    \node[circle, minimum size = 0.5cm] (1221R) [above right = 1.2cm and 0.5cm of 121R.center] {};

    \node[draw, circle, fill = black, minimum size = 0.5cm] (12211R) [above = 1.2cm of 1111R.center] {};
    \node[draw, circle, fill = black, minimum size = 0.5cm] (12212R) [above = 1.2cm  of 1221R.center] {};

    \draw[thick] (111R)--(rootR)--(122R)--(12211R);
    \draw[thick] (rootR)--(121R);
    \draw[thick] (122R)--(12212R);

    \node (LtoM1) [right = 0.25cm of 122L] {};
    \node (LtoM2) [left = 0.25cm of 111M] {};
    \draw[->,very thick] (LtoM1)--(LtoM2);

    \node (MtoR1) [right = 0.25cm of 122M] {};
    \node (MtoR2) [left = 0.25cm of 111R] {};
    \draw[->,very thick] (MtoR1)--(MtoR2);
\end{tikzpicture}}
\caption{\label{fig: bf} Black forest inside the bipartite forest.}
\end{figure}
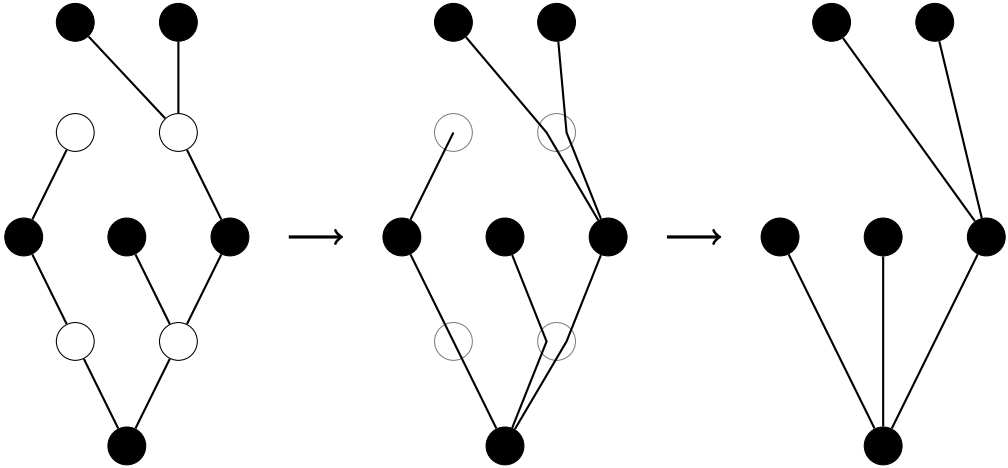

To make the connection between $\cF$ and the previous queuing system, denote by $\cF\rb$ the following forest induced by $\cF$ on the set of black vertices: a black vertex $b$ is the parent of the black vertex $b'$ in $\cF\rb$ if and only if $b$ is the grandparent of $b'$ in $\cF$; see Fig.~\ref{fig: bf}. 
Let us also note that each black vertex $b_i$ appears as a unique $V_k$ during the previous construction. As a result, we have assigned an interval $I(i):=I(V_k)$ used to sample its offspring. 
Meanwhile, consider a queuing system of $n$ clients with respective arrival times $(E\rb_i)_{i\in [n]}$ and the service time for the client arriving at $E\rb_i$ being
\begin{equation}
\label{def: Delta}
\Delta^{\mathrm{bi}}_i:=\sum_{j\in [m]}y_j\mathbf 1_{\{E\rw_j\in I(i)\}}, \quad 1\le i\le n. 
\end{equation}
Call this the {\it bipartite LIFO-queue} and let $\hat\cF$ be the associated forest of rooted ordered trees. 

\begin{lem}
\label{lem: LIFO}
We have $\hat\cF= \cF\rb$. 
\end{lem}

\begin{proof}
Given the queue $(\cA_k)_{k\ge 0}$ for the white vertices, we now define a corresponding queue for the black vertices by merging the white vertices that are offspring of a common black vertex. Indeed, we note that each time new entries are added to $(\cA_k)_{k\ge 1}$, they are added as a collection $\cN(V_k)$, where $V_k$ is some black vertex. Aggregate all the service requests from $\cN(V_k)$ into a single request from the client $V_k$. Operate the queue in an obvious way: if a client from $\cN(V_k)$ is being served at time $t$ in the previous queue, then say $V_k$ is served at time $t$ in the current queue. It is then straightforward to check that the queue for the black vertices coincides with the bipartite LIFO-queue. The conclusion follows. 
\end{proof}

\smallskip
\noindent
{\bf Surplus edges. }
In general, the random graph $B(n, m; \mathbf x, \mathbf y)$ is not a forest. To complete the construction, we need to include some additional edges without affecting the connected component structure of $\cF$. These additional edges will be referred to as {\it surplus edges}. Denote by $K$ the total steps that the construction takes. For $1\le k\le K$, if $V_k$ is defined in that step and $\cA_{k-1}\ne\varnothing$, then set 
\[
S_k = \big\{(V_k, w_j): w_j\in \cA_{k-1}\big\}; 
\]
otherwise, set $S_k=\varnothing$. 
As we will shortly see, $S_{k}$ contains all the candidates for a potential surplus edge at step $k$. 
If $w_j\in \cA_{k-1}$, let us denote $y_j(k)$ the remaining amount of service that $w_j$ has yet to receive when $V_k$ arrives. 
Given $(S_k)_{1\le k\le K}$, let $\{I_{i,j}: i\in [n], j\in [m]\}$
be a collection of independent Bernoulli variables that satisfies 
\begin{equation}
\label{def: sur-prob}
\mathbb E\big[I_{i, j}\,\big|\,(S_k)_{1\le k\le K}\big] = \mathbf 1_{\{\exists\, k\le K: (b_i, w_j)\in S_k\}}\big\{1-\exp(-y_j(k)x_i/z_{n,m})\big\}. 
\end{equation}
Then 
\begin{equation}
\label{def: sur}
\cE_{n,m}:=\big\{(b_i, w_j): I_{i,j}=1,  i\in [n], j\in [m]\big\}
\end{equation}
denotes the set of endpoints of sampled surplus edges. 
Let $B'$ be the {\it simple} graph obtained after adding all the edges from $\cE_{n,m}$ to $\cF$, ignoring multi-edges, and forgetting the roots and vertex ordering in $\cF$. 

\begin{lem}
\label{lem: B-graph}
$B'$ is distributed as $B(n, m; \mathbf x, \mathbf y)$. 
\end{lem}

\begin{proof}
See page~\pageref{pf: B-graph}.
\end{proof}

Thanks to Lemma~\ref{lem: B-graph}, we  will tacitly {\it assume from now on} that $B(n, m; \mathbf x, \mathbf y)$ has been constructed using the exponential variables $(E\rb_i)_{i\in [n]}, (E\rw_j)_{j\in [m]}$ and $\cE_{n,m}$.

\subsection{Graph encoding processes}
\label{sec: xy-graph-coding}

Recall from Section~\ref{sec: explore} the generic LIFO-queue with arrivals at $(t_i)_{1\le i\le N}$ and service requests $(\Delta_i)_{1\le i\le N}$. The {\it service load} at time $t$ refers to the total amount of unfulfilled service time from the clients in the queue at that moment. To track that quantity, let us introduce:
\begin{equation}
\label{def: ZQ}
Z^{Q}_t = -t+\sum_{1\le i\le N} \Delta_i\mathbf 1_{\{t_i\le t\}}, \quad t\ge 0. 
\end{equation}
We note that the negative drift represents the discharge rate of the server, while $Z^Q_t$ jumps upwards of size $\Delta_i$ at each $t_i$, $1\le i\le N$. From this, it is not difficult to convince oneself that the service load at time $t$ is given precisely by 
\[
Z^Q_t - \inf_{u\le t}Z^Q_u.
\]
Recall again from Section~\ref{sec: explore} that we have associated a forest of rooted ordered trees to the queue. Let us observe that in the aforementioned forest the ancestors of a client are those waiting in the queue while the client is being served. Hence, the height of a client in the forest corresponds to the queue length for the duration of its service.\footnote{We define the height of a root vertex as 1. The queue length is the number of clients in the queue, including the one being currently served.} 
To find that quantity, we follow Le Gall and Le Jan~\cite{LGLJ98} to introduce the so-called {\it height process} functional of $Z^Q$: for each $t\ge 0$, let
\begin{equation}
\label{def: ht-queue}
H^Q_t = \#\Big\{s\le t: Z^Q_{s-}< \inf_{u\in [s, t]}Z^Q_u\Big\}=\#\Big\{t_i\le t: 1\le i\le N, Z^Q_{t_i-}< \inf_{u\in [t_i, t]}Z^Q_u\Big\}.
\end{equation}
We claim that $H^Q_{t_i}$ is the height of the client arriving at time $t_i$, $1\le i\le N$. Indeed, due to the Last-In-First-Out rule, the client arriving at time $t_i$ departs from the queue at the first moment after $t_i$ when the service load falls back to the level just prior to its arrival, namely  the moment
\[
\inf\Big\{t>t_i: Z^Q_t \le Z^Q_{t_i-}\Big\} = \inf\Big\{t>t_i: \inf_{u\in[t_i,  t]}Z^Q_u\le Z^Q_{t_i-}\Big\}. 
\]
Therefore, the client arriving at time $t_i$ is still in the queue at time $t$ if and only if $Z^Q_{t_i-}<\inf_{u\in [t_i, t]}Z^Q_u$. 
It follows that $H^Q_{t}$ counts the number of clients in the queue at time $t$. In fact, the previous arguments can be used to establish a stronger identity: denote by $\dgr^{Q}$ the graph distance  of the forest and let $v_t$ be the client being served at time $t$ (set $v_t=\dagger$ if the queue is empty at time $t$); then for any $0\le s\le t$, we have
\begin{equation}
\label{id: tree-dist}
\dgr^{Q}(v_s, v_t) = H^Q_s+H^Q_t-2\min_{s\le u\le t} H^Q_u, 
\end{equation}
if $\inf_{u\le s}Z^Q_u=\inf_{u\le t}Z^Q_u$, indicating the two clients are in the same tree component, and $\dgr^Q(v_s, v_t)=\infty$ otherwise. Underlying the identity~\eqref{id: tree-dist} is the fact that $H^Q$ can be viewed as a time-change of the so-called contour process of the forest; we refer to~\cite{BrDuWa21}, Section 3.2 for further detail. 

\medskip
To introduce the relevant processes for the bipartite graph $B(n,m; \bx, \by)$, let us recall the bipartite LIFO-queue with arrivals at $(E\rb_i)_{i\in [n]}$ and service requests $(\Delta\bi_i)_{i\in [n]}$ from~\eqref{def: Delta}. 
We define for $t>0$, 
\begin{equation}
\label{def: Hdis}
Z^{\bx,\by}_t  = -t + \sum_{i\in [n]}\Delta^{\mathrm{bi}}_i\mathbf 1_{\{E\rb_i\le t\}} \ \text{ and }\  H^{\bx,\by}_t = \#\Big\{s\le t: Z^{\bx,\by}_{s-}< \inf_{u\in [s, t]}Z^{\bx,\by}_u\Big\}.
\end{equation}
Then $H^{\bx, \by}_t$ is the length of the bipartite LIFO-queue at time $t$. 
We recall the forest $\cF$ produced during the LIFO-queue construction of $B(n,m;\bx, \by)$, which is a spanning forest of the graph according to Lemma~\ref{lem: B-graph}. Recalling from Lemma~\ref{lem: LIFO} the connection between $\cF$, the black sub-forest $\cF\rb$ and the forest $\hat\cF$ associated to the bipartite queue, we have
\begin{equation}
\label{id: ht}
2H^{\bx, \by}_{E\rb_i}-1+2\cdot\mathbf 1_{\{\Delta\bi_i=0\}}=\text{height of the black vertex $b_i$ in $\cF$}, \quad 1\le i\le n. 
\end{equation}
In above, the appearance of the indicator function is due to the fact that if $\Delta\bi_i=0$, $E\rb_i$ is not counted in $H^{\bx, \by}_{E\rb_i}$. Note that in that case, $b_i$ has no children and is necessarily a leaf of $\cF$. 
The process $Z^{\bx, \by}$ can be expressed as a more explicit function of $\bx, \by, (E\rb_i)_{ i\in [n]}, (E\rw_j)_{j\in [m]}$. Indeed, let us define
\begin{equation}
\label{def: lambda}
\Lambda^{\bx}(t)=\sum_{i\in [n]}x_i\mathbf 1_{\{E\rb_i\le t\}} \quad \text{and}\quad \Lambda^{\by}(t)=\sum_{j\in [m]}y_j\mathbf 1_{\{E\rw_j\le t\}}. 
\end{equation}

\begin{lem}
\label{lem: Z}
We have
\begin{equation}
\label{id: Z}
Z^{\bx, \by}_t = -t + \Lambda^{\by}\circ\Lambda^{\bx}(t), \quad t\ge 0.
\end{equation}
\end{lem}

\begin{proof}
See page~\pageref{pf: Z}.
\end{proof}

Let $(g_{n, i}, d_{n, i}), 1\le i\le \kappa'_{n,m}$ be the maximal open intervals of the following set
\[
\Big\{t>0: Z^{\bx, \by}_t>\inf_{s\le t}Z^{\bx, \by}_s\Big\},
\]
ranked in increasing order of their left endpoints. 
Meanwhile, we say a graph is {\it non trivial} if it contains at least one edge; 
let $T_{n, i}$ be the $i$-th non trivial tree that appears in the LIFO-queue construction of $B(n,m;\bx, \by)$. We denote by $\bx(T_{n, \ell})$ (resp.~$\by(T_{n, i})$) the total $\mathbf x$-weights of the black vertices in $T_{n, i}$ (resp.~the total $\by$-weights of the white vertices in $T_{n, i}$). Let $\kappa_{n,m}$ stand for the number of non trivial connected components in $B(n, m; \bx, \by)$. 

\begin{lem}
\label{lem: exc}
We have $\kappa_{n,m}=\kappa'_{n,m}$. For $i=1, 2, \dots, \kappa_{n,m}$, we have almost surely
\begin{equation}
\label{id: mass}
d_{n, i}-g_{n, i}= \mathbf y(T_{n, i}) \quad \text{and}\quad \Lambda^{\bx}(d_{n, i})-\Lambda^{\bx}(g_{n, i}-) = \bx(T_{n, i}). 
\end{equation}
\end{lem}

\begin{proof}
See page~\pageref{pf: exc}.
\end{proof}

Let us  note that the previous lemma suggest that we can recover the connected components of the graph by tracking the excursion intervals of $(Z^{\bx, \by}_t)_{t\ge 0}$ above its running infimum. 

\subsection{A coupling with Poisson point processes in the i.i.d.~case}
\label{sec: coupling}

From now on, we are only concerned with the i.i.d.~case. Recall that we take $z_{n,m}=\sqrt{mn}$ in this case. Inspired by the work~\cite{Jo14}, we introduce a pair of Poisson point processes which will generate both the vertex weights and the exponential variables used in the previous construction. 
Recall $F\rb$ and $F\rw$, the c.d.f.~for the respective weight distributions for black and white vertices. Let $\Pi\rb=\{(E\rb_{i}, X_{i}): 1\le i\le N\rb_{n}\}$ be a Poisson point measure on $\R_{+}^{2}$ with intensity measure
\[
\pi\rb(dt, dx) =\sqrt{\tfrac{n}{m}}\, xe^{-xt/\sqrt{mn}}dF\rb(x) dt.
\]
Note in particular that $\int_{\R^2_+}\pi\rb(dt, dx) = n$. 
Similarly, let $\Pi\rw=\{(E\rw_{j}, Y_{j}): 1\le j\le N\rw_{m}\}$ be a Poisson point measure on $\R_{+}^{2}$ with intensity measure
\[
\pi\rw(ds, dy) = \sqrt{\tfrac{m}{n}}\,y e^{-ys/\sqrt{mn}}dF\rw(y) ds,
\]
and independent of $\Pi\rb$. We observe the following distributional properties of $\Pi\rb$ and $\Pi\rw$, whose proof is standard and therefore omitted. 

\begin{lem}
\label{lem: coupling}
The total number of atoms $N\rb_{n}$ of $\Pi\rb$ is a Poisson variable of mean $n$; similarly, $N\rw_{m}$ is a Poisson variable of mean $m$. Given $N\rb_n$, the pairs $(E\rb_i, X_i), 1\le i\le N\rb_n $, are i.i.d.~with the joint distribution $\pi\rb(dt, dx)/n$. In particular, given $N\rb_n=n$, the marginal distribution of $(X_1, X_2, \cdots, X_n)$ is that of an i.i.d.~sequence with common distribution $dF\rb$, and conditional on $(X_1, X_2, \dots, X_n)$, $\{E\rb_i: 1\le i\le n\}$ are independent and each with Exp$(X_i/\sqrt{mn})$ distribution. Similarly, given $N\rw_m=m$, the marginal distribution of $(Y_1, Y_2, \cdots, Y_m)$ is that of an i.i.d.~sequence with common distribution $dF\rw$, and conditional on $(Y_1, Y_2, \dots, Y_m)$, $\{E\rw_j: 1\le j\le m\}$ are independent and each with Exp$(Y_j/\sqrt{mn})$ distribution. 
\end{lem}

Let us consider the bipartite graph with random weight sequences $\tilde{\mathbf X}_n:=\{x_i: 1\le i\le N\rb_n\}$ and $\tilde{\mathbf Y}_m=\{y_j: 1\le j\le  N\rb_m\}$. In particular, on the event $\{N\rb_n= n, N\rw_m= m\}$, the previous lemma ensures that we obtain a version of $B_{n,m}=B(n, m; \mathbf X_n, \mathbf Y_m)$. On the complement of that event, we obtain a graph with a vertex set of a different size; nevertheless, it still makes sense to talk about the graph $B(n,m; \tilde{\mathbf X}_n, \tilde{\mathbf Y}_m)$. 
For this reason, let 
\[
\bP_{n,m} \text{ be the joint distribution of  $\Pi\rb$ and $\Pi\rw$}. 
\] 
The construction in Section~\ref{sec: explore} and Lemma~\ref{lem: coupling} combined allow us to view $B(n,m; \tilde{\mathbf X}_n, \tilde{\mathbf Y}_m)$ as a function of $\tilde{\mathbf X}_n$, $\tilde{\mathbf Y}_m$, $(E\rb_i)_{1\le i\le N\rb_n}$, $(E\rw_j)_{1\le j\le N\rw_m}$ and some extra randomness used to sample the surplus edges that can be assumed to be independent of the previous random variables. 
By enlarging the probability space to accommodate the extra randomness, we will {\it assume from now on} that $B(n,m; \tilde{\mathbf X}_n, \tilde{\mathbf Y}_m)$ is thus defined under  $\bP_{n,m}$. Furthermore, the graph encoding processes introduced in Section~\ref{sec: xy-graph-coding} also have a version under $\bP_{n,m}$: let
\begin{align}
\label{def: Lam}
&\Lambda\rbn(t) = \sum_{1\le i\le N\rb_n}X_i\mathbf 1_{\{E\rb_i\le t\}}, \quad  \Lambda\rwm(t) = \sum_{1\le j\le N\rw_m}Y_j\mathbf 1_{\{E\rw_j\le t\}},\\ \label{def: ZH}
&Z^{n,m}_t = -t +\Lambda\rwm\circ\Lambda\rbn(t), \quad H^{n,m}_t = \#\Big\{s\le t: Z^{n,m}_{s-}< \inf_{u\in [s, t]}Z^{n,m}_u\Big\}.
\end{align}

\subsection{Limit theorems for the graph encoding processes} 
\label{sec: limit-proc}

Recall that we have assumed $m=\lf \theta n\rf$; see~\eqref{hyp: mod-clst}. Recall the c.d.f.~$F\rb, F\rw$ for the weight distributions and the notation $\sigma\rb_r, \sigma\rw_r$ for their respective $r$-th moments.  The critical regime that we are interested in requires  $\sigma\rb_2\cdot \sigma\rw_2=1$; see~\eqref{hyp: crit}. 
We have distinguished three scenarios regarding the tail behaviours of $F\rb$ and $F\rw$: {\it double finite third moments, one dominant heavy tail} and {\it matched heavy tails}. 
We describe here the scaling limit of $(Z^{n,m}, H^{n,m})$, introduced in~\eqref{def: ZH}, in each of these scenarios. For $I\subseteq \R_+$ and $d\in\N$, we denote by $\mathbb D(I, \R^d)$ the space of c\`adl\`ag functions from $I$ to $\R^d$ equipped with the Skorokhod $J_1$-topology. 

\smallskip
\noindent
{\bf Scaling limit of $Z^{n,m}$.} 
For $i\in \{1, 2, 3\}$, let $\cL^{(i)}=(\cL^{(i)}_t)_{t\ge 0}$ be a spectrally positive stable L\'evy process  whose law is characterised by the Laplace transform: 
\begin{equation}
\label{def: L}
\mathbb E[e^{-\lambda \cL^{(i)}_t}] = \exp\big(t \Psi_i(\lambda)\big), \quad \lambda\ge 0,\; t\ge 0,
\end{equation}
where 
\begin{align*}
&\Psi_1(\lambda) =C(2) \,\mathrm C_1  \lambda^2, \quad \mathrm C_1= \sigma\rw_3\sigma\rb_2+\sqrt\theta(\sigma\rw_2)^2\sigma\rb_3,\\ 
&\Psi_2(\lambda) = C(\alpha) \mathrm C_2 \lambda^{\alpha}, \quad \mathrm C_2= C\rb \theta^{\frac{\alpha-1}{2}}(\sigma\rw_2)^{\alpha}, \\ 
&\Psi_3(\lambda) = C(\alpha)\mathrm C_3 \lambda^{\alpha}, \quad \mathrm C_3= C\rw\sigma\rb_2+C\rb\theta^{\frac{\alpha-1}{2}}(\sigma\rw_2)^{\alpha},
\end{align*}
and $C(\alpha) = (\alpha+1)\Gamma(2-\alpha)/(\alpha(\alpha-1))$ for $\alpha\in (1, 2)$ and $C(2) = 1/2$. 
We shall refer to $\Psi_i$ as the {\it Laplace exponent} of $\cL^{(i)}$. 
Now  let $\cZ^{(i)}=(\cZ^{(i)}_t)_{t\ge 0}$ be a stochastic process with c\`adl\`ag sample paths whose distribution satisfies the following absolute continuity relationship: for each $t\ge 0$ and every measurable function $F:  \mathbb D([0, t], \R)\to \R_+$, we have
\begin{equation}
\label{id: girsanov}
\mathbb E\Big[F\Big(\big(\cZ^{(i)}_s\big)_{s\le t}\Big)\Big] = \mathbb E\Big[F\Big(\big(\cL^{(i)}_s\big)_{s\le t}\Big)\cdot \cE^{(i)}_t\Big]
\end{equation}
with 
\begin{equation}
\label{def: cE}
\cE^{(i)}_t = \exp\Big\{-\int_0^t \frac{\sigma\rb_2 s}{\theta}\,d\cL^{(i)}_s -\int_0^t \Psi_i\Big(\frac{\sigma\rb_2s}{\theta}\Big)ds\Big\}.
\end{equation}
Underlying the identity~\eqref{id: girsanov} is the property that $(\cE^{(i)}_t)_{t\ge 0}$ is a unit-mean positive martingale; we refer to Appendix~\ref{sec: Z-prop} for further details. 
We will shortly see that $\cZ^{(i)}$, $i=1, 2, 3$, all appear in the scaling limit of $Z^{n,m}$ under the various assumptions of $F\rb$ and $F\rw$. 

Let us  present some alternative definitions of $\cZ^{(i)}$. For $i=1$, we note that $\cL^{(1)}$ is a Brownian motion. In that case, the identity in~\eqref{id: girsanov} is a special case of Girsanov's Theorem for Brownian motions (see for instance Theorem VIII.1.12 in~\cite{Yorbook}), and we have in fact
\begin{equation}
\label{id: Z1}
\cZ^{(1)} \eqd \Big\{\cL^{(1)}_t - \Big(\frac{\sigma\rb_2}{\theta}\Big)^{-1}\cdot\Psi_1\Big(\frac{\sigma\rb_2 t}{\theta}\Big): t\ge 0\Big\}.
\end{equation}
In words, $\cZ^{(1)}$ is a Brownian motion with a parabolic drift. 
For both $i=2$ and $i=3$, $\cL^{(i)}$ is an $\alpha$-stable process. Let $(\cM_t)_{t\ge 0}$ be a martingale process with c\`adl\`ag sample paths which has independent increments and whose marginal laws are characterised as follows: for all $\lambda\ge 0$ and $t\ge 0$, 
\begin{equation}
\label{def: Malpha}
\mathbb E\big[e^{-\lambda \cM_t}\big] = \exp\Big\{(\alpha+1)\int_0^t\int_0^{\infty}(e^{-\lambda x}-1+\lambda x)e^{-\sigma\rb_2xs/\theta}x^{-\alpha-1}dx\,ds\Big\}. 
\end{equation}
Then a Girsanov-type theorem for L\'evy processes (see Appendix~\ref{sec: Z-prop}) says that 
\begin{equation}
\label{id: cZ-drift}
\cZ^{(i)}\eqd \Big\{\mathrm C_i\cdot \cM_t - \Big(\frac{\sigma\rb_2}{\theta}\Big)^{-1}\cdot\Psi_i\Big(\frac{\sigma\rb_2 t}{\theta}\Big): t\ge 0\Big\}, \quad i=2, 3. 
\end{equation}

\smallskip
\noindent
{\bf Height processes in the limit.} 
For the spectrally positive L\'evy process $\cL^{(i)}, i\in \{1, 2, 3\}$, the notion of height process has been introduced by Le Gall and Le Jan~\cite{LGLJ98}, who show that there exists a stochastic process $(\widehat\cH^{(i)}_t)_{t\ge 0}$ of continuous sample paths so that the following limit exists in probability for all $t\ge 0$:
\begin{equation}
\label{def: Hlevy}
\widehat\cH^{(i)}_t:= \lim_{\epsilon\to 0}\frac{1}{\epsilon} \int_0^t \mathbf 1_{\{ \cL^{(i)}_s\le \inf_{u\in [s, t]}\cL^{(i)}_u+\epsilon\}}ds. 
\end{equation}
The identity~\eqref{id: girsanov} allows us to assert the existence of an analogous height process for $\cZ^{(i)}$. Namely, for $i\in \{1, 2, 3\}$, there exists some process $(\cH^{(i)}_t)_{t\ge 0}$ of continuous sample paths satisfying
\begin{equation}
\label{def: cH}
\cH^{(i)}_t:= \lim_{\epsilon\to 0}\frac{1}{\epsilon} \int_0^t \mathbf 1_{\{ \cZ^{(i)}_s\le \inf_{u\in [s, t]}\cZ^{(i)}_u+\epsilon\}}ds, 
\end{equation}
where the limit above exists in probability for all $t\ge 0$. Furthermore, for $\cZ^{(1)}$, which is a Brownian motion with a parabolic drift, its height process takes a simple form, as almost surely we have 
\begin{equation}
\label{id: cH1}
\cH^{(1)}_t = \frac{2}{\mathrm C_1} \Big(\cZ^{(1)}_t-\inf_{s\le t}\cZ^{(1)}_s\Big), \quad t\ge 0. 
\end{equation}
We refer to~\cite{DuLG02} for further background on the height processes of L\'evy processes and particularly Eq.~(1.7) there for the simplification in the Brownian case. 

\begin{prop}[Convergence of the graph encoding processes] 
\label{prop: cv-ZH}
Assume~\eqref{hyp: mod-clst} and~\eqref{hyp: crit}. 
For $i\in \{1, 2, 3\}$, let $(\cZ^{(i)}_t)_{t\ge 0}$ be the stochastic process defined in~\eqref{id: girsanov} and $(\cH^{(i)}_t)_{t\ge 0}$ be as in~\eqref{def: cH}. 
\begin{enumerate}
\item
{\sf (Double finite third moments)}
Under both~\eqref{hyp: b-third} and~\eqref{hyp: w-third}, the following weak  convergence takes place under $\bP_{n,m}(\cdot\,|\,N\rw_m=m, N\rb_n=n)$: 
\begin{equation}
\label{cv: ZH1}
    \Big\{n^{-\frac13}Z^{n,m}_{n^{2/3}t}, n^{-\frac13}H^{n,m}_{n^{2/3}t}: t\ge 0\Big\} \Longrightarrow \Big\{\cZ^{(1)}_t, \cH^{(1)}_t: t\ge 0\Big\} \text{ in } \mathbb D(\R_+, \R^2). 
\end{equation}

\item
{\sf (One dominant heavy tail)}
Assume that~\eqref{hyp: b-power} holds with $\alpha\in (1, 2)$ and $C\rb\in (0, \infty)$. Assume that either~\eqref{hyp: w-third} is true or~\eqref{hyp: w-power} holds with $\alpha'\in (\alpha, 2)$ and $C\rw\in (0, \infty)$. Then the following weak  convergence takes place under $\bP_{n,m}(\cdot\,|\,N\rw_m=m, N\rb_n=n)$: 
{\small
\begin{equation}
\label{cv: ZH2}
    \Big\{n^{-\frac{1}{\alpha+1}}Z^{n,m}_{ n^{\alpha/(\alpha+1)}t}, n^{-\frac{\alpha-1}{\alpha+1}}H^{n,m}_{n^{\alpha/(\alpha+1)}t}: t\ge 0\Big\}\Longrightarrow \Big\{\cZ^{(2)}_{t}, \cH^{(2)}_t: t\ge 0\Big\} \text{ in } \mathbb D(\R_+, \R^2). 
\end{equation}}

\item
{\sf (Matched heavy tails)}
Assume that both~\eqref{hyp: b-power} and~\eqref{hyp: w-power} hold with $\alpha=\alpha'\in (1, 2)$ and $C\rb\in (0, \infty), C\rw\in (0, \infty)$. Then the following weak  convergence takes place under $\bP_{n,m}(\cdot\,|\,N\rw_m=m, N\rb_n=n)$: 
{\small
\begin{equation}
\label{cv: ZH3}
    \Big\{n^{-\frac{1}{\alpha+1}}Z^{n,m}_{ n^{\alpha/(\alpha+1)}t}, n^{-\frac{\alpha-1}{\alpha+1}}H^{n,m}_{ n^{\alpha/(\alpha+1)}t}: t\ge 0\Big\}\Longrightarrow \Big\{\cZ^{(3)}_{t}, \cH^{(3)}_t: t\ge 0\Big\} \text{ in } \mathbb D(\R_+, \R^2). 
\end{equation}}
\end{enumerate}
\end{prop}

The proof of Proposition~\ref{prop: cv-ZH}, given in Section~\ref{sec: pf-cv-ZH}, is the main focus of the proof section~\ref{sec: proof}.

\subsection{Convergence of the surplus edges}
\label{sec: surplus}

Roughly speaking, Proposition~\ref{prop: cv-ZH} from the previous section implies the convergence of the spanning forest $\cF$ to the spanning forest of the continuum graphs. 
Recall that the finite graph can be obtained from its spanning forest by inserting surplus edges. 
We show here that the collection of surplus edges also has a trackable limit distribution. 
Together with the previous convergence of graph encoding processes,  this will provide the principal arguments of our main results. 

Before presenting our result on the convergence of the surplus edges, we first introduce some modifications on how we sample these edges in $B_{n,m}$. The first modification is deterministic and consists in altering locally the graph distance. The second modification highlights the hidden Poissonian features of the surplus edges. 

\smallskip
\noindent
{\bf Local modifications in the graph distance. } Recall $\cE_{n,m}$, the set of surplus edges of $B_{n,m}$ under $\bP_{n,m}$, and recall that all its elements are of the form $e=(b_i, w_j)$ for which there exists some unique $1\le k\le K$ so that $b_i=V_k$ and $w_j\in \cA_{k-1}$. In words, surplus edges are only formed between the black vertex explored at step $k$ and a white vertex from the queue at that time.  Let $e=(b_i, w_j)\in \cE_{n,m}$ and let $k(e)$ be such that $V_{k(e)}=b_i$. Since the queue $\cA_{k-1}$ is filled by the white neighbours of $(V_{\ell})_{\ell<k}$ that have not been fully explored, we can find a unique $k'(e)<k$ so that $w_j\in \cN(V_{k'(e)})$.  Now let $(V_{k(e)}, V_{k'(e)})$ be the edge obtained from $e$ by replacing the white vertex with its parent in $\cF$ and define the collection
\begin{equation}
\label{def: cE'}
\cE'_{n,m}= \big\{\big(V_{k(e)}, V_{k'(e)}\big): e\in \cE_{n,m}\big\}. 
\end{equation}
Let $B''$ be the graph obtained after adding all the edges from $\cE'_{n,m}$ to $\cF$ and forgetting the roots and vertex ordering in $\cF$. Although $B''$ is no longer a bipartite graph, it provides a convenient approximation of $B_{n,m}$ as shown in the following lemma. Recall the ranked sequence $\{C^{\bx}_{n, (k)}: 1\le k\le \kappa_{n,m}\}$ of the non trivial connected components of $B_{n,m}$. Let $\mathrm{Id}$ be the identity map from the vertex set of $B_{n,m}$ to that of $B''$. Recall that $\dgr$ stands for the graph distance of $B_{n,m}$; let $d''$ be the graph distance of $B''$. 

\begin{lem}
\label{lem: dis}
For $n\in \N$ and $1\le k\le \kappa_{n,m}$,  $\mathrm{Id}(C^{\bx}_{n, (k)})$ is a connected component of $B''$; moreover, denoting for each $a>0$, 
\[
\mathrm{dis}_k(\mathrm{Id}, a)=\sup\{|a\cdot \dgr(x, y) - a\cdot d''(x, y)|: x, y\in C^{\bx}_{n, (k)}\},
\]
we have 
\[
 \mathrm{dis}_k(\mathrm{Id}, a) \le a\cdot s_{n, k}, 
 \]
 where $s_{n,k}$ is the number of surplus edges with both endpoints in $C^{\bx}_{n,(k)}$. 
\end{lem}

Let us observe that in $B_{n,m}$, $V_{k(e)}$ is at distance 2 from $V_{k'(e)}$, as they are both adjacent to a common white vertex. In $B''$ however, the two vertices are neighbours with each other, and this is the main difference between $B_{n,m}$ and $B''$. 
Roughly speaking, the upper bound in Lemma~\ref{lem: dis} is a reflection of the fact that a geodesic in $C^{\bx}_{n, (k)}$ contains at most $s_{n, k}$ surplus edges; we refer to page~\pageref{pf: dis} for a proof of the lemma.

\smallskip
\noindent
{\bf Poissonisation of surplus edges. } We present here an alternative way of sampling surplus edges via a Poisson point process. To that end, we will need an auxiliary function $\Sigma^{n,m}$ introduced below. Roughly speaking, the function transfers $\by$-weights to $\bx$-weights.   Recall from~\eqref{def: Delta} the bipartite LIFO-queue with service requests $(\Delta\bi_k)$. Recall that if $\Delta\bi_k=0$, then Client $k$ will immediately leave the queue the moment it arrives.  Given the bipartite LIFO-queue, we set for $t>0$, 
\[
\mathfrak{D}_t = \big\{1\le k\le N\rb_n: \text{Client $k$ has departed by time $t$ if $\Delta\bi_k>0$ or $E\rb_k\le t$ if $\Delta\bi_k=0$}\big\}, 
\]
and
\begin{equation}
\label{def: Wt}
\mathfrak{W}_t = \big\{1\le k\le N\rb_n: \Delta\bi_k>0 \text{ and Client $k$ is in the queue at time $t$}\big\}.
\end{equation}
If $k\in \mathfrak{W}_t$, denote by $\Delta_k(t)$ the service time that Client $k$ has received by $t$. 
Let $\Sigma^{n,m}: \R_+\to\R_+$ be as follows: 
\begin{equation}
\label{def: Sigma}
\Sigma^{n,m}(t) = \sum_{k\ge 1}X_k\mathbf 1_{\{k\in \mathfrak{D}_t\}}+\sum_{k\ge 1
} \frac{\Delta_k(t)}{\Delta\bi_k} \cdot X_{k}\mathbf 1_{\{k \in \mathfrak{W}_t\}}. 
\end{equation}
For $1\le k\le N\rb_n$, denote by $J_k$ the set of time when Client $k$ is served, with the convention that $J_k=\{E\rb_k\}$ if $\Delta\bi_k=0$. 
Note that $J_k, 1\le k\le N\rb_n$ are disjoint sets. 
For $A\subseteq \R_+$, its image set by $\Sigma^{n,m}$ is defined as follows:
\[
\Sigma^{n,m}(A) =\{y: \exists\, t\in A \text{ s.t. } \Sigma^{n,m}(t-)\le y\le \Sigma^{n,m}(t)\}.
\]
We also write $|A|$ for the Lebesgues measure of a Borel set $A\subseteq \R_+$. 
The following result summarises some useful properties of $\Sigma^{n,m}$. 

\begin{lem}
\label{lem: sigma-p}
The function $t\mapsto \Sigma^{n,m}(t)$ is increasing and strictly increasing on the set of time when the queue is nonempty. Moreover, we have $|\Sigma^{n,m}(J_k)|=X_k$, $1\le k\le N\rb_n$. 
\end{lem}

\begin{proof}
See page~\pageref{pf: sigma-p}.
\end{proof}

\noindent
Let $Q_{n,m}=\{(s_{n,l}, y_{n,l}): 1\le l\le q_{n,m}\}$ be a Poisson point measure with rate $1/\sqrt{mn}$ on 
\[
D_{n,m}:=\big\{(s, y): \exists\,t\ge 0 \text{ s.t. } \Sigma^{n,m}(t-)\le s\le \Sigma^{n,m}(t), 0\le y\le Z^{n,m}_t-\inf_{u\le t}Z^{n,m}_u\big\}. 
\]
Suppose $(s_{n,l}, y_{n,l})\in Q_{n,m}$. By definition, we can find some $t_{n,l}\ge 0$ satisfying $\Sigma^{n,m}(t_{n,l}-)\le s_{n,l}\le \Sigma^{n,m}(t_{n,l})$. Since the queue is not empty at time $t$ if and only if $Z^{n,m}_t>\inf_{u\le t}Z^{n,m}_u$, with probability $1$, the queue is not empty at time $t_{n,l}$. As a result, $\Sigma^{n,m}$ is strictly increasing on a neighbourhood of $t_{n,l}$ and such $t_{n,l}$ is therefore unique a.s. 
Let $b_{n,l}$ be the client being served at time $t_{n,l}$. We set $t'_{n,l}=\sup\{u\le t_{n,l}: Z^{n,m}_u-\inf_{v\le u}Z^{n,m}_v\le y_{n,l}\}$, and  let $b'_{n,l}$ be the client being served at time $t'_{n,l}$. Define
\begin{equation}
\label{def: cE''}
\cS_{n,m}=\big\{(t_{n,l}, t'_{n,l}): 1\le l\le q_{n,m}\big\} \quad\text{and}\quad \cE''_{n,m}=\big\{(b_{n,l}, b'_{n,l}): 1\le l\le q_{n,m}, b'_{n,l}\ne b_{n, l}\big\}. 
\end{equation}

\begin{lem}
\label{lem: dist-sur}
$\cE''_{n,m}$ has the same distribution as $\cE'_{n,m}$. 
\end{lem}

\begin{proof}
See page~\pageref{pf: dist-sur}.
\end{proof}

Thanks to Lemmas~\ref{lem: dis}-\ref{lem: dist-sur},  we can focus on the collection $\cS_{n,m}$ instead of the original surplus edges. We now present the main result concerning its scaling limit. To that end, let us introduce that for $a>0$, 
\[
\mathfrak S[a](\cS_{n,m})=\big\{(t_{n,l}/a, t'_{n,l}/a): 1\le l\le q_{n,m}\big\}. 
\]
Recall $\cZ^{(i)}$ from~\eqref{id: girsanov}. For $i\in \{1, 2,3\}$, given $\cZ^{(i)}$, let $\cQ^{(i)}=\{(s_l, y_l): l\ge 1\}$ be a Poisson point measure with rate $1/\sqrt\theta$ on 
\begin{equation}
\label{def: cQ}
\cD^{(i)}:=\bigg\{(s, y): s\ge 0, 0\le y\le \cZ^{(i)}_{\sqrt\theta s/\sigma\rb_2}-\inf_{u\le  \sqrt\theta s/\sigma\rb_2}\cZ^{(i)}_u\bigg\}. 
\end{equation}
For $i\in \{1, 2, 3\}$ and $l\ge 1$, define $t_l=\sqrt{\theta} s_l/\sigma\rb_2$, 
\begin{equation}
\label{def: cS}
t'_l=\sup\{u\le t_l: \cZ^{(i)}_u-\inf_{v\le u}\cZ^{(i)}_v\le y_l\}, \text{ and let }  \cS^{(i)}=\big\{(t_l, t'_l): l\ge 1\big\}.  
\end{equation}
For $n\in\N\cup\{\infty\}$, let $\pi_n$ be a random locally finite measure of $\R^d$. We say that {\it $\pi_n$ converges in distribution to $\pi_{\infty}$} if $\int f d\pi_n$ converges in distribution to $\int fd\pi_{\infty}$ as $n\to\infty$ for all the continuous function $f: \R^d\to \R_+$ with compact support. 

\begin{prop}[Convergence of the surplus edges] 
\label{prop: cv-surplus}
We assume~\eqref{hyp: mod-clst} and~\eqref{hyp: crit}. 
\begin{enumerate}
\item
{\sf (Double finite third moments)} Under both~\eqref{hyp: b-third} and~\eqref{hyp: w-third}, the rescaled point measure $\mathfrak S[n^{\frac23}](\cS_{n,m})$
under $\bP_{n,m}(\cdot\,|\,N\rb_n=n, N\rw_m=m)$ converges in distribution  to $\cS^{(1)}$, jointly with the convergence in~\eqref{cv: ZH1}. 
\item
{\sf (One dominant heavy tail)} Assume that~\eqref{hyp: b-power} holds with $\alpha\in (1, 2)$ and $C\rb\in (0, \infty)$. Assume that either~\eqref{hyp: w-third} is true or~\eqref{hyp: w-power} holds with $\alpha'\in (\alpha, 2)$ and $C\rw\in (0, \infty)$. Then $\mathfrak S[n^{\frac{\alpha}{\alpha+1}}](\cS_{n,m})$ under $\bP_{n,m}(\cdot\,|\,N\rb_n=n, N\rw_m=m)$ 
converges in distribution to $\cS^{(2)}$, jointly with the convergence in~\eqref{cv: ZH2}. 
\item
{\sf (Matched heavy tails)} Assume that both~\eqref{hyp: b-power} and~\eqref{hyp: w-power} hold with $\alpha=\alpha'\in (1, 2)$ and $C\rb\in (0, \infty), C\rw\in (0, \infty)$. Then $\mathfrak S[n^{\frac{\alpha}{\alpha+1}}](\cS_{n,m})$ under $\bP_{n,m}(\cdot\,|\,N\rb_n=n, N\rw_m=m)$ converges in distribution to $\cS^{(3)}$, jointly with the convergence in~\eqref{cv: ZH3}. 
\end{enumerate}
\end{prop}

\begin{proof}
See Section~\ref{sec: pf: cv-surplus}.
\end{proof}

Let us briefly explain the reason for the slight unusual definitions of $D_{n,m}$ and $\cD^{(i)}$. Roughly speaking,  the time parameter in $Z^{n,m}$ (resp.~$\cZ^{(i)}$) corresponds to the $\by$-weights of the portion of the graph already explored. However, we require $\bx$-weights in order to sample the surplus edges with the correct probabilities.  
For $D_{n,m}$, the function $\Sigma^{n,m}$ handles this issue. Meanwhile, in the proof of Proposition~\ref{prop: cv-surplus}, we show that $\Sigma^{n,m}$, suitably rescaled, converges to the linear function $t\mapsto \sigma\rb_2/\sqrt\theta t$, which explains the definition of $\cD^{(i)}$. 

\subsection{Construction of the limit graphs}

So far, we have seen in Proposition~\ref{prop: cv-ZH} that the pair $(\cZ^{(i)}, \cH^{(i)})$ appears in the scaling limit for the discrete graph encoding processes $(Z^{n,m}, H^{n,m})$, as well as in Proposition~\ref{prop: cv-surplus} that $\cQ^{(i)}$ describes the limit distribution of surplus edges, $i\in\{1, 2, 3\}$.  We explain here how the limit graphs $\cG^{(i)}$ in Theorems~\ref{thm: thm1}-\ref{thm: thm3} are obtained from the triple $(\cZ^{(i)}, \cH^{(i)}, \cQ^{(i)})$. The basic idea is that $\cH^{(i)}$ will be used to construct a forest of continuum random trees, where ``shortcuts'' will then be identified using $ \cQ^{(i)}$. 

\subsubsection{Graphs encoded by real-valued functions}
\label{sec: ghp}

We follow the approach in~\cite{BrDuWa22} to present the construction of a measured metric space from a real-valued function and a finite collection of points. The construction is slightly more general than the one seen for instance in~\cite{ABBrGo12} as it allows the function to have jumps, which is  particularly appealing for us as our discrete height processes $H^{n,m}$ are not continuous. 
We assume that our reader is familiar with the notions of real trees~\cite{Evans}, measured metric spaces~\cite{Gromov}, and Gromov--Hausdorff--Prokhorov topology~\cite{Mi09}. In particular, we adopt the following definition of Gromov--Hausdorff--Prokhorov distance: 
for two measured metric spaces $\cG_{1}=(G_{1}, d_{1}, \mu_{1})$ and $\cG_{2}=(G_{2}, d_{2}, \mu_{2})$, their Gromov--Hausdorff--Prokhorov distance is given by
\[
\dghp(\cG_{1}, \cG_{2}) := \inf \big\{ \dhaus\big(\phi_{1}(G_{1}), \phi_{2}(G_{2})\big) + \dpr\big(\mu_{1}\circ \phi_{1}^{-1}, \mu_{2}\circ \phi_{2}^{-1}\big)\big\},
\]
where the infimum is over all Polish spaces $(E, d_{E})$ and isometric embeddings $\phi_{i}$ from $G_{i}$ to $E$; $\dhaus$ stands for the Hausdorff distance of $E$, and $\dpr$ is the Prokhorov distance for the finite Borel measures on $E$, with $\mu_{i}\circ\phi_{i}^{-1}$ standing for the push-forward of $\mu_{i}$ by $\phi_{i}$, $i=1, 2$. 
We further point out that the space of (equivalence classes) of compact measured metric spaces is a Polish space in the topology induced by $\dghp$ (\cite{Mi09}). 

Let $\zeta_{h}\in (0, \infty)$ and  $h: [0, \zeta_{h}] \to \R_{+}$ be a right-continuous function with left-hand limits (i.e.~c\`adl\`ag). We shall further assume that 
\begin{itemize}
\item
either $h$ is continuous; 
\item
or $h([0, \zeta_h])$ is a finite set.
\end{itemize}
For $(s, t)\in [0, \zeta_h]^2$, let us define
\begin{equation}
\label{def: dh}
d_{h}(s, t) = h(s)+h(t) - 2 b_{h}(s, t),  \text{ where }  b_{h}(s, t) = \min\{h(u): \min(s, t)\le u\le \max(s, t)\}. 
\end{equation}
Say $s\sim_{h} t$ if and only if $d_{h}(s, t)=0$. Then $d_{h}$ induces a distance on the quotient space $T_{h}:=[0, \zeta_{h}]/\sim_{h}$, which we still denote as $d_{h}$. 
The metric space $(T_{h}, d_{h})$ is compact and tree-like as it satisfies the so-called four-point inequality: for all $s_1, s_2, s_3, s_4\in [0, \zeta_h]$, it holds
\[
d_h(s_1, s_2)+d_h(s_3, s_4)\le \max\{d_h(s_1, s_3)+d_h(s_2, s_4), d_h(s_1, s_4)+d_h(s_2, s_3)\}.
\]
On the other hand, $(T_{h}, d_{h})$ is not necessarily connected. Let  $p_{h}$ stand for the canonical projection from $[0, \zeta_{h}]$ to $T_{h}$ and denote by $\mu_{h}$ the push-forward of the Lebesgue measure on $[0, \zeta_{h}]$ by $p_{h}$. Let us define the following measured metric space
\begin{equation}
\label{def: real-tree}
\cT_{h}:= \big(T_{h}, d_{h}, \mu_{h}\big). 
\end{equation}
Informally speaking, the measured metric space that we are aiming for will be obtained from $\cT_h$ by inserting a finite number of ``shortcuts''. 
Despite $(T_h, d_h)$ not being connected in general,  we can introduce a notion of path:  for $x, y\in T_h$, a {\it path from $x$ to $y$ in $T_h$} is a finite sequence $(e_1, e_2, e_3, \dots, e_{p})$, where $p\in \N$, $e_i=(x_i, y_i)\in T_h^2$, $1\le i\le p$, $x_{i+1}=y_i$, $1\le i\le p-1$, and $x_1=x, y_{p}=y$. 
Now let $q\in \N$ and $\Pi=\{(u_i, v_i): 1\le i\le q\}$ with $(u_i, v_i)\in T_h^2$ for $1\le i\le q$. Fix some $\epsilon>0$. For a path $\gamma=(e_1, e_2, e_3, \dots, e_p)$ with $e_i=(x_i, y_i), 1\le i\le p$, its {\it $(\Pi, \epsilon)$-modified length} is defined as
\[
\ell_{\Pi, \epsilon}(\gamma)= \sum_{1\le i\le p} \ell(e_i),
\] 
where $\ell(e_i)=\min\{\epsilon, d_h(x_i, y_i)\}$ if either $(x_i, y_i)$ or $(y_i, x_i)\in \Pi$, and $\ell(e_i)=d_h(x_i, y_i)$ otherwise.  
The following defines a pseudo-distance on $T_{h}$ if $\epsilon=0$ and a distance if $\epsilon>0$: for all $x, y\in T_{h}$, let
\[
d_{h, \Pi, \epsilon}(x, y) = \inf\big\{\ell_{\Pi, \epsilon}(\gamma): \gamma \text{ is a path from $x$ to $y$ in $T_h$}\big\}.
\]
In the case where $\epsilon=0$, we then turn it into a true distance by quotienting the points at $d_{h, \Pi, \epsilon}$-distance 0 from each other, similarly to the way how $\cT_h$ is defined. In both cases, call the resulting metric space $(G_{h, \Pi, \epsilon}, d_{h, \Pi, \epsilon})$ and denote by $p_{h, \Pi, \epsilon}$ the canonical projection from $(T_{h}, d_{h})$ to $(G_{h, \Pi, \epsilon}, d_{h, \Pi, \epsilon})$. Write $\mu_{h, \Pi, \epsilon}$ for the push-forward of $\mu_{h}$ by $p_{h,\Pi, \epsilon}$. The end product of our construction is the following measured metric space
\begin{equation}
\label{def: graph}
\cG(h, \Pi, \epsilon) = \big( G_{h, \Pi, \epsilon}, d_{h, \Pi, \epsilon},  \mu_{h, \Pi, \epsilon}\big). 
\end{equation}
Denote by $a\cdot h$ the function $a\cdot h(x) =ah(x), x\in [0, \zeta]$. 
Standard arguments (see for instance Appendix C in~\cite{BrDuWa21}) yield the following bound: for $a>0$, 
\begin{equation}
\label{bd: shortcut-ghp}
\dghp\big(\cG(a\cdot h, \Pi, \epsilon), \cT_{a\cdot h}\big)\le  (2\epsilon+a)\cdot \#\Pi+a. 
\end{equation}
In the sequel, it is sometimes more convenient to define the endpoints of shortcuts from their pre-images in $[0, \zeta_h]$. Accordingly, let $q\in \N$ and $\varpi=\{(s_i, t_i): 1\le i\le q\}$ with $s_i, t_i\in [0, \zeta_h]$ for $1\le i\le q$. We then define $\Pi(\varpi)=\{(p_h(s_i), p_h(t_i)): 1\le i\le q\}$ and set
\begin{equation}
\label{def: graph'}
\mathfrak G(h, \varpi, \epsilon) = \cG(h, \Pi(\varpi), \epsilon). 
\end{equation}

\smallskip
\noindent
{\bf Comparison of two function-encoded graphs.}
As our main strategy for proving convergence of graphs relies upon their coding functions, we will utilise the following result from \cite{BrDuWa21}. Let $q\in \N$.  For $j\in \{1, 2\}$, suppose that $h_j: [0, \zeta_{j}]\to \R_{+}$ is a c\`adl\`ag function that is either continuous or $h_j([0,\zeta_j])$ is a finite set, and $\varpi_{j}=\{(s_{j, i}, t_{j, i}): 1\le i\le q\}$ is a collection of points with $0\le s_{j, i}\le t_{j, i}<\zeta_{j}$ for each $i\le q$. Suppose further that there is some $\delta>0$ verifying
\[
\max_{1\le i\le q}|s_{1, i}-s_{2, i}| \le \delta, \quad \max_{1\le i\le q}|t_{1, i}-t_{2, i}| \le \delta. 
\]
Let $\epsilon_1, \epsilon_2>0$. Then we have (\cite{BrDuWa21}, Lemma 2.7)
\begin{align}
\label{eq: ghp}
&\dghp\big(\mathfrak G(h_{1}, \varpi_{1}, \epsilon_1), \mathfrak G(h_{2}, \varpi_{2}, \epsilon_2)\big) \\ \notag
&\qquad\qquad\qquad \le 6 (q+1)\big(\|\hat h_{1}-\hat h_{2}\|_{\infty}+\omega_{\delta}(\hat h_{1})\big) +3q\max\{\epsilon_1, \epsilon_2\}+ |\zeta_{1}-\zeta_{2}|,
\end{align}
where $\hat h_{j}$ is the extension of $h_{j}$ to $\R_{+}$ by setting $\hat h_{j}(x)=0$ for all $x\ge \zeta_{j}$, and $\omega_{\delta}(\hat h_{1})=\sup\{|\hat h_{1}(s)-\hat h_{1}(t)|: |s-t|\le \delta\}$ is the $\delta$-modulus of continuity of $\hat h_{1}$. 

\subsubsection{The three limit graphs}
\label{sec: limit-gr}
 
We present here the construction of the limit graphs $\cG^{(i)}, i=\{1,2,3\}$, that appear in the main theorems.
 
\smallskip
\noindent
{\bf The double finite third moment scenario.}
Recall from~\eqref{id: girsanov} the limit process $\cZ^{(1)}$ in the double finite third moments scenario and note that~\eqref{id: Z1} says that $\cZ^{(1)}$ is a Brownian motion with parabolic drift. 
Results from Aldous~\cite{Al97} (see in particular Section 5 there) show that the excursions of $\cZ^{(1)}$ above its running infimum process can be ranked in decreasing order of their lengths. For $k\in \N$, let $(g_{k}, d_{k})$ be the $k$-th longest excursion of $\cZ^{(1)}$ above its running infimum. 
Recall from~\eqref{def: cH} as well as~\eqref{id: cH1} the height process $\cH^{(1)}$ for $\cZ^{(1)}$. 
We note that $(g_{k}, d_{k})$ is also the $k$-th longest excursion of $\cH^{(1)}$ above 0. Recall $\rho=\sigma\rb_2/\sqrt\theta$. Let $\cH^{1, k}$ be as follows:  
\[
\cH^{1, k}_t = \cH^{(1)}_{g_{k}+t/\rho}, \quad 0\le t\le \rho(d_{k}-g_{k}).
\]
The $\rho^{-1}$ factor in the time parameter is due to the fact that our encoding $H^{n,m}$ explores connected components in their $\by$-weights and the ratio of $\bx$-weight and $\by$-weight of a large component tends to $\rho$ in the limit. 
Recall from~\eqref{def: cS} the collection $\cS^{(1)}=\{(t_l, t'_l): l\ge 1\}$.  
Let us now define
\[
\varpi_k = \big\{\big(\rho(t_l- g_k), \rho(t'_l-g_k)\big): (t_l, t'_l)\in \cS^{(1)}, t_l\in (g_k,  d_k)\big\}, \quad k\ge 1. 
\]
Recall from~\eqref{def: graph'} the measured metric space $\mathfrak G(h, \varpi_h, \epsilon)$. Then $\cG^{(1)}=\{(\cC^{(1)}_k, d^{(1)}_k, \mu^{(1)}_k): k\in \N\}$ in Theorem~\ref{thm: thm1} is given by 
\begin{equation}
\label{def: cC1}
\big(\cC^{(1)}_k, d^{(1)}_k, \mu^{(1)}_k\big) = \mathfrak G\big(\cH^{1, k}, \varpi_{k}, 0), \quad k\in\N.
\end{equation}

\smallskip
\noindent
{\bf The one dominant heavy tail scenario.}
In this case, the limit process $\cZ^{(2)}$ is a semi-martingale with independent increments which further satisfies the absolute continuity relationship~\eqref{id: girsanov} with an $\alpha$-stable process. We point out that $\cZ^{(2)}$ has already been used to describe the scaling limit of critical configuration models with i.i.d.~degrees; see~\cite{Jo14} and~\cite{CG23}. In particular, in the notation of 
Theorem 8.1 in~\cite{Jo14}, we have 
\[
\big\{\tfrac{1}{\alpha+1}\cZ^{(2)}_t: t\ge 0\big\} \eqd \{X^{\nu}_t+A^{\nu}_t: t\ge 0\}
\]
with $\gamma=\alpha+2, \mu=\theta/\sigma\rb_2$ and $c=(\alpha+1)\mathrm{C}_2\cdot \mu$ in the notation there. 
As a result, Theorem 8.3 of~\cite{Jo14} applies, and  we can therefore rank the excursion intervals of $\cZ^{(2)}$ above its running infimum in decreasing order of their lengths. As in the previous case, let $(g_{k}, d_{k})$ be the $k$-th longest excursion of $\cZ^{(2)}$ above its running infimum, $k\in \N$. 
The height process $\cH^{(2)}$ is introduced in~\eqref{def: cH}. The definition there implies that $(g_{k}, d_{k})$ is also the $k$-th longest excursion of $\cH^{(2)}$ above 0.
The rest of the construction is done in exactly the same way as in the double finite third moments regime, except that we replace $\cZ^{(1)}$ with $\cZ^{(2)}$, $\cH^{(1)}$ with $\cH^{(2)}$ and $\cS^{(1)}$ with $\cS^{(2)}$.

\smallskip
\noindent
{\bf The matched heavy tail scenario.} The semi-martingale representation~\eqref{id: cZ-drift} implies that 
\[
\big\{\mathrm{C}_2^{-1}\cZ^{(2)}_t: t\ge 0\big\} \eqd \big\{\mathrm{C}_3^{-1}\cZ^{(3)}_t: t\ge 0\big\}. 
\]
Therefore, the construction in the previous regime can be easily adapted to the current one; we omit the detail. 

\subsubsection{Connection with the Binomial model}

Let us take $X_i\equiv 1, 1\le i\le n$ and $Y_j\equiv 1, 1\le j\le m$. As mentioned  in the Introduction, in this case, $B_{n,m}$ corresponds to the Binomial bipartite graph with $n$ black vertices, $m$ white vertices and each edge present with probability $p_n=1-\exp(-1/\sqrt{mn})$, studied in~\cite{W25}. In particular, Theorem 2.1 combined with Proposition 2.4 there establishes a limit theorem for $(C^{\bx}_{n, (k)}, \tfrac12 n^{-1/3}\dgr, n^{-2/3}\mu^{\bx}_{n, k}), k\ge 1$. Let us argue that the limit object obtained in~\cite{W25} has the same distribution as $\cG^{(1)}$. More specifically, we will show that both are encoded by the same Brownian motion with parabolic drift. On the one hand, Proposition 2.8 (taking $\lambda=0$) combined with (2.19) in~\cite{W25} says that the height process encoding the limit graph is given by 
\[
\cH^{0, \theta}_t=\frac{2}{\nu_{\theta}}\big(\cS^{0, \theta}_t-\mathcal I^{0, \theta}_t), \quad \text{with } \cS^{0, \theta}_t=\sqrt{\nu_{\theta}}B_t-\tfrac12 \nu_{\theta}t^2, \ \cI^{0, \theta}_t=\min_{s\le t}\cS^{0, \theta}_t, \ \nu_{\theta}=1+\theta^{-\frac12}
\]
In above, $(B_t)_{t\ge 0}$ denotes a standard linear Brownian motion. Writing $\gamma_{\theta}=1+\sqrt\theta = \sqrt\theta\nu_{\theta}$, we find 
\[
\cH^{0, \theta}_t=2\big(\widehat{\cS}_t-\widehat{\mathcal I}_t), \quad \text{with } \widehat{\cS}_t = \theta^{\frac14}\gamma_\theta^{-\frac12}B_t - \tfrac12 t^2\  \text{ and } \ \widehat{\cI}_t = \min_{s\le t}\widehat{\cS}_s. 
\]
On the other hand, from~\eqref{def: L} and~\eqref{id: Z1}, we find that $\cZ^{(1)}_t$ has the same distribution as
$\sqrt{\gamma_\theta}B_t - \gamma_\theta t^2/(2\theta)$,  
since $\sigma\rb_r=\sigma\rw_r=1$ for all $r\ge 1$. It then follows from~\eqref{id: cH1} that 
\[
\cH^{(1)}_t \eqd 2\big(\widehat{\cS}_{t/\sqrt\theta}-\widehat{\mathcal I}_{t/\sqrt\theta}). 
\]
Finally, we note that from the construction in Section~\ref{sec: limit-gr} that the limit graph $\cG^{(1)}$ is encoded by $(\cH^{(1)}_{t/\rho})_{t\ge 0}$ with $\rho=1/\sqrt\theta$ in the current case. The conclusion follows.

\section{Proof}
\label{sec: proof}

This proof section is organised as follows. In Section~\ref{sec: x-graph}, we prove the various properties concerning the LIFO-queue construction. Section~\ref{sec: pf-cv-ZH} contains the proof of Proposition~\ref{prop: cv-ZH} about the convergence of $(Z^{n,m}, H^{n,m})$. We then proceed to the convergence of the surplus edges in Section~\ref{sec: cv-surplus}. Finally, in Section~\ref{sec: pf-main}, we gather all the ingredients and give the proof of our main theorems.  

\subsection{Some results on the bipartite graph with fixed weights}
\label{sec: x-graph}

Recall from Introduction the bipartite graph $B(n,m;\bx, \by)$ with the fixed weight sequences $\bx=(x_i)_{1\le i\le n}$ and $\by=(y_j)_{1\le j\le m}$. We give here the proofs of Lemmas~\ref{lem: B-graph}-\ref{lem: exc} that concern this model.  
Clearly, the same properties also hold for $B_{n,m}$, the model with i.i.d.~weights. 

\begin{proof}[Sketch proof for Lemma~\ref{lem: B-graph}]
\label{pf: B-graph}
The main idea is similar to the one used in~\cite{BrDuWa22} for rank-1 models. A formal argument in that case can be found in Section 3 of~\cite{BrDuWa22}. We sketch here an informal argument, primarily based upon the memoryless properties of the exponential distribution. 
Let $\hat I_{i,j}$ denote the indicator for the event that the black vertex $b_i$ is adjacent to the white vertex $w_j$ in the graph $B'$. By definition, it suffices to show that $\hat I_{i,j}, i\in [n], j\in [m]$ are independent Bernoulli variables with respective means $p_{i, j}$ as defined in~\eqref{def: edge-prob}. 
Say that at step $k$, $V_k=b_i$ has been assigned. Among the $m$  white vertices, we distinguish two subsets: $\mathcal W$, which consists of those $j$ satisfying $E\rw_j>\tau\rw_{k-1}$, and the complement of $\mathcal W$. Conditional on $w_j\in \mathcal W$, by the memoryless property, the probability that $w_j\in \cN(V_k)$ is then precisely $p_{i, j}$.  Therefore, conditional on $w_j\in \mathcal W$ and $V_k=b_i$, $\hat I_{i, j}$ occurs with probability $p_{i,j}$. 
If $w_j\notin \mathcal W$, then either $w_j$ has left the queue prior to the arrival of $V_k$, or $w_j$ is still in the queue, i.e.~$w_j\in \cA_{k-1}$. 
In the first case, we must have $w_j=U_{k'}$ for some $k'<k$ and $w_j$ receives its full service by the end of step $k'$. This means that all the neighbours of $w_j$ have been identified by then and $V_k$ is not among them. Since the total service time of $w_j$ is $y_j$, this occurs with probability $1-p_{i,j}$. 
In the second case, either $w_j=U_k$ or $w_j\in \cA_{k-1}\setminus\{U_k\}$. If $w_j=U_k$, then $V_k$ has been identified as a neighbour of $w_j$, which occurs with probability $p_{i, j}$. 
If $w_j\ne U_k$, then 
recall that $y_j(k)$ denotes the remaining amount of service that $w_j$ has yet to receive when $V_k$ arrives. During the $\Delta:=y_j-y_j(k)$ amount of service that $w_j$ received, $V_k$ has not been identified as a neighbour,  which is an event of probability $\exp(-x_i\cdot \Delta/z_{n,m})$. The conditional distribution of $\hat I_{i,j}$ given the event is then precisely the one given in~\eqref{def: sur-prob}. Finally, let us point out in our definition of $S_k$, we do not exclude the pair $(V_k, U_k)$. Therefore,  $\cE_{n,m}$ can contain $(V_k, U_k)$ with a positive probability. However, since the edge already exists in $\cF$ and $B'$ is a simple graph, this does not change the distribution of $B'$. 
\end{proof}

\begin{proof}[Proof of Lemma~\ref{lem: Z}]
\label{pf: Z}
Let us write $\hat Z_t = -t + \Lambda^{\by}\circ\Lambda^{\bx}(t)$. Since both have the same linear drift, it suffices to show that $Z^{\bx, \by}$ and $\hat Z$ share the same jumps. 
To that end, we first note that as $\Lambda^{\bx}$ only increases through jumps, the jumping times of $\hat Z_t$ are necessarily jumping times of $\Lambda^{\bx}$, which are $E\rb_i, 1\le i\le n$. It remains to show that for each $i\in [n]$, 
\[
\hat\Delta_i:=\hat Z_{E\rb_i}-\hat Z_{E\rb_i-} = \Delta^{\mathrm{bi}}_i.
\] 
By definition, we have 
\[
\hat \Delta_i = \sum_{j\in [m]} y_j\mathbf 1_{\{\Lambda^{\bx}(E\rb_i-)<E\rw_j\le \Lambda^{\bx}(E\rb_i)\}}.
\]
Recall from the LIFO-queue construction that the black vertex $b_i$ appears as certain $V_k$ during the construction. In other words, there is a unique $k=k(i)$ such that $V_{k}=b_i$. Recall also the way that the white dial is updated: we have $\tau\rw_j-\tau\rw_{j-1}=x_{V_j}$ if $V_j$ is defined and $\tau\rw_j-\tau\rw_{j-1}=0$ if no $V_j$ is appointed at step $j$. It follows that 
\[
\tau\rw_{k}=\sum_{j\le k}x_{V_{j}}\mathbf 1_{\{V_{j} \text{ is defined}\}}. 
\]
Since the order in which $V_j, j\ge 1$, appear coincides with the ranking of $(E\rb_i)_{i\in [n]}$, we deduce that  
\[
\Lambda^{\bx}(E\rb_i-) =\tau\rw_{k-1} \quad\text{and}\quad \Lambda^{\bx}(E\rb_i) = \tau\rw_{k}. 
\]
Hence, we have 
\[
\hat \Delta_i = \sum_{j\in [m]}y_j\mathbf 1_{\{E\rw_j\in (\tau\rw_{k-1}, \tau\rw_k]\}} = \sum_{j\in [m]}y_j\mathbf 1_{\{E\rw_j\in I(V_k)\}}=\sum_{j\in [m]}y_j\mathbf 1_{\{E\rw_j\in I(i)\}}=\Delta^{\mathrm{bi}}_i,
\]
according to the definition~\eqref{def: Delta}. This completes the proof. 
\end{proof}

\begin{proof}[Proof of Lemma~\ref{lem: exc}]
\label{pf: exc}
Recall $Z^Q_t$ from~\eqref{def: ZQ}, which tracks the service load. In particular, the server is idle at time $t$ if and only if $Z^Q_t = \inf_{s\le t} Z^Q_s$. 
Moreover, after the first customer of the queue arrives, the server will only become idle after having served out all the clients in the first tree. 
It follows that if 
$(g^Q_1, d^Q_1)$ is the first excursion interval of $Z^Q$ above its running infimum, then $d^Q_1-g^Q_1$ is equal to the total service requests from the customers in the first tree. 
Iterating this argument for the subsequent excursion intervals, we obtain a bijection between the excursion lengths of $Z^Q$ and the total service requests from each tree. 
We now apply this to $Z^{\bx,\by}$. To that end, first note from~\eqref{def: Delta} that the total service request from each tree correspond to its $\by$-weight, which is non zero if and only if the tree is non trivial. Therefore, we have $\kappa_{n,m}=\kappa'_{n,m}$. Furthermore, applying the previous bijection to $Z^{\bx, \by}$ yields the first identity in~\eqref{id: mass}. 
For the second, observe that $\Lambda^{\bx}$ and $Z^{\bx, \by}$ share the same jump times. From the definition of the bipartite LIFO-queue, we see that the jumps of $Z^{\bx, \by}$ that appear in $[g_{n,i}, d_{n,i})$ are precisely the black vertices of $T_{n, i}$. The second identity in~\eqref{id: mass} now follows from the definition of $\Lambda^{\bx}$.
\end{proof}

\subsection{Proof of Proposition~\ref{prop: cv-ZH}}
\label{sec: pf-cv-ZH}

We prove in this subsection Proposition~\ref{prop: cv-ZH}, which asserts the convergence of the graph encoding processes $(Z^{n,m}, H^{n,m})$ in the three different scenarios. 
Our approach here is inspired by the work of Conchon-Kerjan and Goldschmidt~\cite{CG23}, although in our case $Z^{n,m}$ is itself a composition of exploration-type processes. On the other hand, the exploration processes we introduce here have more regular distributional properties, which will significantly simplify the calculations for their convergence. 
Let us first lay out the main steps of the proof.

\subsubsection{An overview of the proof} 

Recall from Section~\ref{sec: coupling} the coupling of the LIFO-queue construction with a pair of Poisson point measures. In particular, we have 
\[
B_{n,m} \eqd B(n,m; \tilde{\mathbf X}_n, \tilde{\mathbf Y}_n) \text{ under } \bP_{n,m}(\cdot\,|\, N\rb_n=n, N\rw_m=m).
\]
The following lemma says that in fact we can remove the conditioning on $N\rb_n$ and $N\rw_m$.  Roughly speaking, this is because we are primarily concerned with large connected components, which appear early in the exploration and only rely upon the first $o(n)$ atoms in the Poisson point measures. Recall $\Lambda\rbn$ from~\eqref{def: Lam} and $(Z^{n,m}, H^{n,m})$ from~\eqref{def: ZH}. 
 
\begin{lem}
\label{lem: cond}
Assume~\eqref{hyp: mod-clst}. 
Let $(b_n)_{n\in \N}$ be a sequence of positive real numbers that satisfies $b_n\to\infty$ and $b_n=o(n)$ as $n\to\infty$. Let $t_0\in (0, \infty)$. 
Assume that the sequence $b_n^{-1}\Lambda\rbn_{b_nt_0}, n\ge 1$, is tight. 
Denote by $\mu_n$ the law of $(Z^{n,m}_{b_n t}, H^{n,m}_{b_n t})_{t\le t_0}$ under $\bP_{n,m}$ and by $\hat\mu_n$ the law of the same pair under $\bP_{n,m}(\cdot\,|\, N\rb_n=n, N\rw_m=m)$. Then
\[
\dTV(\hat\mu_n, \mu_n):=\sup_{A}\big|\hat\mu_n(A)-\mu_n(A)\big| \to 0, \quad \text{as } n\to\infty,
\]
where  the supremum is over all the Borel sets $A$ of the Skorokhod space $\mathbb D([0,t_0], \R^2)$. 
\end{lem}

\begin{proof}
See Section~\ref{sec: pf-lem-cond}.
\end{proof}

Thanks to Lemma~\ref{lem: cond}, we only need to prove the claimed convergences in Proposition~\ref{prop: cv-ZH} under $\bP_{n,m}$. We divide the remaining proof into the following steps. 

\smallskip
\noindent
{\bf Step 1: Radon--Nikodym derivatives for the graph exploration processes.}
We note that under $\bP_{n,m}$, $Z^{n,m}$ has independent but inhomogeneous increments. We establish here an analogue of~\eqref{id: girsanov} for $Z^{n,m}$. To that end, let us first introduce two mutually independent compound Poisson processes $(L\rb_t)_{t\ge 0}$ and $(L\rw_t)_{t\ge 0}$ under $\bP_{n,m}$ with respective L\'evy measures $\sqrt{n/m}\,xdF\rb(x)$ and $\sqrt{m/n}\,ydF\rb(y)$. In other words, we have for all $t\ge 0, \lambda\ge 0$: 
\begin{equation}
\label{def: Lbw}
\mathbf E_{n,m}[e^{-\lambda L\rb_t}]=\exp(t\,\varphi\rb(\lambda)), \quad \mathbf E_{n,m}[e^{-\lambda L\rw_t}]=\exp(t\,\varphi\rw(\lambda)),
\end{equation}
where the Laplace exponents $\varphi\rb$ and $\varphi\rw$ are given by
{\small\begin{equation}
\label{def: varphibw}
\varphi\rb(\lambda) = \sqrt{\frac{n}{m}} \int_{(0,\infty)}(e^{-\lambda x}-1) x dF\rb(x), \quad \varphi\rw(\lambda) =\sqrt{\frac{m}{n}} \int_{(0,\infty)}(e^{-\lambda x}-1) x dF\rw(x). 
\end{equation}}
We then set 
\begin{equation}
\label{def: Lnm}   
L^{n,m}_t := -t + L\rw\circ L\rb(t):=-t+L\rw_{L\rb_t}, \quad t\ge 0.
\end{equation}
We point out that $(L^{n,m}_t)_{t\ge 0}$ is itself a L\'evy process, by standard results on the subordination of L\'evy processes (see Ch.~6 in~\cite{Sato}).  

\begin{lem}
\label{lem: girsa}
Let $n,m\in \N$. The following identity holds for all $t\ge 0$ and all measurable function $F:\mathbb D([0, t], \R)\to \R_+$: 
\begin{equation}
\label{id: gir-dis}
\bE_{n,m}\Big[F\Big(\big(Z^{n,m}_s\big)_{s\le t}\Big)\Big] = \bE_{n,m}\Big[F\Big(\big(L^{n,m}_s\big)_{s\le t}\Big)\cdot E^{n,m}_t\Big],
\end{equation}
where 
{\small\begin{equation}
\label{def: Enm}
E^{n,m}_t = \exp\Big\{-\int_0^{L\rb_t} \!\!\!\!\!\frac{s}{\sqrt{mn}} \,d L\rw_s - \int_0^{L\rb_t} \!\!\!\!\!\varphi\rw\Big(\frac{s}{\sqrt{mn}}\Big)ds-\int_0^{t}\!\!\! \frac{s}{\sqrt{mn}} \,d L\rb_s- \int_0^{t}\!\!\! \varphi\rb\Big(\frac{s}{\sqrt{mn}}\Big)ds\Big\}.
\end{equation}}
\end{lem}

\begin{proof}
See Section~\ref{sec: pf-girsa}.
\end{proof}

\smallskip
\noindent
{\bf Step 2: Convergence of the L\'evy processes and their height processes.} 
In Section~\ref{sec: xy-graph-coding}, we have introduced the height process $H^{n,m}_t$ as the following functional of $Z^{n,m}$:
\[
H^{n,m}_t = \#\Big\{s\le t: Z^{n,m}_{s-} < \inf_{s\le u\le t}Z^{n,m}_u\Big\}.
\]
We now introduce its counterpart for $L^{n,m}$ using the same functional: 
\[
\widehat H^{n,m}_t = \#\Big\{s\le t: L^{n,m}_{s-} < \inf_{s\le u\le t}L^{n,m}_u\Big\}, \quad t\ge 0.
\]

\begin{lem}
\label{lem: cv-levy}
Assume~\eqref{hyp: mod-clst} and~\eqref{hyp: crit}. 
For $i\in \{1, 2, 3\}$, let $(\cL^{(i)}_t)_{t\ge 0}$ be the L\'evy process specified in~\eqref{def: L} and $(\widehat\cH^{(i)}_t)_{t\ge 0}$ be as in~\eqref{def: Hlevy}. 
\begin{enumerate}
\item
{\sf (Double finite third moments)}
Under both~\eqref{hyp: b-third} and~\eqref{hyp: w-third}, the following weak convergence takes place under $\bP_{n,m}$: 
\begin{equation}
\label{cv: LH1}
    \Big\{n^{-\frac13}L^{n,m}_{n^{2/3}t}, n^{-\frac13}\widehat H^{n,m}_{n^{2/3}t}: t\ge 0\Big\} \Longrightarrow \Big\{\cL^{(1)}_t, \widehat\cH^{(1)}_t: t\ge 0\Big\} \text{ in } \mathbb D(\R_+, \R^2). 
\end{equation}

\item
{\sf (One dominant heavy tail regime)}
Assume that~\eqref{hyp: b-power} holds with $\alpha\in (1, 2)$ and $C\rb\in (0, \infty)$. Assume that either~\eqref{hyp: w-third} is true or~\eqref{hyp: w-power} holds with $\alpha'\in (\alpha, 2)$ and $C\rw\in (0, \infty)$. Then  the following weak convergence takes place under $\bP_{n,m}$: 
\begin{equation}
\label{cv: LH2}
    \Big\{n^{-\frac{1}{\alpha+1}}L^{n,m}_{ n^{\alpha/(\alpha+1)}t}, n^{-\frac{\alpha-1}{\alpha+1}}\widehat H^{n,m}_{ n^{\alpha/(\alpha+1)}t}: t\ge 0\Big\}\Longrightarrow \Big\{\cL^{(2)}_{t}, \widehat\cH^{(2)}_t: t\ge 0\Big\} \text{ in } \mathbb D(\R_+, \R^2). 
\end{equation}

\item
{\sf (Matched heavy tails regime)}
Assume that both~\eqref{hyp: b-power} and~\eqref{hyp: w-power} hold with $\alpha=\alpha'\in (1, 2)$ and $C\rb\in (0, \infty), C\rw\in (0, \infty)$. The following weak convergence holds under $\bP_{n,m}$: 
\begin{equation}
\label{cv: LH3}
    \Big\{n^{-\frac{1}{\alpha+1}}L^{n,m}_{ n^{\alpha/(\alpha+1)}t}, n^{-\frac{\alpha-1}{\alpha+1}}\widehat H^{n,m}_{ n^{\alpha/(\alpha+1)}t}: t\ge 0\Big\}\Longrightarrow \Big\{\cL^{(3)}_{t}, \widehat\cH^{(3)}_t: t\ge 0\Big\} \text{ in } \mathbb D(\R_+, \R^2). 
\end{equation}
\end{enumerate}
\end{lem}

Our proof of Lemma~\ref{lem: cv-levy}, given in Section~\ref{sec: pf-cv-levy}, relies heavily on Duquesne--Le Gall's Theorem that asserts the convergence for the contour functions of critical Bienaym\'e trees;  see Theorem 2.4.1 in~\cite{DuLG02}. However, for our purpose here, we require a modified version--stated in Appendix~\ref{sec: DuLGthm}--of the original theorem. 

\smallskip
\noindent
{\bf Step 3: Convergence of the Radon--Nikodym derivatives.}
The final ingredient for the proof of Proposition~\ref{prop: cv-ZH} concerns the convergence of $E^{n,m}_t$ in~\eqref{def: Enm}. Recall $\cE^{(i)}_t$ from~\eqref{def: cE}  and $\Lambda\rbn$ from~\eqref{def: Lam}.
\begin{lem}
\label{lem: cv-dens}
For each $t\ge 0$, the following statements hold true under $\bP_{n,m}$: 
\begin{enumerate}
\item
Under the assumptions of Lemma~\ref{lem: cv-levy} (1), $E^{n,m}_{n^{2/3}t}$ is uniformly integrable and converges in distribution to $\cE^{(1)}_t$, jointly with the convergence in~\eqref{cv: LH1}. Furthermore, for all $t_0\ge 0$, we have the following convergence in probability: 
\begin{equation}
\label{LLN1}
\sup_{t\le t_0}\Big|n^{-\frac23}\Lambda\rbn(n^{\frac23} t) - \theta^{-\frac12}\sigma\rb_2 t\Big|\to 0.
\end{equation}

\item
Under the assumptions of Lemma~\ref{lem: cv-levy} (2), $E^{n,m}_{n^{\alpha/(\alpha+1)}t}$ is uniformly integrable and converges in distribution to $\cE^{(2)}_t$, jointly with the convergence in~\eqref{cv: LH2}. Furthermore, for all $t_0\ge 0$, we have the following convergence in probability: 
\begin{equation}
\label{LLN2}
\sup_{t\le t_0}\Big|n^{-\frac{\alpha}{\alpha+1}}\Lambda\rbn(n^{\frac{\alpha}{\alpha+1}} t) - \theta^{-\frac12}\sigma\rb_2 t\Big|\to 0.
\end{equation}

\item
Under the assumptions of Lemma~\ref{lem: cv-levy} (3), $E^{n,m}_{ n^{\alpha/(\alpha+1)}t}$ is uniformly integrable and converges in distribution to $\cE^{(3)}_t$, jointly with the convergence in~\eqref{cv: LH3}. Furthermore, for all $t_0\ge 0$, we have the following convergence in probability: 
\begin{equation}
\label{LLN3}
\sup_{t\le t_0}\Big|n^{-\frac{\alpha}{\alpha+1}}\Lambda\rbn(n^{\frac{\alpha}{\alpha+1}} t) - \theta^{-\frac12}\sigma\rb_2 t\Big|\to 0.
\end{equation}

\end{enumerate}
\end{lem}

\begin{proof}
See Section~\ref{sec: pf-cv-dens}.
\end{proof}

To conclude this overview, we explain how the previous results will yield Proposition~\ref{prop: cv-ZH}. 

\medskip

\begin{proof}[Proof of Proposition~\ref{prop: cv-ZH} subject to Lemmas~\ref{lem: cond}-\ref{lem: cv-dens}]
We use $a_n, b_n$ to denote the appropriate spacial and temporal scaling for $Z^{n,m}$ in each of the three scenarios, namely $a_n=n^{1/3}, b_n=n^{2/3}$ under the double finite third moments assumptions and $a_n =  n^{1/(\alpha+1)}, b_n =  n^{\alpha/(\alpha+1)}$ in the other two scenarios. 
Let $t_0\ge 0$ and $F: \mathbb D([0, t_0], \R^2)\to \R$ be continuous and bounded. 
Lemma~\ref{lem: girsa} applies to yield that
\begin{align*}
&\bE_{n,m}\Big[F\Big(\big(\tfrac{1}{a_n} Z^{n,m}_{b_n t}\big)_{t\le t_0}, \big(\tfrac{a_n}{b_n}H^{n,m}_{b_nt}\big)_{t\le t_0}\Big)\Big] \\
&\qquad\qquad = \bE_{n,m}\Big[F\Big(\big(\tfrac{1}{a_n} L^{n,m}_{b_n t}\big)_{t\le t_0}, \big(\tfrac{a_n}{b_n}\widehat H^{n,m}_{b_nt}\big)_{t\le t_0}\Big)\cdot E^{n,m}_{b_n t} \Big].
\end{align*}
Let $K\in (0, \infty)$. Thanks to Lemmas~\ref{lem: cv-levy} and~\ref{lem: cv-dens} , we find that for $i\in \{1, 2, 3\}$ and under the relevant assumptions, 
\begin{align*}
&\lim_{K\to\infty}\lim_{n\to\infty} \bE_{n,m}\Big[F\Big(\big(\tfrac{1}{a_n} L^{n,m}_{b_n t}\big)_{t\le t_0}, \big(\tfrac{a_n}{b_n}\widehat H^{n,m}_{b_nt}\big)_{t\le t_0}\Big)\cdot E^{n,m}_{b_n t}\wedge K \Big]\\ 
&\qquad\qquad\qquad = \mathbb E\Big[F\Big(\big(\cL^{(i)}_t\big)_{t\le t_0}, \big(\widehat \cH^{(i)}_t\big)_{t\le t_0}\Big)\cdot \cE^{(i)}_t \Big]
 = \mathbb E\Big[F\Big(\big(\cZ^{(i)}_t\big)_{t\le t_0}, \big(\cH^{(i)}_t\big)_{t\le t_0}\Big)\Big],
\end{align*}
where in the last identity we have used~\eqref{id: girsanov} and the observation that $\widehat\cH^{(i)}_t$ and $\cH^{(i)}_t$ are the same functional of $\cL^{(i)}$ and $\cZ^{(i)}$ respectively. Combined with the previous arguments, this proves the convergences~\eqref{cv: ZH1}-\eqref{cv: ZH3} under $\bP_{n,m}$. Lemma~\ref{lem: cond} then allows us to conclude that the same convergences also hold under $\bP_{n,m}(\cdot\,|\,N\rb_n=n, N\rw_m=m)$, noting that the tightness assumption in the lemma is ensured by \eqref{LLN1}, \eqref{LLN2}, \eqref{LLN3} respectively. 
\end{proof}

\subsubsection{Proof of Lemma~\ref{lem: cond}} 
\label{sec: pf-lem-cond}

We start with an observation on Poisson distributions. 
\begin{lem}
\label{lem: cond'}
Suppose that $(r_n)_{n\in\N}$ and $(s_n)_{n\in \N}$ are two sequences of positive real numbers that satisfy $r_n\to\infty$ and $r_n= o(s_n)$ as $n\to\infty$. For each $n\in \N$, let $N^n_1$ and $N^n_2$ be two independent Poisson random variables of respective means $r_n$ and $s_n$ under $\mathbb P$.   
Then
\[
\sum_{k\ge 0}\Big|\mathbb P\big(N^n_1=k\big) - \mathbb P\big(N^n_1=k\,\big|\, N^n_1+N^n_2=\lfloor r_n+s_n\rfloor \big)\Big|\to 0, \quad \text{as } n\to\infty.
\]
\end{lem}

\begin{proof}
For $K\in (0, \infty)$, let us denote the set 
\[
\Delta_{n, K} = \big\{k\in \Z_+: |k-r_n|\le K\sqrt{r_n}\big\}.
\]
Thanks to the Central Limit Theorem, we have
\[
\lim_{K\to\infty}\liminf_{n\to\infty}\mathbb P\big(N^n_1\in\Delta_{n, K}\big)=1.
\]
In view of this, we note that the conclusion will follow once we show that for each $K>0$, 
\begin{equation}
\label{eqbd: unifev}
\limsup_{n\to\infty}\sum_{k\in\Delta_{n,K}} \Big|\mathbb P\big(N^n_1=k\big)-\mathbb P\big(N^n_1=k\,\big|\,N^n_1+N^n_2=\lfloor r_n+s_n\rfloor \big)\Big|= 0.
\end{equation}
Let us denote $t_n=r_n+s_n$ and $N^n_3 = N^n_1+N^n_2$, which follows a Poisson distribution of mean $t_n$. 
To the end of proving~\eqref{eqbd: unifev}, let us first argue that
\begin{equation}
\label{eq: unifk}
\sup_{k\in \Delta_{n, K}}\Big|\frac{\mathbb P(N^n_2=\lfloor t_n\rfloor -k)}{\mathbb P(N^n_3=\lfloor t_n\rfloor)}-1\Big|\to 0, \quad\text{as } n\to\infty.
\end{equation}
To see why this is true, let us recall Stirling's formula (\cite{Feller-vol1}, II.9.15 on p.54), which says 
\begin{equation}
\label{id: stir}
\forall\, n\in \N:\quad n! =\sqrt{2\pi n} \exp\Big(n\log \frac{n}{e}+\epsilon_n\Big), \quad\text{with }\  0\le \epsilon_n\le \frac{1}{12n}. 
\end{equation}
Let $z=\lfloor t_n\rfloor-k-s_n=r_n-k+\lfloor t_n\rfloor-t_n$. Then $k\in \Delta_{n,K}$ implies that $z\in \hat\Delta_{n, K}:=\{x\in \R: |x|\le K\sqrt{r_n}+1\}$. We deduce from \eqref{id: stir} that 
\begin{align*}
&\mathbb P(N^n_2=\lfloor t_n\rfloor-k)=\mathbb P(N^n_2=s_n+z)=e^{-s_n}\frac{s_n^{s_n+z}}{(s_n+z)!} \\
&\qquad\qquad=  \frac{1}{\sqrt{2\pi(s_n +z)}} \exp\Big\{-s_n + (s_n +z)\log s_n - (s_n+z)\log \frac{s_n +z}{e}-\epsilon_n(z)\Big\}\\
&\qquad\qquad= \frac{1}{\sqrt{2\pi(s_n +z)}} \exp\Big\{z-(s_n+z)\log \Big(1+\frac{z}{s_n}\Big)-\epsilon_n(z)\Big\},
\end{align*}
where $0\le \epsilon_n(z)\le(12(s_n+z))^{-1}$; in particular, $\sup_{z\in \hat\Delta_{n,K}} |\epsilon_n(z)|\to 0$. 
Let us write $\log(1+x) = x -\frac12 x^2+ r(x)x^3$, where the function $r$ is uniformly bounded on $[-1/2, 1/2]$. 
We find
\[
\exp\Big\{z-(s_n+z)\log \Big(1+\frac{z}{s_n}\Big)\Big\}=
 \exp\Big\{-\frac12\frac{z^2}{s_n}+\frac12 \frac{z^3}{s_n^2}-r\Big(\frac{z}{s_n}\Big)\Big(\frac{z^3}{s_n^2}+\frac{z^4}{s_n^3}\Big)\Big\}.
\]
With the assumption that $r_n=o(s_n), n\to\infty$, we deduce that 
\[
\sup_{z\in \hat\Delta_{n, K}}\Big|\exp\Big\{z-(s_n+z)\log \Big(1+\frac{z}{s_n}\Big)-\epsilon_n(z)\Big\}-1\Big|\to 0, \quad n\to\infty.
\]
Since $\sqrt{2\pi (s_n+z)}\sim \sqrt{2\pi s_n}$ uniformly for $z\in \hat\Delta_{n, K}$, 
it follows that 
\begin{equation}
\label{as: nom}
\sup_{k\in\Delta_{n, K}}\Big|\mathbb P(N^n_2=\lfloor t_n\rfloor-k) -\frac{1}{\sqrt{2\pi s_n}}\Big|\to 0, \quad n\to\infty. 
\end{equation}
A similar but simpler argument shows that
\[
\mathbb P\big(N^n_3=\lfloor t_n\rfloor\big) \sim \frac{1}{\sqrt{2\pi t_n}}, \quad n\to\infty. 
\]
Together with~\eqref{as: nom} and the fact that $t_n=r_n+s_n\sim s_n$ as $n\to\infty$, we deduce the convergence in~\eqref{eq: unifk}. 
Back to~\eqref{eqbd: unifev}, we note that the independence between $N^n_1$ and $N^n_2$ implies that  
\[
\!\!\sum_{k\in \Delta_{n, K}}\!\!\!\Big|\mathbb P\big(N^n_1=k\big)-\mathbb P\big(N^n_1=k\,\big|\,N^n_3=\lfloor t_n\rfloor \big)\Big| 
= \!\!\sum_{k\in \Delta_{n, K}}\!\!\!\!\mathbb P\big(N^n_1=k\big)\cdot \bigg|1-\frac{\mathbb P(N^n_2=\lfloor t_n\rfloor -k)}{\mathbb P(N^n_3=\lfloor t_n\rfloor)}\bigg|.
\]
In view of~\eqref{eq: unifk}, the above converges to 0 by the dominated convergence theorem. 
By the previous arguments, the proof is complete. 
\end{proof}

Recall $\Lambda\rbn$ and $\Lambda\rwm$ from~\eqref{def: Lam}. Recall that $\dTV$ stands for the total variation distance between probability measures. As an application of Lemma~\ref{lem: cond'}, let us show the following. 

\begin{lem}
\label{lem: Lambda-cond}
Assume~\eqref{hyp: mod-clst}. 
Let $(b_n)_{n\in \N}$ be a sequence of positive real numbers that satisfies $b_n\to\infty$ and $b_n=o(n)$ as $n\to\infty$. 
Let $t_0\in (0, \infty)$. Denote by $\nu_n\rb$ the law of the collection
\[
\cN_{t_0}:=\big\{(E\rb_i, X_i): 1\le i\le N\rb_n, E\rb_i\le b_nt_0\big\}
\]
under $\bP_{n,m}$ and by $\hat\nu_n\rb$ the law of $\cN_{t_0}$ under $\bP_{n,m}(\cdot\,|\,N\rb_n=n, N\rw_m=m)$. Similarly, let $\nu_m\rw$ be the law of the collection $\{(E\rw_j, Y_j): 1\le j\le N\rw_m, E\rw_j\le b_n t_0\}$ under $\bP_{n,m}$, and let  $\hat\nu_m\rw$ be its law under $\bP_{n,m}(\cdot\,|\,N\rb_n=n, N\rw_m=m)$. 
Then
\begin{equation}
\label{dtv-b}
\dTV(\nu_n\rb, \hat\nu_n\rb)\to 0\quad \text{and}\quad \dTV(\nu\rw_m, \hat\nu\rw_m)\to 0, \quad\text{as } n\to\infty.
\end{equation}
\end{lem}

\begin{proof}
We focus on the first convergence, as the other one can be similarly argued. By the independence between the two Poisson point measures, $\hat\nu\rb_n$ is in fact the law of $\cN_{t_0}$ conditioned only on $N\rb_n=n$. 
Let us denote 
\[
N^n_1=\#\mathcal N_{t_0} \; \text{ and }\; N^n_2 = N\rb_n-N^n_1.
\]
We note that conditional on $N^n_1=k\in \N$, $\cN_{t_0}$ is distributed as a collection of $k$ i.i.d.~pairs of random variables with the joint distribution $\pi\rb(dt, dx)/n$. The proof boils down to showing 
\begin{equation}
\label{bd: pbth}
\dTV\Big(\mathbb P(N^n_1\in \cdot),\;  \mathbb P(N^n_1\in\cdot \,|\, N\rb_n=n)\Big) \to 0, \quad n\to\infty. 
\end{equation}
We note that $N^n_1, N^n_2$ are two independent Poisson variables with respective means: 
\[
r_n:=\mathbb E[N^n_1] = n\int_{(0, \infty)}\Big(1-e^{-\frac{xb_n t_0}{\sqrt{mn}}}\Big)dF\rb(x)\le \sqrt{\frac{n}{m}}\,b_n t_0\sigma\rb_1, \  
s_n:= \mathbb E[N^n_2]=n-r_n,
\]
where we have used the elementary inequality $1-e^{-x}\le x$ for $x\ge 0$. Our assumption on $(b_n)_{n\in\N}$ implies that $r_n\to \infty$ and $r_n=o(s_n)$, as $n\to\infty$. Lemma~\ref{lem: cond'} applies to yield~\eqref{bd: pbth}. The conclusion follows. 
\end{proof}

\begin{proof}[Proof of Lemma~\ref{lem: cond}]
We note that $(\Lambda\rbn(t))_{t\le b_n t_0}$ is determined by the collection $\cN_{t_0}$ and is independent of $N\rw_m$. 
Thanks to Lemma~\ref{lem: Lambda-cond}, we can find two processes $(\mathtt{\Lambda}\rb(t))_{t\le b_nt_0}$ and $(\hat{\mathtt{\Lambda}}\rb(t))_{t\le b_nt_0}$ defined on  some probability space $(\Omega, \mathbb P, \mathcal F)$ so that $(\mathtt{\Lambda}\rb(t))_{t\le b_nt_0}$ has the same distribution as $(\Lambda\rbn(t))_{t\le b_n t_0}$ under $\bP_{n,m}$, $(\hat{\mathtt{\Lambda}}\rb(t))_{t\le b_nt_0}$ has the same distribution as $(\Lambda\rbn(t))_{t\le b_n t_0}$ under $\bP(\cdot\,|\,N\rb_n=n, N\rw_m=m)$, and
\[
\mathbb P\Big(\exists\, t\le b_n t_0: \mathtt{\Lambda}\rb(t)\ne \hat{\mathtt{\Lambda}}\rb(t)\Big)\to 0.
\]
By the tightness assumption, for each $\epsilon>0$, we can find $K=K(\epsilon, t_0)\in (0, \infty)$ satisfying
\[
\liminf_{n\to\infty}\mathbb P(\mathtt\Lambda\rb(b_n t_0)\le Kb_n) \ge 1-\epsilon.
\]
By enlarging the probability space $(\Omega, \mathbb P, \mathcal F)$, we apply Lemma~\ref{lem: Lambda-cond} again to find  two processes $(\mathtt{\Lambda}\rw(t))_{t\le Kb_n}$ and $(\hat{\mathtt{\Lambda}}\rw(t))_{t\le Kb_n}$ defined on it that have respective distributions $\nu_m\rw$ and $\hat\nu_m\rw$ and satisfy
\[
\mathbb P\Big(\exists\, t\le K b_n: \mathtt{\Lambda}\rw(t)\ne \hat{\mathtt{\Lambda}}\rw(t)\Big)\to 0.
\]
It follows that 
\[
\liminf_{n\to\infty}\mathbb P\Big(\mathtt\Lambda\rb(b_n t_0)\le Kb_n \text{ and } \mathtt{\Lambda}\rw\circ\mathtt{\Lambda}\rb(t)=\hat{\mathtt{\Lambda}}\rw\circ\hat{\mathtt{\Lambda}}\rb(t) \text{ for all } t\le b_n t_0\Big) \ge 1-\epsilon.
\]
This shows that the total variation distance between the law of $(\Lambda\rwm\circ\Lambda\rbn(t))_{t\le b_n t_0}$ under $\bP_{n,m}$ and the same process under $\bP_{n,m}(\cdot\,|\, N\rb_n=n, N\rw_m=m)$ tends to $0$. 
Since $Z^{n,m}_t = -t + \Lambda\rwm\circ\Lambda\rbn(t)$ and $H^{n,m}_t$ is a measurable function of $(Z^{n,m}_s)_{s\le t}$, the conclusion follows. 
\end{proof}

\subsubsection{Proof of Lemma~\ref{lem: girsa}}
\label{sec: pf-girsa}

Our proof of Lemma~\ref{lem: girsa} relies upon a Girsanov-type theorem for spectrally positive L\'evy processes, recalled in Appendix~\ref{sec: Z-prop}. 
 
\begin{proof}[Proof of Lemma~\ref{lem: girsa}]
Recall from~\eqref{def: Lbw} the compound Poisson processes $L\rb$ and $L\rw$ and from~\eqref{def: varphibw} their Laplace exponents $\varphi\rb$ and $\varphi\rw$. Recall $\Lambda\rbn(t)$ and $\Lambda\rwm(t)$ from~\eqref{def: Lam}, which are respective functions of the Poisson point measures $\Pi\rb$ and $\Pi\rw$. The exponential formula for Poisson point measures yields that 
\begin{align*}
&\bE_{n,m}[e^{-\lambda \Lambda\rwm(t)}]=\mathbb E\Big[\exp\Big(-\!\!\sum_{1\le j\le N\rw_m}\!\!\lambda Y_j\mathbf 1_{\{E\rw_j\le t\}}\Big)\Big]\\
&\qquad=\exp\Big\{-\int_{\R^2_+}\big(1-e^{-\lambda y\mathbf 1_{\{s\le t\}}}\big)\pi\rw(dy, ds) \Big\}\\
&\qquad = \exp\bigg\{-\sqrt{\frac{m}{n}}\int_0^t\int_{(0, \infty)}(1-e^{-\lambda y}) ye^{-ys/\sqrt{mn}}dF\rw(y) ds\bigg\}\\
&\qquad= \exp\bigg\{\sqrt{\frac{m}{n}}\int_0^t\int_{(0, \infty)}\Big\{(e^{-y(\lambda +s/\sqrt{mn})}-1)-(e^{-ys/\sqrt{mn}}-1)\Big\}ydF\rw(y)ds\bigg\}\\
&\qquad = \exp\Big\{\int_0^t \Big\{\varphi\rw\Big(\lambda+\frac{s}{\sqrt{mn}}\Big)-\varphi\rw\Big(\frac{s}{\sqrt{mn}}\Big)\Big\}ds \Big\}.
\end{align*}
We note that $\Lambda\rwm$ has independent increments. 
Comparing the previous calculations with~\eqref{id: lapQ}, we deduce the following absolute continuity relationship between $\Lambda\rwm$ and $L\rw$: for all $t\ge 0$ and measurable function $F:\mathbb D([0, t], \R_+)$, we have
\begin{equation}
\label{id: denLrw}
\bE_{n,m}\Big[F\Big(\big(\Lambda\rwm(s)\big)_{s\le t}\Big)\Big] = \bE_{n,m}\Big[F\Big(\big(L\rw_s\big)_{s\le t}\Big)\cdot E\rw_t\Big],
\end{equation}
where 
\[
E\rw_t = \exp\Big\{-\int_0^t \frac{s}{\sqrt{mn}}\,dL\rw_s-\int_0^t \varphi\rw\Big(\frac{s}{\sqrt{mn}}\Big)ds\Big\}.
\]
A similar identity holds for $\Lambda\rbn$ and $L\rb$:
\begin{equation}
\label{id: denLrb}
\bE_{n,m}\Big[F\Big(\big(\Lambda\rbn(s)\big)_{s\le t}\Big)\Big] = \bE_{n,m}\Big[F\Big(\big(L\rb_s\big)_{s\le t}\Big)\cdot E\rb_t\Big],
\end{equation}
with
\[
E\rb_t = \exp\Big\{-\int_0^t \frac{s}{\sqrt{mn}}\,dL\rb_s-\int_0^t \varphi\rb\Big(\frac{s}{\sqrt{mn}}\Big)ds\Big\}.
\]
Since $\Lambda\rbn(t)<\infty$ almost surely and its distribution is independent of $\Lambda\rwm$ and $L\rw$, we apply~\eqref{id: denLrw} and~\eqref{id: denLrb} to find that 
\begin{align*}
\bE_{n,m}[e^{-\lambda \Lambda\rwm\circ\Lambda\rbn(t)}]&=\int_{(0,\infty)}\bP_{n,m}\big(\Lambda\rbn(t)\in ds\big)\bE_{n,m}\big[e^{-\lambda \Lambda\rwm(s)}\big]\\
&=\int_{(0,\infty)}\bP_{n,m}\big(\Lambda\rbn(t)\in ds\big)\bE_{n,m}\big[e^{-\lambda L\rw_s}\cdot E\rw_s\big]\\
& = \int_{(0,\infty)}\bE_{n,m}\big[\mathbf 1_{\{L\rb_t\in ds\}}E\rb_t\big]\cdot \bE_{n,m}\big[e^{-\lambda L\rw_s}\cdot E\rw_s\big]\\
& = \bE_{n,m}\Big[e^{-\lambda L\rw\circ L\rb(t)}\cdot E\rw_{L\rb_t}\cdot E\rb_t\Big]\\
& =\bE_{n,m}\Big[e^{-\lambda L\rw\circ L\rb(t)}\cdot E^{n,m}_t\Big]. 
\end{align*}
Since $Z^{n,m}_t = -t +\Lambda\rwm\circ\Lambda\rbn(t)$ and $L^{n,m}_t = -t + L\rw\circ L\rb(t)$, we deduce that 
\[
\bE_{n,m}[e^{-\lambda Z^{n,m}_t}]=\bE_{n,m}[e^{-\lambda L^{n,m}_t}\cdot E^{n,m}_t], \quad \lambda\ge 0.
\]
The above shows an absolute continuity relationship between $Z^{n,m}_t$ and $L^{n,m}_t$. 
As both $(Z^{n,m}_t)_{t\ge 0}$ and $(L^{n,m}_t)_{t\ge 0}$ have independent increments, the distributions of both processes are determined by their marginal laws. The conclusion follows. 
\end{proof}

\subsubsection{Proof of Lemma~\ref{lem: cv-levy}}
\label{sec: pf-cv-levy}

We first recall from Whitt~\cite{Whitt} the following result on the convergence of the composed functions. 

\begin{lem}[Theorem 3.1 in~\cite{Whitt}]
\label{lem: cv-comp}
For each $n\in \N$, let $\mathbf x_n\in\mathbb D(\R_+, \R)$ and $\lambda_n\in \mathbb D(\R_+, \R_+)$. Suppose further that $t\mapsto \lambda_n(t)$ is non-decreasing. Let $\mathbf x\in \mathbb D(\R_+, \R)$ and $\lambda\in \mathbb C(\R_+, \R_+)$; moreover, suppose that $\lambda$ is  strictly increasing. Suppose further that $\mathbf x_n\to \mathbf x$ in $\mathbb D(\R_+, \R)$ and $\lambda_n\to\lambda$ in $\mathbb C(\R_+, \R)$. Then $\mathbf x_n\circ \lambda_n\to \mathbf x\circ\lambda$ in $\mathbb D(\R_+, \R)$. 
\end{lem}

\smallskip
\noindent
{\bf Part I: Proof of Lemma~\ref{lem: cv-levy} under the double finite third moments assumptions}

\smallskip
\noindent
{\bf Step 1: Convergence of $L^{n,m}$.}
It is well-known that the functional convergence of spectrally positive L\'evy processes is equivalent to the convergence of their Laplace exponents; see for instance Lemma A.3 combined with Theorem B.8 in~\cite{BrDuWa21}. 
Recall from~\eqref{def: Lbw} that $L\rb$ is a compound Poisson process of Laplace exponent $\varphi\rb$ under $\bP_{n,m}$. 
Let us introduce 
\begin{equation}
\label{def: hatLbw}
\hat L\rb_t := L\rb_t-\bE_{n,m}[L\rb_t]=L\rb_t-\sqrt{\tfrac{n}{m}}\,\sigma\rb_2 t, \quad t\ge 0,
\end{equation}
so that for all $\lambda\ge 0$, 
\[
\bE_{n,m}[e^{-\lambda \hat L\rb_1}]=e^{\hat\varphi\rb(\lambda)}
\]
with
\begin{equation}
\label{def: hatphi'}
\hat\varphi\rb(\lambda) =\varphi\rb(\lambda)+\sqrt{\tfrac{n}{m}}\sigma\rb_2\lambda = \sqrt{\tfrac{n}{m}}\int_{(0,\infty)} (e^{-\lambda x}-1+\lambda x) x dF\rb(x). 
\end{equation}
Similarly, let us set
\begin{equation}
\label{def: hatL}
\hat L\rw_t := L\rw_t-\bE_{n,m}[L\rw_t]=L\rw_t-\sqrt{\tfrac{m}{n}}\,\sigma\rw_2 t, \quad t\ge 0,
\end{equation}
and 
\begin{equation}
\label{def: hatphi}
\hat\varphi\rw(\lambda) =\log\bE_{n,m}[e^{-\lambda \hat L\rw_1}] = \sqrt{\tfrac{m}{n}}\int_{(0,\infty)} (e^{-\lambda x}-1+\lambda x) x dF\rw(x). 
\end{equation}
Under the assumptions~\eqref{hyp: b-third},~\eqref{hyp: w-third} and~\eqref{hyp: mod-clst}, Lemma~\ref{lem-cv: psi-br} asserts that 
\begin{equation}
\label{cv: varphi}
n^{\frac23}\hat\varphi\rb(n^{-\frac13}\lambda) \xrightarrow{n\to\infty} \tfrac12\theta^{-\frac12} \sigma\rb_3\, \lambda^2, \quad
n^{\frac23}\hat\varphi\rw(n^{-\frac13}\lambda) \xrightarrow{n\to\infty} \tfrac12\theta^{\frac12} \sigma\rw_3\, \lambda^2.
\end{equation}
Combined with the independence between $(\hat L\rb_t)_{t\ge 0}$ and $(\hat L\rw_t)_{t\ge 0}$, this implies:
\begin{equation}
\label{cv: L-bw-br}
\Big\{n^{-\frac13}\hat L\rb_{n^{2/3}t}, n^{-\frac13}\hat L\rw_{n^{2/3}t}: t\ge 0\Big\} \Longrightarrow \big\{\cL\rb_t, \cL\rw_t: t\ge 0\big\} \quad \text{in }\mathbb D(\R_+, \R^2),
\end{equation}
where $(\cL\rb, \cL\rw)$ is a pair of independent Brownian motion with respective Laplace exponents $\Psi\rb, \Psi\rw$ given by 
\begin{equation}
\label{def: cLrb-br}
\Psi\rb(\lambda) = \tfrac12 \theta^{-\frac12}\sigma\rb_3\lambda^2, \quad \Psi\rw(\lambda) =  \tfrac12 \theta^{\frac12}\sigma\rw_3\lambda^2, \quad \lambda\ge 0.  
\end{equation}
We observe that the previous convergence of $\hat L\rb_t$ also implies that for any $t_0\ge 0$, 
\begin{equation}
\label{cv: L-b}
\sup_{t\le t_0}n^{-\frac23}\hat L\rb_{n^{2/3}t} \to 0 \text{ in prob., so that } 
\sup_{t\le t_0} \Big|n^{-\frac23} L\rb_{n^{2/3}t} - \theta^{-\frac12}\sigma\rb_2 t\Big| \to 0 \text{ in prob}.
\end{equation}
Recall from~\eqref{def: Lnm} the definition of $L^{n,m}_t$. 
Using the criticality assumption~\eqref{hyp: crit}, we can write
\begin{align*}
L^{n,m}_t &= -t + L\rw\circ L\rb(t) = -t + \hat L\rw\circ L\rb(t) + \sqrt{\tfrac{m}{n}}\sigma\rw_2 L\rb_t \\
&= \hat L\rw\circ L\rb(t) + \sqrt{\tfrac{m}{n}}\sigma\rw_2 \hat L\rb_t - \Big(1-\sqrt{\tfrac{m}{n}}\sqrt{\tfrac{n}{m}}\,\sigma\rw_2\sigma\rb_2\Big)t\\
& = \hat L\rw\circ L\rb(t) + \sqrt{\tfrac{m}{n}}\sigma\rw_2 \hat L\rb_t. 
\end{align*}
It follows that
\[
n^{-\frac13}L^{n,m}_{n^{2/3}t} = n^{-\frac13}\hat L\rw\Big(n^{\frac23}\cdot n^{-\frac23} L\rb(n^{\frac23}t)\Big) + \sqrt{\tfrac{m}{n}}\sigma\rw_2 n^{-\frac13}\hat L\rb_{n^{2/3}t}.
\]
Using the joint convergence in~\eqref{cv: L-bw-br} and the convergence in probability from~\eqref{cv: L-b}, we apply Lemma~\ref{lem: cv-comp} to find that 
\begin{equation}
\label{cv: Lnm-br}
\Big\{n^{-\frac13}L^{n,m}_{n^{2/3}t}: t\ge 0\Big\} \Longrightarrow \Big\{\cL\rw_{\sigma\rb_2 t/\sqrt\theta}+\sqrt\theta\sigma\rw_2\cL\rb_t: t\ge 0\Big\} \quad \text{in }\mathbb D(\R_+, \R),
\end{equation}
jointly with the convergence in~\eqref{cv: L-bw-br}. Using the scaling property of Brownian motions, we find that the right-hand side in~\eqref{cv: Lnm-br} is a Brownian motion with the quadratic variation $\sigma\rw_3\sigma\rb_2+\sqrt\theta(\sigma\rw_2)^2\sigma\rb_3$, which has the same distribution as $(\cL^{(1)}_t)_{t\ge 0}$ in~\eqref{def: L}. This confirms the desired convergence from $L^{n,m}$.

\smallskip
\noindent
{\bf Step 2: Convergence of $H^{n,m}$.}
We intend to apply Theorem~\ref{thm: DuLG} with $a_n=n^{1/3}$ and $b_n=n^{2/3}$ to $X^{(n)}=L^{n,m}$ and $X=\cL^{(1)}$. 
To check the assumptions of the theorem, we recall that $(L^{n,m}_t+t)_{t\ge 0}$ is a compound Poisson process. 
Denoting $\pi_n$ for the L\'evy measure of this compound Poisson process and writing $\varphi$ for the Laplace exponent of $L^{n,m}$, we have 
\begin{equation}
\label{eq: vph}
\varphi(\lambda)-\lambda = \int_{(0,\infty)} (e^{-\lambda x}-1) \pi_n(dx) = -\int_0^{\lambda} \int_{(0, \infty)}e^{-ux}x\pi_n(dx),
\end{equation}
 where we have used Fubini's theorem in the last identity.
Meanwhile, the subordination of L\'evy processes implies that 
\begin{equation}
\label{id: subord}
\varphi(\lambda) -\lambda =  \varphi\rb \Big(- \varphi\rw(\lambda)\Big), \quad \lambda\ge 0,
\end{equation}
where $\varphi\rw$ and $\varphi\rb$ are  the respective Laplace exponents of $L\rw$ and $L\rb$. 
Taking $\lambda$ to infinity in~\eqref{eq: vph}, we deduce from the monotone convergence theorem that 
$q_n:=\pi_n(\R_+) = -\lim_{\lambda\to \infty}(\varphi(\lambda)-\lambda)$. 
In view of~\eqref{id: subord} and the definition~\eqref{def: varphibw} of $\varphi\rw$, we have
\[
q_n = - \lim_{\lambda\to \infty}\varphi\rb\Big(-\varphi\rw(\lambda) \Big)= -\varphi\rb\Big(\sqrt{\tfrac{m}{n}} \sigma\rw_1\Big).
\]
It follows that
\begin{equation}
\label{lim: qn}
\lim_{n\to\infty} q_n = \theta^{-\frac12} \int_{(0,\infty)}\big(1-e^{-\sqrt\theta \sigma\rw_1 x}\big)\, xdF\rb(x)=: q_0\in (0,\infty).
\end{equation}
Similarly for the uniform integrability of $(\pi_n)_{n\in \N}$, we differentiate~\eqref{eq: vph} and~\eqref{id: subord} to find that
\[
\int_{(0, \infty)}e^{-\lambda x}x\pi_n(dx) = \big(\varphi\rb\big)'\Big(-\varphi\rw(\lambda)\Big) \cdot (\varphi\rw)'(\lambda).
\]
We deduce that  $\phi(\lambda):=\lim_{n\to\infty}\int e^{-\lambda x}x\pi_n(dx)$ exists and is continuous at 0; hence, $(\pi_n)_{n\in \N}$ is uniformly integrable.
On the other hand, the weak convergence of $n^{-1/3}L^{n,m}_{n^{2/3}}$ towards $\cL^{(1)}_1$ is implied by~\eqref{cv: Lnm-br}. It remains to check~\eqref{condH}. We shall prove that the condition~\eqref{cond: H} is satisfied in this case. 
Recall $g_n$ from~\eqref{def: gn} and $\Psi_n$ from~\eqref{def: Psin}. Comparing them with $\varphi$ in~\eqref{eq: vph}, we find that
\[
g_n(s) -1 = \frac{1}{q_n}\int_{(0,\infty)} (e^{-(1-s)q_nx}-1)\pi_n(dx) = \frac{1}{q_n}\Big(\varphi\big((1-s)q_n\big)-(1-s)q_n\Big),
\]
so that 
\begin{align*}
\Psi_n(u) & = n^{\frac23}q_n^{-1}\varphi(un^{-\frac13}q_n) =  n^{\frac23}q_n^{-1}\Big( \varphi\rb \big(- \varphi\rw(un^{-\frac13}q_n)\big)+un^{-\frac13}q_n\Big)\\
& =n^{\frac23}q_n^{-1}\Big( \hat\varphi\rb \big(- \varphi\rw(un^{-\frac13}q_n)\big)+\sqrt{\tfrac{n}{m}}\sigma\rb_2\hat\varphi\rw(un^{-\frac13}q_n)\Big). 
\end{align*}
where we have used~\eqref{id: subord} in the second identity, the definitions~\eqref{def: hatphi'},~\eqref{def: hatphi} of $\hat\varphi\rb$, $\hat\varphi\rw$ and the criticality assumption~\eqref{hyp: crit} in the third. 
Using the inequality $e^{-x}-1+x\ge \frac12 x^2e^{-x}$ for all $x\ge 0$ and the convergence of $q_n\to q_0$ as seen in~\eqref{lim: qn}, we find that for all $u\in [0, n^{1/3}]$ and $n$ sufficiently large, 
\begin{align*}
\hat\varphi\rw(un^{-\frac13}q_n)&\ge \tfrac12\sqrt{\tfrac{m}{n}}(q_n u n^{-\frac13})^2\int_{(0, \infty)} e^{-un^{-1/3}q_n x}x^3dF\rw(x)\\
&\ge \tfrac12\sqrt{\tfrac{m}{n}}(q_n u n^{-\frac13})^2\int_{(0, \infty)} e^{-q_n x}x^3dF\rw(x) \ge \tfrac14 \sqrt{\tfrac{m}{n}}\, q_n^2\, \eta\, u^2 n^{-\frac23},
\end{align*}
with  $\eta:=\int_{(0,\infty)}x^3 e^{-q_0x}dF\rw(x)\in (0,\infty)$. Combined with the fact that $\hat{\varphi}\rb$ only takes non negative values, this yields that for sufficiently large $n$ and all $u\in [0, n^{1/3}]$: 
\[
\Psi_n(u)\ge n^{\frac23}q_n^{-1}\sqrt{\tfrac{n}{m}}\sigma\rb_2\hat\varphi\rw(un^{-\frac13}q_n)\ge \tfrac14\sigma\rb_2 q_n \eta\, u^2, 
\]
from which~\eqref{cond: H} follows. In consequence, Theorem~\ref{thm: DuLG} applies to yield the convergence in~\eqref{cv: LH1}. 

\smallskip
\noindent
{\bf Part II: Proof of Lemma~\ref{lem: cv-levy} under the one dominant heavy tail assumptions.}
The proof here follows the same steps as in the previous one. We only point out the main differences and omit the details. For the convergence of $L^{n,m}$, introduce $\hat L\rb$ and $\hat L\rw$ as in~\eqref{def: hatLbw} and~\eqref{def: hatL}, with their respective Laplace exponents $\hat\varphi\rb$ and $\hat\varphi\rw$. 
Under~\eqref{hyp: b-power} and~\eqref{hyp: mod-clst}, Lemma~\ref{lem-cv: psi-power} implies that
\begin{equation}
\label{eqcv: Lrb-power}
\Big\{n^{-\frac{1}{\alpha+1}}\hat L\rb_{ n^{\alpha/(\alpha+1)}t}: t\ge 0\Big\} \Longrightarrow \big\{\cL\rb_t: t\ge 0\big\}\quad \text{in }\mathbb D(\R_+, \R),
\end{equation}
where in this case $\cL\rb$ is a spectrally positive $\alpha$-stable process satisfying  
\begin{equation}
\label{def: Lrb-stable}
\mathbb E[e^{-\lambda \cL\rb_1}]=e^{\Psi\rb(\lambda)}, \quad \text{with } \Psi\rb(\lambda) = \theta^{-\frac12}C\rb C(\alpha)\lambda^{\alpha}, \quad \lambda\ge 0. 
\end{equation}
In above, the constant $C(\alpha)=(\alpha+1)\Gamma(2-\alpha)/\alpha(\alpha-1)$. 
On the other hand, assuming in the first instance that~\eqref{hyp: w-third} holds. Then  Lemma~\ref{lem-cv: psi-br} says that for all $\lambda\ge 0$,  
\begin{equation}
\label{cv: Lw-one}
n^{\frac{\alpha}{\alpha+1}}\hat{\varphi}\rw(n^{-\frac{1}{\alpha+1}}\lambda) \to 0, \text{ so that } n^{-\frac{1}{\alpha+1}}\sup_{t\le t_0}|\hat L\rw_{ n^{\alpha/(\alpha+1)}t}| \to 0 \quad\text{in prob.~for all }t_0\ge 0.
\end{equation}
If, instead,~\eqref{hyp: w-power} holds for some $\alpha'\in (\alpha, 2)$, replacing Lemma~\ref{lem-cv: psi-br} with Lemma~\ref{lem-cv: psi-power} yields the same conclusion. Meanwhile, the same arguments leading to~\eqref{cv: L-b} yield in the current case: 
\begin{equation}
\label{cv: bk}
\sup_{t\le t_0} \Big|n^{-\frac{\alpha}{\alpha+1}}L\rb_{ n^{\alpha/(\alpha+1)}t} -\theta^{-\frac12}\sigma\rb_2 t\Big|\to 0 \quad\text{in probability for all }t_0\ge 0.
\end{equation}
Arguing as in the previous part, we deduce from the critical assumption~\eqref{hyp: crit} that
\[
n^{-\frac{1}{\alpha+1}}L^{n,m}_{n^{\alpha/(\alpha+1)}t} = n^{-\frac{1}{\alpha+1}}\hat L\rw\Big(n^{\frac{\alpha}{\alpha+1}}\cdot n^{-\frac{\alpha}{\alpha+1}}L\rb(n^{\frac{\alpha}{\alpha+1}}t)\Big) + \sqrt{\tfrac{m}{n}} \sigma\rw_2 n^{-\frac{1}{\alpha+1}}\hat L\rb_{ n^{\alpha/(\alpha+1)}t}.
\]
Lemma~\ref{lem: cv-comp} combined with~\eqref{cv: Lw-one},~\eqref{cv: bk} and~\eqref{eqcv: Lrb-power} implies that
\begin{equation}
\label{cv: Lnm-one}
\Big\{n^{-\frac{1}{\alpha+1}}L^{n,m}_{ n^{\alpha/(\alpha+1)}t}: t\ge 0\Big\} \Longrightarrow \Big\{\sqrt\theta\sigma\rw_2\cL\rb_t: t\ge 0\Big\} \quad \text{in }\mathbb D(\R_+, \R),
\end{equation}
jointly with the convergence in~\eqref{eqcv: Lrb-power}. We readily check that the limit process above is distributed as $\cL^{(2)}$ in~\eqref{def: L}. For the convergence of the height process $H^{n,m}$, we rely on Theorem~\ref{thm: DuLG} with $a_n=n^{1/(\alpha+1)}$ and $b_n=n^{\alpha/(\alpha+1)}$. The beginning of the proof is identical to the previous part. In particular,~\eqref{lim: qn} still holds, with $q_n$ standing for the total mass of the L\'evy measure of $L^{n,m}$. To check the condition~\eqref{cond: H} in the current case, let us show that we can find some $C\in (0,\infty)$ so that for all $u\in [0, n^{1/(\alpha+1)}]$ and sufficiently large $n$,
\[
\Psi_n(u) = n^{\frac{\alpha}{\alpha+1}}q_n^{-1} \Big( \hat\varphi\rb \big(- \varphi\rw(un^{-\frac{1}{\alpha+1}}q_n)\big)+\sqrt{\tfrac{n}{m}}\sigma\rb_2\hat\varphi\rw(un^{-\frac{1}{\alpha+1}}q_n)\Big)\ge Cu^{\alpha}.
\]
Indeed, from the inequality $1-e^{-x}\ge xe^{-x}$ for $x\ge 0$, we obtain that for all $u\in [0, n^{1/(\alpha+1)}]$ and sufficiently large $n$, 
\[
-\varphi\rw(un^{-\frac{1}{\alpha+1}}q_n) = \sqrt{\tfrac{m}{n}}\int_{(0,\infty)} (1-e^{-u n^{-\frac{1}{\alpha+1}}q_nx}) xdF\rw(x)\ge  \tfrac12\sqrt\theta q_0\,\eta 'un^{-\frac{1}{\alpha+1}},
\]
with $\eta':=\int_{(0,\infty)}x^2 e^{-q_0x}dF\rw(x)\in (0,\infty)$. 
Applying the monotonicity of $\hat\varphi\rb$ and then Lemma~\ref{lem: psi-bound}, we can find some positive constant $c$ so that for all $u\in [0,  n^{1/(\alpha+1)}]$ and sufficiently large $n$, 
 \[
\hat\varphi\rb \big(- \varphi\rw(un^{-\frac{1}{\alpha+1}}q_n)\big)\ge \hat\varphi\rb\big(\tfrac12 \sqrt\theta q_0\eta' un^{-\frac{1}{\alpha+1}}\big) \ge c u^{\alpha}n^{-\frac{\alpha}{\alpha+1}},
\]
from which it follows that
\[
\Psi_n(u) \ge n^{\frac{\alpha}{\alpha+1}}q_n^{-1}\hat\varphi\rb \big(- \varphi\rw(un^{-\frac{1}{\alpha+1}}q_n)\big)\ge  c q_n^{-1} u^{\alpha}. 
\]
Hence, condition~\eqref{cond: H} holds in this case.  The joint convergence in~\eqref{cv: LH2} then follows as a result of Theorem~\ref{thm: DuLG}. 

\smallskip
\noindent
{\bf Part III: Proof of Lemma~\ref{lem: cv-levy} under the matched heavy tails assumptions.} This is very similar to the previous two cases. The main difference lies in the counterpart of~\eqref{cv: L-bw-br}: under the current assumptions, we have instead
\begin{equation}
\label{cv: L-bw-matched}
\Big\{n^{-\frac{1}{\alpha+1}}\hat L\rb_{n^{\alpha/(\alpha+1)}t}, n^{-\frac{1}{\alpha+1}}\hat L\rw_{n^{\alpha/(\alpha+1)}t}: t\ge 0\Big\} \Longrightarrow \big\{\cL\rb_t, \cL\rw_t: t\ge 0\big\} \quad \text{in }\mathbb D(\R_+, \R^2),
\end{equation}
where $(\cL\rb, \cL\rw)$ is a pair of independent $\alpha$-stable processes with respective Laplace exponents $\Psi\rb, \Psi\rw$ given by 
\begin{equation}
\label{def: cLbw-stable}
\Psi\rb(\lambda) = \theta^{-\frac12}C\rb C(\alpha)\lambda^{\alpha}, \quad\Psi\rw(\lambda) = \theta^{\frac12}C\rw C(\alpha)\lambda^{\alpha}, \quad \lambda\ge 0.
\end{equation}
As a result, we have
\begin{equation}
\label{eqcv: L-3}
\Big\{n^{-\frac{1}{\alpha+1}}L^{n,m}_{n^{\alpha/(\alpha+1)}t}: t\ge 0\Big\} \Longrightarrow \Big\{\cL\rw_{\sigma\rb_2 t/\sqrt\theta}+\sqrt\theta\sigma\rw_2\cL\rb_t: t\ge 0\Big\} \quad \text{in }\mathbb D(\R_+, \R),
\end{equation}
jointly with the convergence in~\eqref{cv: L-bw-matched}. We then observe that the limit process in~\eqref{eqcv: L-3} has the same distribution as $\cL^{(3)}$ introduced in~\eqref{def: L}. The condition~\eqref{cond: H} can be checked similarly as in Part II. 

\medskip
Before ending this subsection, let us point out that we have in fact shown that under the  respective sets of assumptions in Lemma~\ref{lem: cv-levy}, 
\begin{equation}
\label{eqcv: Lnm-uni}
\Big\{\tfrac{1}{a_n}\hat L\rb_{b_n t}, \tfrac{1}{a_n}\hat L\rw_{b_n t}, \tfrac{1}{a_n}L^{n,m}_{b_n t}: t\ge 0\Big\} \Longrightarrow \Big\{\cL\rb_t, \cL\rw_t, \cL\rw_{\sigma\rb_2 t/\sqrt\theta}+\sqrt\theta\sigma\rw_2\cL\rb_t: t\ge 0\Big\}
\end{equation}
in $\mathbb D(\R_+, \R^2)\times \mathbb D(\R_+, \R)$, as well as
\begin{equation}
\label{eqcv: Lb-uni}
\sup_{t\le t_0} \Big|\tfrac{1}{b_n}L\rb_{b_n t} -\theta^{-\frac12}\sigma\rb_2 t\Big|\to 0 \quad\text{in probability for all }t_0\ge 0.
\end{equation}
where 
\begin{itemize}
\item
under the double finite third moments assumptions,  $a_n=n^{1/3}, b_n=n^{2/3}$ and $(\cL\rb, \cL\rw)$ is a pair of independent Brownian motions introduced in~\eqref{def: cLrb-br};
\item
under the one dominant heavy tail assumptions, $a_n = n^{1/(\alpha+1)}, b_n = n^{\alpha/(\alpha+1)}$, $\cL\rw$ is the null process while $\cL\rb$ is an $\alpha$-stable process that satisfies~\eqref{def: Lrb-stable};
\item
under the matched heavy tails assumptions, $a_n = n^{1/(\alpha+1)}, b_n = n^{\alpha/(\alpha+1)}$, and $(\cL\rb, \cL\rw)$ is a pair of independent $\alpha$-stable processes satisfying~\eqref{def: cLbw-stable};
\end{itemize}
Under the relevant assumptions, the convergence in~\eqref{eqcv: Lnm-uni} corresponds to respectively~\eqref{cv: L-bw-br} and~\eqref{cv: Lnm-br} from Part I,~\eqref{eqcv: Lrb-power},~\eqref{cv: Lw-one},~\eqref{cv: Lnm-one} from Part II, and~\eqref{cv: L-bw-matched},~\eqref{eqcv: L-3} from Part III. The convergence in~\eqref{eqcv: Lb-uni} is shown as~\eqref{cv: L-b} in Part I, as~\eqref{cv: bk} in Part II, and can be shown in the same way under the matched heavy tails assumptions.

\subsubsection{Proof of Lemma~\ref{lem: cv-dens}}
\label{sec: pf-cv-dens}

We will need the following lemma, whose proof is elementary and thus omitted. 

\begin{lem}
\label{lem: int}
Let $\bx_n\in \mathbb D(\R_+, \R)$ for each $n\in\N$, and $\bx\in \mathbb D(\R_+, \R)$. Assume that $\bx_n\to \bx$ in $\mathbb D(\R_+, \R)$ as $n\to\infty$. Then for any positive sequence $a_n\to t\in (0, \infty)$, we have 
\[
\int_0^{a_n} \bx_n(s) ds \xrightarrow{n\to\infty}  \int_0^t \bx(s)ds=\int_0^t \bx(s-)ds.
\]
\end{lem}

\begin{proof}[Proof of Lemma~\ref{lem: cv-dens}]
We adopt the notation from~\eqref{eqcv: Lnm-uni}. We note that the convergences of $(\hat L\rb, \hat L\rw)$ therein imply the following convergences of their Laplace exponents: for all $\lambda\ge 0$, 
\begin{equation}
\label{eqcv: lap}
b_n\hat\varphi\rb(\lambda/a_n) \to \Psi\rb(\lambda), \quad b_n\hat\varphi\rw(\lambda/a_n)\to \Psi\rw(\lambda), 
\end{equation}
where 
under the double finite third moments assumptions, $\Psi\rb, \Psi\rw$ are given in~\eqref{def: cLrb-br}; 
under the one dominant heavy tail assumptions, $\Psi\rb$ is defined in~\eqref{def: Lrb-stable} and $\Psi\rw\equiv 0$; 
under the matched heavy tails assumptions, $\Psi\rb, \Psi\rw$ are introduced in~\eqref{def: cLbw-stable}. In all three cases, we have $a_n=n^{1/(\alpha+1)}$ and $b_n=n^{\alpha/(\alpha+1)}$ with the understanding that $\alpha=2$ for the double finite third moments assumptions. Let us denote
\[
\rho = \sigma\rb_2/\sqrt\theta, \quad\text{so that}\quad \rho^{-1} = \sqrt\theta\sigma\rw_2,
\]
thanks to the criticality assumption~\eqref{hyp: crit}. Recall that $\cL^{(i)}$ are the respective limit of $L^{n,m}$ under the relevant sets of assumptions and $\Psi_i$ is its Laplace exponent, for $i\in \{1, 2, 3\}$. In view of~\eqref{eqcv: Lnm-uni}, we can assume from now on that 
\begin{equation}
\label{def': cL}
\cL^{(i)}_t = \cL\rw_{\rho t} + \rho^{-1} \cL\rb_t \quad \text{and} \quad \Psi_i(\lambda) = \rho\Psi\rw(\lambda) + \Psi\rb(\lambda/\rho).
\end{equation}
Thanks to Skorokhod's Representation Theorem, we are able to assume that the convergences in~\eqref{eqcv: Lnm-uni} and~\eqref{eqcv: Lb-uni} both take place {\it almost surely}. Recall from~\eqref{def: Enm} the density $E^{n,m}_t$, which we split as follows: 
\[
-\log E^{n,m}_t =D\rw_t + D\rb_t, 
\]
where
{\small\[
D\rw_t = \int_0^{L\rb_t} \frac{s}{\sqrt{mn}} \,dL\rw_s + \int_0^{L\rb_t} \varphi\rw\Big(\frac{s}{\sqrt{mn}}\Big) ds, \ 
D\rb_t =\int_0^{t} \frac{s}{\sqrt{mn}}\, dL\rb_s + \int_0^{t} \varphi\rb\Big(\frac{s}{\sqrt{mn}}\Big) ds. 
\]}
Recall from~\eqref{def: hatLbw} and~\eqref{def: hatL} the compensated processes $\hat L\rw, \hat L\rb$,  with respective Laplace exponents $\hat\varphi\rw, \hat\varphi\rb$. We have
\begin{align*}
D\rw_t &= \int_0^{L\rb_t} \frac{s}{\sqrt{mn}} \,d\Big(\hat L\rw_s+\sqrt{\tfrac{m}{n}}\,\sigma\rw_2 s\Big) + \int_0^{L\rb_t} \Big\{\hat\varphi\rw\Big(\frac{s}{\sqrt{mn}}\Big) - \sqrt{\tfrac{m}{n}}\,\sigma\rw_2\frac{s}{\sqrt{mn}}\Big\}ds\\
& = \int_0^{L\rb_t} \frac{s}{\sqrt{mn}} \,d\hat L\rw_s + \int_0^{L\rb_t} \hat\varphi\rw\Big(\frac{s}{\sqrt{mn}}\Big) ds\\
&=\frac{L\rb_t}{\sqrt{mn}}\hat L\rw\circ L\rb(t)-\frac{1}{\sqrt{mn}}\int_0^{L\rb_t}\hat L\rw_sds +\int_0^{L\rb_t} \hat\varphi\rw\Big(\frac{s}{\sqrt{mn}}\Big)ds\\
& = \frac{1}{\sqrt{mn}}\int_0^{L\rb_t} \big(\hat L\rw\circ L\rb(t)-\hat L\rw_s\big) ds + \int_0^{L\rb_t} \hat\varphi\rw\Big(\frac{s}{\sqrt{mn}}\Big)ds,
\end{align*}
where we have applied an  integration by parts for stochastic integration in the penultimate line (\cite{Protter}, Section II.6, Corollary 2). 
A change of variables $s=b_nu$ then yields
\[
D\rw_{b_n t} = \frac{b_n}{\sqrt{mn}}\int_0^{L\rb_{b_nt}/b_n} \Big(\hat L\rw\circ L\rb(b_nt)-\hat L\rw_{b_n u}\Big) du + \int_0^{L\rb_{b_nt}/b_n} b_n\hat\varphi\rw\Big(\frac{b_n u}{\sqrt{mn}}\Big)du.
\]
We note that $a_n b_n =n$. 
In view of Lemmas~\ref{lem: int} and~\ref{lem: cv-comp}, the convergences in~\eqref{eqcv: Lnm-uni},~\eqref{eqcv: Lb-uni} along with~\eqref{eqcv: lap} imply that 
\[
D\rw_{b_nt} \xrightarrow{\text{a.s.}} \frac{1}{\sqrt\theta} \int_0^{ \rho t} \big(\cL\rw_{\rho t} - \cL\rw_{u}\big) du + \int_0^{\rho t} \Psi\rw\Big(\frac{u}{\sqrt\theta}\Big) du. 
\]
We then apply a change of variables, followed by  an integration by parts to find that the previous limit is equal to 
\[
D\rw_{b_n t}\xrightarrow{\text{a.s.}}\frac{1}{\sqrt\theta}  \int_0^{t} \rho s\, d\cL\rw_{\rho s} + \int_0^{t} \rho \Psi\rw\Big(\frac{\rho s}{\sqrt\theta}\Big) ds.
\]
By a similar argument, we can show that
\begin{equation}
\label{cv: Db}
D\rb_{b_n t}\xrightarrow{\text{a.s.}}\frac{1}{\sqrt\theta}\int_0^t s\,d\cL\rb_s +\int_0^t \Psi\rb\Big(\frac{s}{\sqrt\theta}\Big)\,ds. 
\end{equation}
Combining the two and noting~\eqref{def': cL}, we conclude that 
\[
-\log E^{n,m}_{b_nt} = D\rw_{b_nt}+D\rb_{b_nt} \xrightarrow{\text{a.s.}}  \int_0^t \frac{\rho s}{\sqrt\theta}\,d\cL^{(i)}_s +  \int_0^t \Psi_i\Big(\frac{\rho s}{\sqrt\theta}\Big) ds = -\log\cE^{(i)}_t,
\]
jointly with the convergence in~\eqref{eqcv: Lnm-uni}. 
As $E^{n,m}_t$ takes non negative values and $\bE_{n,m}[E^{n,m}_t]=1=\mathbb E[\cE^{(1)}_t]$, which is implied by~\eqref{id: gir-dis} and~\eqref{id: girsanov}, the previous convergence also implies that $E^{n,m}_t, n\ge 1$, is uniformly integrable; see for instance Lemma 3.11 in~\cite{Kal}. 

For the convergence of $\Lambda\rbn$, recall $E\rb_t$ from~\eqref{id: denLrb} and note that $E\rb_t = \exp(-D\rb_t)$. The convergence in~\eqref{cv: Db}, combined with the fact that $E\rb_t$ is a positive random variable with unit mean, implies that $E\rb_{b_nt}, n\ge 1$, is uniformly integrable. 
We then deduce from~\eqref{id: denLrb} and~\eqref{eqcv: Lb-uni} that for any $\epsilon>0$, 
\[
\bP_{n,m}\Big(\sup_{s\le t}\Big|\frac{1}{b_n}\Lambda\rbn(b_n s)-\rho s\Big|\ge \epsilon\Big) =\bE_{n,m}\Big[\mathbf 1_{\{\sup_{s\le t}|L\rb_{b_n s}/b_n-\rho s|\ge \epsilon\}}\cdot E\rb_{b_n t}\Big] \xrightarrow{n\to\infty} 0,
\]
which completes the proof. 
\end{proof}

\subsection{Proof for the convergence of surplus edges}
\label{sec: cv-surplus}

\subsubsection{Preliminaries}

We give here the proof of Lemmas~\ref{lem: dis}-\ref{lem: dist-sur}. 

\begin{proof}[Proof of Lemma~\ref{lem: dis}]
\label{pf: dis}
The first statement is clear since for $e\in\cE'_{n,m}$, its counterpart $e'=(V_{k(e)}, V_{k'(e)})$ is in the same connected component as $e$.  
For the upper bound of $\mathrm{dis}_k$, recall from Section~\ref{sec: ghp} the notion of path and $(\Pi, \epsilon)$-modified length. Suppose that $\cE_k=\{e_l=(b_{i_l}, w_{i_l}): 1\le l\le s_{n,k}\}$ is the subset of $\cE_{n,m}$ containing those edges with ends in $C^{\bx}_{n, (k)}$. Write $\cE'_k = \{e'_l=(V_{k(e_l)}, V_{k'(e_l)}): 1\le l\le s_{n, k}\}$ its counterpart in $\cE'_{n,m}$. Let $x, y\in C^{\bx}_{n,(k)}$ and let $(f_1, f_2, \dots, f_p)$ be the shortest path in $B''$ between $x$ and $y$, where $f_i = (x_i, y_i)$, $1\le i\le p$, are edges of $B''$. Denote $\overleftarrow{f_i}=(y_i, x_i)$. 
Since $\cE'_{k}$ contains at most $s_{n, k}$ elements, there are at most $s_{n, k}$ elements among $(f_1, f_2, \dots, f_p)$ for which either $f_i$ or $\overleftarrow{f_i}$ belongs to $\cE'_k$, as otherwise the path can be shortened by removing cycles. For those $f_i$, their lengths in $d''$ are $1$, while their lengths in $\dgr$ is 2. The remaining $f_i$ are edges that are present in both $B_{n,m}$ and $B''$. It follows that 
\[
d''(x, y) = \sum_{1\le i\le p}d''(x_i, y_i) \ge\sum_{1\le i\le p}\dgr(x_i, y_i) - s_{n,k}\ge \dgr(x, y)-s_{n,k}, 
\]
for any $x, y\in C^{\bx}_{n, (k)}$. A similar argument for $\dgr$ leads to the bound in the other direction. 
The conclusion then follows. 
\end{proof}

\begin{proof}[Proof of Lemma~\ref{lem: sigma-p}]
\label{pf: sigma-p}
For the first statement, note that we can write 
\[
\Sigma^{n,m}(t)=\sum_{1\le k\le N\rb_n}\Sigma_k(t), \quad\text{where}\quad \Sigma_k(t)=X_k\mathbf 1_{\{k\in \mathfrak D_t\}}+\frac{\Delta_k(t)}{\Delta\bi_k}X_k\mathbf 1_{\{k\in \mathfrak W_t\}}.
\]
Observe that $t\mapsto \Sigma_k(t)$ is non-decreasing. Indeed, denote by $\tau_k$ the moment when Client $k$ leaves the queue, with the convention $\tau_k=E\rb_k$ if $\Delta\bi_k=0$; then $\Sigma_k(t)=0$ for all $t<E\rb_k$, $\Sigma_k(t)=X_k$ for all $t\ge \tau_k$, and $\Sigma_k$ is  non-decreasing between $E\rb_k$ and $\tau_k$. Moreover, if the queue is not empty at time $t$, then the server is serving a unique client, say $k$. 
If $\Delta\bi_k=0$, then we have $t=E\rb_k$ and $\Sigma_k$ jumps upwards at $t$. If $\Delta\bi_k>0$, we can find a small neighbourhood of $t$ contained in $J_k$, so that $\Sigma_k$ is strictly increasing on this neighbourhood. 
The first statement then follows. 
For the second, we first consider the case where $\Delta\bi_k>0$. In that case, $J_k$ has Lebesgue measure $\Delta\bi_k$, since in the bipartite queue Client $k$ requests $\Delta\bi_k$ service time. We then observe that when Client $k$ is being served, the slopes of both $\Sigma_k$ and $\Sigma^{n,m}$ are precisely $X_k/\Delta\bi_k$. 
Integrating over $J_k$ yields $|\Sigma^{n,m}(J_k)|=X_k$. 
Now if $\Delta\bi_k=0$, then $\Sigma^{n,m}$ has a jump of size $X_k$ at $E\rb_k$. Since $\Sigma^{n,m}(J_k)=\Sigma^{n,m}(\{E\rb_k\})=[\Sigma^{n,m}(E\rb_k-), \Sigma^{n,m}(E\rb_k)]$ by our definition of the image set, the conclusion follows.  
\end{proof}

\begin{proof}[Proof of Lemma~\ref{lem: dist-sur}]
\label{pf: dist-sur}
Recall that all surplus edges are of the form $(b_i, w_j)$ for some $k\ge 1$ so that $b_i=V_k$, $w_j\in \cA_{k-1}$ and 
\[
\bP_{n,m}((b_i, w_j)\in \cE_{n,m}\,|\,b_i=V_k, w_j\in \cA_{k-1})= 1-\exp(-X_{i} Y_{j}(k)/\sqrt{mn}). 
\]
Let us first show that the Poisson point measure $Q^{n,m}$ will produce an atom corresponding to the pair $(b_i, w_j)$ with the same probability. To that end, first recall from the LIFO-queue construction that $Y_j(k)$ is the remaining service time of $w_j$ when $V_k$ arrives, provided $w_j\in \cA_{k-1}$. Moreover, during the time when $\cN(V_k)$, the neighbours of $V_k$ are being served, all members of $\cA_{k-1}$ are waiting in the queue. 
Suppose that at time $t$, the queue consists of $w_{j_1}, w_{j_2}, \dots, w_{j_M}$ with $M=M(k)\ge 1$. 
Since $Z^{n,m}_t-\inf_{u\le t}Z^{n,m}_u$ equals the total service request from the clients in the queue, if some member of $\cN(V_k)$ is receiving service at time $t$, then $\mathrm Z_t:=Z^{n,m}_t-\inf_{u\le t}Z^{n,m}_u\ge \sum_{1\le q\le M}Y_{j_q}(k)$. As a result, we can partition the interval $[0, \mathrm Z_t]$ into $M+1$ disjoint subintervals $I_{t, 1}, I_{t, 2}, \dots, I_{t, M+1}$ with respective lengths $I_{t, q}=Y_{j_q}(k), 1\le q\le M$, and $I_{t, M+1}=\mathrm Z_t-\sum_{1\le l\le M}Y_{j_l}(k)$, which represents the service request from $\cN(V_k)$. Imagine these subintervals are stacked in order with $I_{t, 1}$ at the bottom, so that throughout the service time of $\cN(V_k)$, the bottom $M$ subintervals do not change at all. 
Next, recall that $J_k$ stands for the set of time when $V_k$ receives its service in the {\it bipartite LIFO-queue}. By definition, this is the set of time when members of $\cN(V_k)$ receive their service in the original LIFO-queue. 
For $1\le k\le K$ and $1\le q\le M(k)$, let us set
\[
D(V_k, w_{j_q})=\big\{(s, y): \exists\, t\in J_k \text{ s.t. } \Sigma^{n,m}(t-)\le s\le \Sigma^{n,m}(t), y\in I_{t, q}\big\}. 
\]
The previous arguments combined with Lemma~\ref{lem: sigma-p} yield that the Lebesgue measure of $D(V_k, w_{j_q})$ is $X_{i} Y_{j_q}(k)$. It follows that the probability of $Q^{n,m}$ containing at least one atom in $D(V_k, w_{j_q})$ is $ 1-\exp(-X_{i} Y_{j_q}(k)/\sqrt{mn})$. 

Next, let us argue that an atom from $D(V_k, w_{j_q})$ will yield an edge $e'\in \cE''_{n,m}$ between $V_k$ and the parent of $w_{j_q}$. Indeed, recall from~\eqref{def: cE''} the way $t_{n, l}$ and $t'_{n,l}$ are defined. Clearly, if $(s_{n,l}, y_{n,l})\in D(V_k, w_{j_q})$, then we have $\Sigma^{n,m}(t-)\le s_{n,l}\le \Sigma^{n,m}(t)$ for some $t\in J_k$; therefore $t_{n,l}\in J_k$ and $V_k$ is the client being served at time $t_{n,l}$. To identify $t'_{n,l}$, we note that by definition, this is the last moment before $t_{n,l}$ that the service load  is below the level $y_{n,q}$. By the LIFO-rule, this is the moment when the parent of $w_{j_q}$ arrived in the queue. 

To conclude, we note that $D(V_k, w_{j_q}), 1\le q\le M(k), 1\le k\le K$, are disjoint, so that the different elements of $\cE''_{n,m}$ are independent. Secondly, $Q^{n,m}$ can also contain atoms that correspond to $V_k$ and a member of $\cN(V_k)$, which will lead to the edge $(V_k, V_k)$ that is excluded from $\cE''_{n,m}$. 
\end{proof}

\subsubsection{Proof of Proposition~\ref{prop: cv-surplus}}
\label{sec: pf: cv-surplus}

We start with an observation: the convergences of $(\Lambda\rbn(t))_{t\ge 0}$ under $\bP_{n,m}$ in Lemma~\ref{lem: cv-dens} combined with Lemma~\ref{lem: Lambda-cond} immediately yield the same convergences under $\bP_{n,m}(\cdot\,|\,N\rb_n=n, N\rw_m=m)$, stated in the following lemma.

\begin{lem}
\label{lem: cvLambda}
Assume~\eqref{hyp: mod-clst} and~\eqref{hyp: crit}. 
\begin{enumerate}
\item
Under the assumptions of Proposition~\ref{prop: cv-ZH} (1), for all $t_0\ge 0$, we have
{\small\begin{equation}
\label{eqcv: Lambda-br}
\sup_{t\le t_0} \Big|n^{-\frac23}\Lambda\rbn(n^{\frac23} t)-\theta^{-\frac12}\sigma\rb_2 t\Big|\to 0  \text{ in prob.~under } \bP_{n,m}(\cdot\,|\,N\rb_n=n, N\rw_m=m). 
\end{equation}}
\item
Under the assumptions of Proposition~\ref{prop: cv-ZH} (2) or (3), for all $t_0\ge 0$, we have
{\small\begin{equation}
\label{eqcv: Lambda-power}
\sup_{t\le t_0} \Big|n^{-\frac{\alpha}{\alpha+1}}\Lambda\rbn(n^{\frac{\alpha}{\alpha+1}} t)-\theta^{-\frac12}\sigma\rb_2 t\Big|\to 0  \text{ in prob.~under } \bP_{n,m}(\cdot\,|\,N\rb_n=n, N\rw_m=m). 
\end{equation}}
\end{enumerate}
\end{lem}

Recall $\Sigma^{n,m}$ from~\eqref{def: Sigma}. As a first step in proving Proposition~\ref{prop: cv-surplus}, let us show the following result. 
\begin{lem}
\label{lem: cv-Sigma}
Assume~\eqref{hyp: mod-clst} and~\eqref{hyp: crit}. 
\begin{enumerate}
\item
Under the assumptions of Proposition~\ref{prop: cv-ZH} (1), for all $t_0\ge 0$, we have
{\small\begin{equation}
\label{eqcv: Sigma-br}
\sup_{t\le t_0} \Big|n^{-\frac23}\Sigma^{n,m}_{n^{2/3} t}-\theta^{-\frac12}\sigma\rb_2 t\Big|\to 0 \quad \text{in prob.~under } \bP_{n,m}(\cdot\,|\,N\rb_n=n, N\rw_m=m). 
\end{equation}}
\item
Under the assumptions of Proposition~\ref{prop: cv-ZH} (2) or (3), for all $t_0\ge 0$, we have
{\small\begin{equation}
\label{eqcv: Sigma-power}
\sup_{t\le t_0} \Big|n^{-\frac{\alpha}{\alpha+1}}\Sigma^{n,m}_{ n^{\alpha/(\alpha+1)} t}-\theta^{-\frac12}\sigma\rb_2 t\Big|\to 0 \quad \text{in prob.~under } \bP_{n,m}(\cdot\,|\,N\rb_n=n, N\rw_m=m). 
\end{equation}}
\end{enumerate}
\end{lem}

\begin{proof}
To ease the notation, let us denote $\hat \bP=\bP_{n,m}(\cdot\,|\,N\rb_n=n,N\rw_m=m)$. 

\noindent
{\bf Proof of~\eqref{eqcv: Sigma-br}.} 
Recall from~\eqref{def: Wt} the set $\mathfrak{W}_t$ of clients waiting in the bipartite queue at time $t\ge 0$. Let 
\[
\hat\Sigma_t :=  \sum_{1\le i\le N\rb_n} X_i\mathbf 1_{\{i\in \mathfrak{W}_t\}},
\]
Since the customers that have either departed before $t$ or are currently waiting in the queue must have arrived by $t$, we deduce that for all $t\ge 0$, 
\begin{equation}
\label{eqin: sigma}
\Lambda\rbn(t)\ge \Sigma^{n,m}_t \ge \Lambda\rbn(t)-\hat\Sigma_t. 
\end{equation}
We claim that under the assumptions of Proposition~\ref{prop: cv-ZH} (1), 
\begin{equation}
\label{eqbd: sigma}
 n^{-\frac23}\sup_{t\le n^{2/3} t_0}\hat\Sigma_t \to 0 \quad \text{in probability under } \hat\bP, 
\end{equation}
which will then imply~\eqref{eqcv: Sigma-br}, thanks to~\eqref{eqin: sigma} and the convergence in~\eqref{eqcv: Lambda-br}. 
To show~\eqref{eqbd: sigma}, let us recall that $H^{n,m}_t$ corresponds to the queue length at time $t$. Take $\epsilon>0$ and $L>0$. We have
\begin{equation}
\label{bd: sig}
\hat\bP\Big(\sup_{t\le n^{2/3} t_0}\hat\Sigma_t \ge \epsilon n^{\frac23}\Big)\le \hat \bP(\sup_{t\le n^{2/3}t_0} H^{n,m}_t\ge L n^{\frac13}\Big) + \hat \bP\Big(X^{\ast}(n^{2/3}t_0)\ge L^{-1}\epsilon n^{\frac13}\Big),
\end{equation}
where $X^{\ast}(n^{2/3}t_0)=\max_{1\le i\le N\rb_n}X_i\mathbf 1_{\{E\rb_i\le n^{2/3}t_0\}}$. To control this quantity, let us introduce for $\delta>0$: 
\[
J(\delta)=\#\big\{1\le i\le N\rb_n: E\rb_i\le n^{\frac23}t_0, X_i\ge \delta n^{\frac13}\big\}.
\]
Then $X^{\ast}(n^{2/3}t_0)\ge \delta n^{1/3}$ if and only if $J(\delta)\ge 1$. Meanwhile, $J(\delta)$ under $\bP_{n,m}$ is a Poisson random variable with mean
\begin{align*}
\alpha_n &=\sqrt{\frac{n}{m}}\int_0^{n^{2/3}t_0}\!\!\!\int_{[\delta n^{1/3},\infty)} xe^{-xs/\sqrt{mn}}dF\rb(x)ds \\
&\le \sqrt{\frac{n}{m}}\int_0^{\infty}\int_{[\delta n^{1/3},\infty)}xe^{-xs/\sqrt{mn}}dF\rb(x)ds \\
& = n\int_{[\delta n^{1/3},\infty)} dF\rb(x)\le \delta^{-3}\int_{[\delta n^{1/3},\infty)}x^3 dF\rb(x),
\end{align*}
which tends to $0$ as $n\to\infty$, since $dF\rb$ has finite third moment. This shows that $J(\delta)\to 0$ in probability under $\bP_{n,m}$, which in turn implies that 
$\bP_{n,m}(X^{\ast}(n^{2/3}t_0)\ge \delta n^{\frac13})\to 0$,
for each fixed $\delta>0$. Since $X^{\ast}(n^{2/3}t_0)$ is a measurable function of the collection $\cN_{t_0}$, Lemma~\ref{lem: Lambda-cond} yields that 
\[
\hat\bP\Big(X^{\ast}(n^{2/3}t_0)\ge \delta n^{\frac13}\Big)\to 0, \quad\text{as } n\to\infty.
\]
Plugging this into~\eqref{bd: sig}, we find that
\[
\limsup_{n\to\infty}\hat\bP\Big(\sup_{t\le n^{2/3} t_0}\hat\Sigma_t \ge \epsilon n^{\frac23}\Big)\le \limsup_{n\to\infty}\hat \bP(\sup_{t\le n^{2/3}t_0} H^{n,m}_t\ge L n^{\frac13}\Big).
\]
According to Proposition~\ref{prop: cv-ZH} (1), $(n^{-1/3}H^{n,m}_{n^{2/3}t})_{t\ge 0}$, $n\ge 1$ is tight. Then~\eqref{eqbd: sigma} follows by taking $L\to\infty$. 

\medskip
\noindent
{\bf Proof of~\eqref{eqcv: Sigma-power}.} We introduce $\hat\Sigma_t$ as in the previous case. Note that~\eqref{eqin: sigma} still holds. In view of the convergence in~\eqref{eqcv: Lambda-power},  it remains to show that
\begin{equation}
\label{eqbd: sigma'}
 n^{-\frac{\alpha}{\alpha+1}}\sup_{t\le n^{\alpha/(\alpha+1)} t_0}\hat\Sigma_t \to 0 \quad \text{in probability under } \hat\bP. 
\end{equation}
Let $X_{\sigma(i)}$ be the $i$-th largest among $\{X_i: 1\le i\le N\rb_n, E\rb_i\le n^{\alpha/(\alpha+1)}t_0\}$ (breaking ties arbitrarily). Clearly, for each $N\in \N$, we have
\[
\hat\Sigma_t \le \sum_{1\le i\le N}X_{\sigma(i)}+ \sum_{N<i\le N\rb_n} X_{\sigma(i)}\mathbf 1_{\{\text{$\sigma(i)$ is in the queue at time $t$}\}}.
\]
It follows that 
\[
\hat\bP\Big(\sup_{t\le n^{\alpha/(\alpha+1)} t_0}\hat\Sigma_t \ge \epsilon n^{\frac{\alpha}{\alpha+1}}\Big)\le P_1+P_2+P_3,
\]
where
\begin{align*}
P_1&=\hat\bP\Big(\sum_{1\le i\le N}X_{\sigma(i)}\ge\epsilon n^{\frac{\alpha}{\alpha+1}}\Big),\ 
P_2= \hat \bP\Big(\sup_{t\le  n^{\alpha/(\alpha+1)}t_0} H^{n,m}_t\ge L n^{\frac{\alpha-1}{\alpha+1}}\Big),\\
P_3 &= \hat \bP\Big(X^{\ast}_N(n^{\frac{\alpha}{\alpha+1}}t_0)\ge L^{-1}\epsilon n^{\frac{1}{\alpha+1}}\Big)\quad \text{with}\quad  X^{\ast}_N(n^{\frac{\alpha}{\alpha+1}}t_0)=\max_{N<i\le N\rb_n}X_{\sigma(i)}\mathbf 1_{\{E\rb_{\sigma(i)}\le n^{\frac{\alpha}{\alpha+1}}t_0\}}. 
\end{align*}
This time let us set for each $a>0$, 
\[
J(a )=\#\{1\le i\le N\rb_n: E\rb_i\le n^{\frac{\alpha}{\alpha+1}}t_0, X_i\ge a \},
\]
which, under $\bP_{n,m}$, has a Poisson distribution of mean
\[
\alpha_n(a) = \sqrt{\frac{n}{m}}\int_0^{ n^{\alpha/(\alpha+1)}t_0} \int_{[a, \infty)}xe^{-xs/\sqrt{mn}}dF\rb(x)ds. 
\]
In the first instance, we take $a=\delta  n^{\alpha/(\alpha+1)}$ with some fixed $\delta>0$. Note that 
\begin{align*}
\alpha_n(\delta n^{\frac{\alpha}{\alpha+1}}) 
&\le \sqrt{\frac{n}{m}}n^{\frac{\alpha}{\alpha+1}}t_0 \int_{[\delta  n^{\alpha/(\alpha+1)}, \infty)} xdF\rb(x) \le \delta^{-1} t_0 \sqrt{\frac{n}{m}}  \int_{[\delta n^{\alpha/(\alpha+1)}, \infty)} x^2dF\rb(x). 
\end{align*}
Since $m=\lfloor \theta n\rfloor$ and $dF\rb$ has finite second moment, we deduce from the above  that $\alpha_n(\delta  n^{\alpha/(\alpha+1)})\to 0$ as $n\to\infty$, which implies that for each $\delta>0$, 
\[
\bP_{n,m}\Big(X_{\sigma(1)}\ge \delta n^{\frac{\alpha}{\alpha+1}}\Big) \to 0, \quad n\to\infty.
\]
Once again, we can substitute in above $\bP_{n,m}$ with $\hat\bP$ thanks to Lemma~\ref{lem: Lambda-cond}. It then follows that for each fixed $N$, $P_1\to 0$ as $n\to\infty$.  Next, let us take $a= \delta  n^{1/(\alpha+1)}$, which yields 
\[
\alpha_n(\delta n^{\frac{1}{\alpha+1}})\le n\int_{[ \delta  n^{1/(\alpha+1)},\infty)} dF\rb(x)= n\big(1-F(\delta n^{\frac{1}{\alpha+1}})\big)\le  2C\rb \delta^{-1-\alpha} 
\]
for sufficiently large $n$, 
where we have used that under~\eqref{hyp: b-power}, we have $1-F\rb(x)\le 2C\rb x^{-1-\alpha}$ for sufficiently large $x$. We note that
\begin{align*}
\bP_{n,m}\big(X^{\ast}_N(n^{\frac{\alpha}{\alpha+1}}t_0)\ge \delta n^{\frac{1}{\alpha+1}}\big) &= \bP_{n,m}\big(J(\delta n^{\frac{1}{\alpha+1}}) > N)\\
&\le N^{-1} \bE_{n,m}\big[J(\delta n^{\frac{1}{\alpha+1}})\big] = N^{-1}\alpha_n(\delta n^{\frac{1}{\alpha+1}}). 
\end{align*}
Thanks to the previous bound, we deduce that 
\[
\limsup_{N\to\infty}\limsup_{n\to\infty}\bP_{n,m}\big(X^{\ast}_N(n^{\frac{\alpha}{\alpha+1}}t_0)\ge \delta n^{\frac{1}{\alpha+1}}\big) =0.
\]
Lemma~\ref{lem: Lambda-cond} then allow us to conclude that $\limsup_{N\to\infty}\limsup_{n\to\infty}P_3 = 0$. Finally, the respective convergences in~\eqref{cv: ZH2} and~\eqref{cv: ZH3} imply that $\lim_{L\to\infty}\limsup_{n\to\infty} P_2=0$. This completes the proof of~\eqref{eqbd: sigma'} and then~\eqref{eqcv: Sigma-power} as well.
\end{proof}

\begin{proof}[Proof of Proposition~\ref{prop: cv-surplus}]
Let $a_n= n^{1/(\alpha+1)}$, $b_n= n^{\alpha/(\alpha+1)}$ with the understanding that $\alpha=2$ under the double finite third moments assumptions. 
By Skorokhod's Representation theorem, we can assume that the convergences of $(Z^{n,m}, H^{n,m})$ in Proposition~\ref{prop: cv-ZH} and $\Sigma^{n,m}$ in Lemma~\ref{lem: cv-Sigma} take place almost surely under the relevant assumptions. Let us denote $\Sigma^{-1}(s)=\inf\{t: \Sigma^{n,m}(t)> s\}$. Then Lemma~\ref{lem: cv-Sigma} combined with the convergence of inverse functions (see for instance Theorem 7.2 in~\cite{Whitt}) implies that for any $s_0\ge 0$, 
\begin{equation}
\label{eqcv: inv-Sigma}
\sup_{s\le b_n s_0}\Big|b_n^{-1}\Sigma^{-1}(b_n s)-\rho^{-1}s \Big|\to 0 \quad\text{a.s.}
\end{equation}
where we have denoted $\rho=\sigma\rb_2/\sqrt\theta$. 
Let us show that the Poisson point measure
\[
\widetilde Q^{n,m}:=\{(b_n^{-1}s_{n, l}, a_n^{-1}y_{n,l}): 1\le l\le q_{n,m}\}
\]
converges in distribution to $\cQ^{(i)}$ on any compact set of $\R^2$, under the relevant assumptions for $i\in \{1, 2, 3\}$. 
Since Poisson point measures converge if and only if their intensity measures converge accordingly, it suffices to show that for any continuous function $f: \R_+^2\to \R_+$, we have for $i\in \{1, 2, 3\}$, 
\begin{equation}
\label{cv: pp-int}
\frac{1}{\sqrt{mn}}\int_{\R^2_+} f(b_n^{-1}s, a_n^{-1}y)\mathbf 1_{D_{n,m}}(s, y)dsdy \xrightarrow[n\to\infty]{\text{a.s.}}
\frac{1}{\sqrt\theta}\int_{\R^2_+} f(s, y)\mathbf 1_{\cD^{(i)}}(s, y)ds dy.
\end{equation}
where we recall $D_{n,m}$ and $\cD^{(i)}$ are the respective supports of $Q^{n,m}$ and $\cQ^{(i)}$. 
With a change of variables and noting $a_n b_n=n$, we can rewrite the left-hand side of~\eqref{cv: pp-int} as follows:
\begin{align*}
&\quad \frac{1}{\sqrt{mn}}\int_{\R^2_+} f(b_n^{-1}s, a_n^{-1}y)\mathbf 1_{D_{n,m}}(s, y)dsdy  =\frac{n}{\sqrt{mn}}\int_{\R^2_+} f(s', y')
\mathbf 1_{D_{n,m}}(b_n s', a_n y')ds' dy'\\
& = \frac{n}{\sqrt{mn}}\int_{\R^2_+} f(s', y')
\mathbf 1_{\{a_ny'\le Z^{n,m}\circ \Sigma^{-1}(b_ns')-\inf_{u\le b_ns'}Z^{n,m}\circ\Sigma^{-1}(u)\}}ds' dy',\\
& \xrightarrow[n\to\infty]{\text{a.s.}} \frac{1}{\sqrt\theta}\int_{\R^2_+}f(s', y')\mathbf 1_{\{y'\le \cZ^{(i)}_{\rho^{-1}s'}-\inf_{v\le s'/\rho}\cZ^{(i)}_v\}}ds'dy',
\end{align*}
which is the right-hand side of~\eqref{cv: pp-int}. In the second line above, we have implicitly used the fact that for almost every $(s, y)\in D_{n,m}$, if $\Sigma^{n,m}(t-)\le s\le \Sigma^{n,m}(t)$, then $t=\Sigma^{-1}(s)$; in the final line, we have applied the assumption $m=\lfloor \theta n\rfloor$, the aforementioned convergence of $(Z^{n,m}, H^{n,m})$,~\eqref{eqcv: inv-Sigma} and Lemma~\ref{lem: cv-comp}. 

Suppose that the atoms $\{(s_l, y_l), l\ge 1\}$ of $\cQ^{(i)}$ are ranked in increasing order of the first coordinate, and similarly let $\{(s_{n,l}, y_{n,l}), l\ge 1\}$ be the atoms of $Q^{n,m}$, ranked in increasing order of the first coordinate. The previous convergence of $Q^{n,m}$ then implies that for any $L\ge 1$, we have
\[
\{(b_n^{-1}s_{n, l}, a_n^{-1}y_{n, l}): 1\le l\le L\} \Longrightarrow \{(s_l, y_l): 1\le l\le L\}.
\]
Let us deduce the convergence of $(b_n^{-1}t_{n, l}, b_n^{-1}t'_{n,l}), 1\le l\le L$. Indeed, recall that $t_{n,l}$ satisfies $\Sigma^{n,m}(t_{n,l}-)\le s_{n,l}\le \Sigma^{n,m}(t_{n,l})$ and a.s.~$\Sigma^{n,m}$ is strictly increasing on a neighbourhood of $t_{n,l}$, so that $t_{n,l}=\Sigma^{-1}(s_{n,l})$ a.s. 
Combining this with~\eqref{eqcv: inv-Sigma}, we obtain that $\{b_n^{-1}t_{n,l}: 1\le l\le L\}$ converges in distribution to $\{t_l=\rho^{-1}s_l: 1\le l\le L\}$. Recall that
\[
t'_{n,l}=\sup\Big\{u\le t_{n,l}: Z^{n,m}_u-\inf_{v\le u}Z^{n,m}_v\le y_{n,l}\Big\}, \quad t'_l = \sup\Big\{u\le t_l: \cZ^{(i)}_u-\inf_{v\le u}\cZ^{(i)}_v\le y_l\Big\}. 
\]
The respective convergences of $t_{n,l}$ to $t_l$, of $y_{n,l}$ to $y_l$, of $Z^{n,m}$ to $\cZ^{(i)}$, combined with the sample path properties of $\cZ^{(i)}$, imply that
\[
\{(b_n^{-1}t_{n, l}, b_n^{-1}t'_{n, l}): 1\le l\le L\} \Longrightarrow \{(t_l, t'_l): 1\le l\le L\}.
\]
This completes the proof.
\end{proof}

\subsection{Proof of the main theorems}
\label{sec: pf-main}

Having shown Proposition~\ref{prop: cv-ZH} for the convergence of graph encoding processes and Proposition~\ref{prop: cv-surplus} for the convergence of surplus edges, we explain here how we obtain our main results Theorems~\ref{thm: thm1}-\ref{thm: thm3} and Proposition~\ref{prop: rnk}. We will focus on the double finite third moment scenario, as the other two cases can be argued in a similar fashion.

Recall that $(g_k, d_k)$ is the $k$-th longest excursion interval of $\cZ^{(1)}$ above its running infimum. 
Recall that $C^{\bx}_{n,(k)}$ is the $k$-th largest connected components of $B_{n,m}$ ranked in $\bx$-weights. We note from Lemma~\ref{lem: exc} that there is a unique excursion interval of $Z^{n,m}$ above its running infimum, say $(g_{n, (k)}, d_{n, (k)})$,  that is associated with $C^{\bx}_{n, (k)}$ and satisfies 
\begin{equation}
\label{def: gnk}
d_{n, (k)}-g_{n, (k)}=\by\big(C^{\bx}_{n,(k)}\big) \text{ and }\Lambda\rbn\big(d_{n, (k)}\big)-\Lambda\rbn\big(g_{n, (k)}-\big)=\bx\big(C^{\bx}_{n, (k)}\big). 
\end{equation}
We start with the following convergences. 

\begin{lem}
\label{lem: cv-int}
Under the assumptions of Proposition~\ref{prop: cv-ZH} (1), the following waek convergences take place jointly with the convergence in~\eqref{cv: ZH1}: we have 
\begin{equation}
\label{cv: l2-x}
\Big\{n^{-\frac23}\bx(C^{\bx}_{n, (k)}): 1\le k\le \kappa_{n,m}\Big\} \Longrightarrow \Big\{\tfrac{\sigma\rb_2}{\sqrt\theta}(d_k-g_k): k\ge 1\Big\}
\end{equation}
with respect to the $\ell^2$-topology, and 
\begin{equation}
\label{cv: ex-int}
\Big\{\Big(n^{-\frac23}g_{n,(k)}, n^{-\frac23}d_{n,(k)}\Big): 1\le k\le \kappa_{n,m}\Big\} \Longrightarrow \{(g_k, d_k): k\ge 1\}
\end{equation}
with respect to the product topology of $(\R^2)^{\N}$. 
\end{lem}

\begin{proof}
We will use the shorthand notation $\hat\bP=\bP_{n,m}(\cdot\,|\,N\rb_n=n, N\rw_m=m)$. 
By Skorokhod's Representation Theorem, we can assume that the convergence~\eqref{cv: ZH1} of $(Z^{n,m}, H^{n,m})$ in Proposition~\ref{prop: cv-ZH} and the convergence~\eqref{eqcv: Lambda-br} of $\Lambda\rbn$ both take place {\it almost surely}. 
The idea of the proof is to use Aldous' size-biased point process technique. However, in our LIFO-queue construction, connected components appear in size-biased order with respect to their $\bx$-weights, while the excursion lengths in $Z^{n,m}$ correspond to their $\by$-weights. We therefore need a slight adaptation. Let us consider all the connected components of $B_{n,m}$ that contain at least one black vertex, ranked in decreasing order of their $\bx$-weights: $\{C'_{n, (k)}: 1\le k\le \hat\kappa_{n,m}\}$. In particular, $\{C^{\bx}_{n, (k)}: 1\le k\le \kappa_{n,m}\}$ appears as a subsequence. For each $C'_{n, (k)}$, we can associate a unique and possibly empty interval $[g'_{n, (k)}, d'_{n, (k)})$ so that
\[
\Lambda\rbn\big(d'_{n, (k)}\big)-\Lambda\rbn\big(g'_{n, (k)}-\big)=\bx\big(C'_{n, (k)}\big); 
\]
and if $C'_{n,(k)}$ only contains a single black vertex, then $d'_{n, (k)}=g'_{n, (k)}$ is a jump time of $\Lambda\rbn$. Let us first show that 
\begin{equation}
\label{cv: l2-x'}
\Big\{n^{-\frac23}\bx(C'_{n, (k)}): 1\le k\le \hat\kappa_{n,m}\Big\} \Longrightarrow \Big\{\tfrac{\sigma\rb_2}{\sqrt\theta}(d_k-g_k): k\ge 1\Big\}
\end{equation}
with respect to the $\ell^2$-topology. To that end, let us consider the following point process: 
\[
\Xi_{n,m}=\Big\{\Big(n^{-\frac23}\Lambda\rbn\big(g'_{n, (k)}-\big), n^{-\frac23}\Lambda\rbn\big(d'_{n,(k)}\big)-n^{-\frac23}\Lambda\rbn\big(g'_{n,(k)}-\big)\Big): 1\le k\le \hat\kappa_{n,m}\Big\}.
\]
Note that the second coordinate corresponds to the rescaled $\bx$-weight of $C'_{n, (k)}$. Furthermore, $\Lambda\rbn(t-)$ counts the $\bx$-weights of all the black vertices that have arrived prior to time $t$. 
Therefore, $\Xi_{n,m}$ is a size-biased point process in the sense of Aldous~\cite{Al97}.
For its continuum counterpart, denote $\rho = \sigma\rb_2/\sqrt\theta$ and let 
\[
\Xi= \big\{\big(\rho g_k, \rho(d_k-g_k) \big): k\in\N\big\}.
\]
Properties of Brownian motion then imply that almost surely $\Xi$ has  finitely many elements on any compact of $\R_+\times (0, \infty)$. 

We claim that $\Xi_{n,m}$ converges to $\Xi$ in the vague topology of $\R_+\times (0, \infty)$. This boils down to proving the following convergence for any fixed $T, \epsilon>0$. 
Let $(g^T_j, d^T_j), 1\le j\le J$ be those among $(g_k, d_k), k\ge 1$ that satisfy $g_k\le \rho^{-1}T$ and $d_k-g_k>\rho^{-1}\epsilon$, ranked in increasing order of their left endpoints; 
similarly, let  $(g^T_{n, j}, d^T_{n,j}), 1\le j\le J_n$ be  those among $(g'_{n, (k)}, d'_{n, (k)}), 1\le k\le \hat\kappa_{n,m}$ that satisfy $\Lambda\rbn(g'_{n, (k)}-)\le Tn^{2/3}$ and $\Lambda\rbn(d'_{n, (k)})-\Lambda\rbn(g'_{n, (k)}-)>\epsilon n^{2/3}$, ranked in increasing order of their left endpoints; 
we then have 
\begin{equation}
\label{cv: sbpp}
J_n\to J  \ \text{ and for }1\le j\le J:  n^{-\frac23}\Lambda\rbn(g^T_{n, j}-)\to \rho g^T_j, \quad n^{-\frac23}\Lambda\rbn(d^T_{n, j})\to \rho d^T_j.  
\end{equation}
To see why~\eqref{cv: sbpp} is true, we note that 
thanks to the uniform convergence~\eqref{eqcv: Lambda-br} of $\Lambda\rbn$, for $\delta>0$ and $n$ taken sufficiently large, $\{(g^T_{n,j}, d^T_{n, j}): 1\le j\le J_n\}$ is contained in the sub-collection of those among $(g'_{n, (k)}, d'_{n, (k)}), 1\le k\le \hat\kappa_{n,m}$ that satisfy $g'_{n, (k)}\le (1+\delta)\rho^{-1} Tn^{2/3}$ and $d'_{n, (k)}-g'_{n, (k)}>(1-\delta)\rho^{-1}\epsilon n^{2/3}$. 
Since almost surely $d_k-g_k, k\ge 1$ take distinct values, we can choose the previous $\delta$ so that with probability exceeding $1-\epsilon$, $(g^T_j, d^T_j), 1\le j\le J$ are  exactly those among $(g_k, d_k), k\ge 1$ that satisfy $g_k\le (1+\delta) \rho^{-1}T$ and $d_k-g_k>(1-\delta)\rho^{-1}\epsilon $. 
Now the convergences~\eqref{cv: ZH1} and~\eqref{eqcv: Lambda-br} combined with standard arguments (see for instance Lemma 7 in~\cite{Al97}) imply that almost surely
\begin{equation}
\label{cv: gd}
J_n\to J \ \text{ and for }1\le j\le J:   n^{-\frac23} g^T_{n,j} \to g^T_j, n^{-\frac23} d^T_{n,j} \to d^T_j.
\end{equation}
With another application of the convergence~\eqref{eqcv: Lambda-br}, the convergence in~\eqref{cv: sbpp} now follows. 
Aldous' theory of size-biased point processes (Proposition 15 in~\cite{Al97}) allows us to deduce the $\ell^2$-convergence of the second coordinates of $\Xi_{n,m}$, which is precisely~\eqref{cv: l2-x'}.

Next, let us show that for each $K\ge 1$, 
\begin{equation}
\label{eq: xrk}
\hat{\bP}\Big(C'_{n, (k)}= C^{\bx}_{n, (k)}, 1\le k\le K\Big) \xrightarrow{n\to\infty} 1.
\end{equation} 
Indeed, by standard properties of the Brownian motion, we can  choose $T, \delta>0$ so that $\mathbb P(\max_{1\le k\le K}g_{k}\le T, d_K-g_K\ge \delta)\ge 1-\epsilon$. Together with the $\ell^2$-convergence in~\eqref{cv: l2-x'}, ~\eqref{cv: gd} then implies that 
\[
\limsup_{n\to\infty}\mathbb P\Big( \min_{1\le k\le K}\big(d'_{n, (k)}-g'_{n, (k)}\big) <\tfrac12\delta n^{2/3}\Big)\le \epsilon. 
\]
Since $d'_{n, (k)}-g'_{n, (k)}=\by(C'_{n, (k)})$, this shows that with high probability, $C'_{n, (k)}$ is nontrivial for each $1\le k\le K$, from which~\eqref{eq: xrk} follows. As a consequence, we deduce that the convergence in~\eqref{cv: l2-x} holds with respect to the product topology. Meanwhile, the left-hand side of~\eqref{cv: l2-x} is tight in the $\ell^2$-topology, since it is a subsequence of $\{n^{-2/3}\bx(C'_{n, (k)}): 1\le k\le \hat \kappa_{n,m}\}$, whose $\ell^2$-convergence we have seen in~\eqref{cv: l2-x'}. The $\ell^2$-convergence in~\eqref{cv: l2-x} readily follows. 

Finally, let us show~\eqref{cv: ex-int}. Appealing to Skorokhod's Representation Theorem again, we can assume that the convergence in~\eqref{cv: l2-x'} also takes place {\it almost surely}.  
In particular, this implies that for all $\varepsilon>0, K\ge 1$, we can find such $\delta>0$ that 
\begin{equation}
\label{eq: gevent}
\limsup_{n\to\infty}\hat{\bP}\Big(\bx(C^{\bx}_{n,(K)})<\delta n^{\frac23}\Big) < \frac{\varepsilon}{2}.
\end{equation}
On the other hand,  note that $\chi:=\sum_{k\ge 1}\rho^2(d_k-g_k)^2<\infty$ a.s.~and recall that $\{(g^T_j, d^T_j), 1\le j\le J\}$ is the sub-collection of the excursions that satisfy $g_k\le \rho^{-1}T$ and $d_k-g_k>\rho^{-1}\epsilon$. We can choose $T$ and $\epsilon$ so that 
\[
\mathbb P\Big(\sum_{1\le j\le J}\rho^2(d^T_j-g^T_j)^2\le \chi-\tfrac{\delta}{4}\Big)\le\frac{\varepsilon}{2}. 
\]
Combining this with~\eqref{cv: sbpp}, we find that
\[
\limsup_{n\to\infty}\hat{\bP}\Big(n^{-\frac43}\sum_{1\le j\le J_n}\big(\Lambda\rbn(d^T_{n,j})-\Lambda\rbn(g^T_{n,j}-)\big)^2\le \chi-\tfrac{\delta}{2}\Big)\le \frac{\varepsilon}{2}.
\]
Applying the almost sure $\ell^2$-convergence in~\eqref{cv: l2-x'} yields  
\begin{equation}
\label{bd: expbyT}
\limsup_{n\to\infty}\hat{\bP}\Big(n^{-\frac43}\sum_{j> J_n}\big(\Lambda\rbn(d^T_{n,j})-\Lambda\rbn(g^T_{n,j}-)\big)^2\ge \delta\Big)\le \frac{\varepsilon}{2}. 
\end{equation}
Comparing this with~\eqref{eq: gevent}, we conclude that there exists some $T$ so that 
\[
\limsup_{n\to\infty}\hat\bP\big(\exists\, k\le K: g_{n, (k)}>\rho^{-1}Tn^{2/3}\big)\le \varepsilon. 
\]
The convergence in~\eqref{cv: ex-int} is now a consequence of this and~\eqref{cv: gd}. 
\end{proof}

Recall that $C^{\by}_{n,(k)}$ is the $k$-th largest nontrivial connected components of $B_{n,m}$ in $\by$-weights. 
Using the symmetry of the bipartite model in black and white vertices, we next show the following statement.
\begin{lem}
\label{lem: rkcons}
Let $K\in\N$. Under the assumptions of Proposition~\ref{prop: cv-ZH} (1), we have 
\[
\limsup_{n\to\infty}\mathbf P_{n,m}\Big(\exists\, k\le K: C^{\by}_{n, (k)}\ne C^{\bx}_{n, (k)}\,\Big|\,N\rb_n=n, N\rw_m=m\Big)=0.
\] 
\end{lem}

\begin{proof}
We denote $\hat\bP=\bP_{n,m}(\cdot\,|\,N\rb_n=n, N\rw_m=m)$. Recall that $\zeta_k=d_k-g_k$ is the length of the $k$th longest excursion of $\cZ^{(1)}$. 
As a first step, let us show that 
\begin{equation}
\label{cv: comp}
\big\{n^{-2/3}\by(C^{\bx}_{n, (k)}): k\ge 1\big\}\Rightarrow \big\{\zeta_k, k\ge 1\big\},\  \big\{n^{-2/3}\by(C^{\by}_{n, (k)}): k\ge 1\big\}\Rightarrow \big\{\zeta_k, k\ge 1\big\}
\end{equation}
The first convergence is immediate from~\eqref{def: gnk} and~\eqref{cv: ex-int}. To prove the second convergence, let $\tilde B_{n,m}$ be the graph obtained from $B_{n,m}$ by swapping the vertex colours; then $\by(C^{\by}_{n, (k)})$ corresponds to the total black weights of the $k$-th largest component in the black weights in $\tilde B_{n,m}$. On the other hand, $\tilde B_{n,m}$ satisfies the conditions of Lemma~\ref{lem: cv-int}, as it is a random bipartite graph on $m$ black vertices and $n$ white vertices, where the black weights (resp.~white weights) are sampled from $F\rw$ (resp.~from $F\rb$). Lemma~\ref{lem: cv-int} combined with the previous colour swapping argument implies that 
\begin{equation}
\label{cv: l2-ywt}
\Big\{m^{-\frac23}\by(C^{\by}_{n, (k)}): 1\le k\le \kappa_{n,m}\Big\} \Longrightarrow \big\{\rho^{-1}\tilde \zeta_k, k\ge 1\big\},
\end{equation}
where $\rho=\sigma\rb_2/\sqrt\theta$ and $\tilde\zeta_k$ is the length of the $k$-th longest excursion of 
$\tilde Z_t - \inf_{s\le t}\tilde Z_s, t\ge 0$. Here, $(\tilde Z_t)_{t\ge 0}$ is the Brownian motion with parabolic drift given by 
\[
\mathbb E[e^{-\lambda \tilde Z_t}]=\exp\big\{\tfrac{\mathrm C_1}{2}(t\lambda^2\rho+\sqrt\theta \lambda t^2)\big\}, \quad \lambda\ge 0.
\]
It follows that 
\[
\{\tilde Z_t, t\ge 0\}\eqd \big\{\rho\theta^{-\frac13}Z^{(1)}_{\theta^{2/3}t/\rho}, t\ge 0\big\}, \quad \text{so that} \quad \tilde \zeta_k \eqd \rho\,\theta^{-\frac23}\zeta_k.
\]
Combined with $n^{-2/3}/m^{-2/3}\to \theta^{2/3}$, this proves~\eqref{cv: comp}. 

Next, let us prove by induction that $C^{\by}_{n, (k)}=C^{\bx}_{n, (k)}$ with high probability. For $k=1$, we have seen in~\eqref{cv: comp} that $n^{-2/3}\by(C^{\by}_{n, (2)})$ converges in distribution to $\zeta_2$. 
Standard properties of the Brownian motion imply that $\zeta_1>\zeta_2$ a.s. Then we can find some $z\in \R_+$ and $\epsilon>0$ so that 
\[
\mathbb P(\zeta_1\ge z) \ge \mathbb P(\zeta_2\ge z) +\epsilon. 
\]
On the event $\{C^{\by}_{n, (1)}\ne C^{\bx}_{n, (1)}\}$, we have $\by(C^{\bx}_{n, (1)})\le \by(C^{\by}_{n, (2)})$. It follows that 
\[
\limsup_{n\to\infty} \hat\bP\big(\by(C^{\bx}_{n, (1)})\ge n^{2/3}z\,\big|\, C^{\by}_{n, (1)}\ne C^{\bx}_{n, (1)}\big) \le \mathbb P(\zeta_2\ge z)\le \mathbb P(\zeta_1\ge z)-\epsilon. 
\]
Comparing this with~\eqref{cv: comp}, we conclude that $\hat\bP(C^{\by}_{n, (1)}\ne C^{\bx}_{n, (1)})\to 0$ as $n\to\infty$. 

Suppose that we have proved the claim up to $k-1\ge 1$. Then on the event $\{\forall\,j\le k-1: C^{\by}_{n, (j)}= C^{\bx}_{n, (j)}\}\cap\{C^{\by}_{n, (k)}\ne C^{\bx}_{n, (k)}\}$, we will have $\by(C^{\bx}_{n, (k)})\le \by(C^{\by}_{n, (k+1)})$. Together with the fact that $\zeta_k>\zeta_{k+1}$ a.s., we can argue as previously and complete the induction step. 
\end{proof}

\smallskip
\noindent
{\bf Encoding the bipartite graph.}
Lemmas~\ref{lem: cv-int} and~\ref{lem: rkcons} allow us to identify the excursion intervals that correspond to $C^{\bx}_{n, (k)}$ and $C^{\by}_{n, (k)}$. We now follow the encoding introduced in Section~\ref{sec: ghp} to introduce the coding functions for the components. 
Recall from~\eqref{def: ZH} the pair $(Z^{n,m}, H^{n,m})$. In particular, $H^{n,m}$ only takes  finitely many values on any compact set of $\R_+$, and therefore satisfies the conditions of $h$ in Section~\ref{sec: ghp}. 
Recall from~\eqref{def: gnk} that $(g_{n, (k)}, d_{n, (k)})$ is the excursion interval associated with $C^{\bx}_{n, (k)}$. 
Let $\mathrm{H}^{n, k}$ denote the portion of the rescaled $H^{n,m}$ running on this interval: 
\[
\mathrm{H}^{n, k}(t)= n^{-\frac13}H^{n,m}_{tn^{2/3}+g_{n,(k)}}, \quad 0\le t\le n^{-\frac23}(d_{n,(k)}-g_{n,(k)}).
\]
Recall the forest $\cF$ obtained from the LIFO-queue construction in Section~\ref{sec: explore}.
Let $T_{n, (k)}$ be the connected component of $\cF$ that has the same vertex set as $C^{\bx}_{n, (k)}$, so that $T_{n,(k)}$ is a spanning tree of $C^{\bx}_{n,(k)}$. Recall the measure $\mu^{\by}_{n, k}$, which assigns an atom of size $Y_j$ to each of the white vertex $w_j$ contained in $C^{\bx}_{n, (k)}$. 
Recall from~\eqref{def: Delta} the quantity $\Delta\bi_i$ and let $\hat\mu^{\by}_{n, k}$ be the measure which assigns an atom of size $\Delta\bi_i$ to each black vertex $b_i$ in $C^{\bx}_{n, (k)}$. Put differently, $\hat\mu^{\by}_{n, k}$ can be obtained from $\mu^{\by}_{n, k}$ by transferring the total $\by$-weights of their offspring to each black vertex.  
Let $\dgr^{\cF}$ denote the graph distance in $\cF$. Denote by $T^{(b)}_{n, (k)}$ the subset that contains all the black vertices  of $T_{n, (k)}$. 
Comparing~\eqref{id: ht} and~\eqref{id: tree-dist} with the definition~\eqref{def: dh} of $d_h$, we find that
\begin{equation}
\label{eq: ghp-bd1}
\dghp\Big(\big(T^{(b)}_{n, (k)}, n^{-\frac13} \dgr^{\cF}, n^{-\frac23}\hat\mu^{\by}_{n,k}\big),  \cT_{2\mathrm{H}^{n, k}}\Big)\le 2n^{-\frac13}, \quad 1\le k\le \kappa_{n,m}. 
\end{equation}
Meanwhile, since in $T_{n, (k)}$ each white vertex is at distance $1$ from a black vertex, we have 
\begin{equation}
\label{eq: ghp-bd2}
\dghp\Big(\big(T^{(b)}_{n, (k)}, n^{-\frac13}\dgr^{\cF}, n^{-\frac23}\hat\mu^{\by}_{n,k}\big), \big(T_{n, (k)}, n^{-\frac13}\dgr^{\cF}, n^{-\frac23}\mu^{\by}_{n,k}\big)\Big)\le 2n^{-\frac13}. 
\end{equation}
Next, recall from~\eqref{def: cE''} the set $\cS_{n,m}$, which can be used to create surplus edges. These surplus edges may include self loops, which nevertheless do not affect graph distance. 
We extract from~$\cS_{n,m}$ those elements that belong to the interval $(g_{n, (k)}, d_{n, (k)})$: 
\[
\mathrm S^{n,k}=\big\{\big(n^{-\frac23}t_{n,l}, n^{-\frac23}t'_{n, l}\big): g_{n, (k)}<t_{n,l}<d_{n, (k)}, 1\le l\le q_{n,m}\big\}.
\]
Recall from~\eqref{def: graph'} the metric measured space $\mathfrak G(h, \varpi, \epsilon)$. 
Taking into account~\eqref{bd: shortcut-ghp},~\eqref{eq: ghp-bd1},~\eqref{eq: ghp-bd2} and Lemma~\ref{lem: dis}, we conclude that  
\begin{equation}
\label{eq: ghp-bd3}
\dghp\Big(\cG(2\mathrm{H}^{n, k}, \mathrm S^{n,k}, n^{-\frac13}), \big(C^{\bx}_{n,(k)}, n^{-\frac13} \dgr, n^{-\frac23}\mu^{\by}_{n,k}\big)\Big) \le (3s_{n,k}+5)n^{-\frac13}, 
\end{equation}
with $s_{n,k}$ standing for the cardinality of $\mathrm S^{n, k}$, $1\le k\le \kappa_{n,m}$. 

We also need an encoding of $C^{\bx}_{n, (k)}$ equipped with the measure $\mu^{\bx}_{n,k}$. To that end, recall from~\eqref{def: Sigma} the function $\Sigma^{n,m}$, and recall from Lemma~\ref{lem: sigma-p} that $\Sigma_{n,m}$ is strictly increasing on $(g_{n, (k)}, d_{n, (k)})$. Denote by $f^{n,k}$ the inverse function of $t\mapsto \Sigma^{n,m}(t+g_{n, (k)})$ defined on $0\le t\le d_{n, (k)}-g_{n, (k)}$. Note that the domain of $f^{n,k}$ is $[0, \bx(C^{\bx}_{n, (k)})]$. 
Set
\[
\hat{\mathrm H}^{n,k}(s) = \mathrm H^{n,k}\big(n^{-\frac23} f^{n,k}(n^{\frac23}s)\big), \quad 0\le s\le n^{-\frac23}\bx(C^{\bx}_{n, (k)}). 
\]
Since $\Sigma^{n,m}$ transfers the $\by$-weights to $\bx$-weights in the LIFO-order, we have in fact
\[
\big(T^{(b)}_{n, (k)}, n^{-\frac13}\dgr^{\cF}, n^{-\frac23}\mu^{\bx}_{n,k}\big) \text{ isometric to } \cT_{2\hat{\mathrm H}^{n, k}}.
\]
Following a similar argument as previously, we find that
\begin{equation}
\label{eq: ghp-bd4}
\dghp\Big(\cG(2\hat{\mathrm{H}}^{n, k}, \mathrm S^{n,k}, n^{-\frac13}), \big(C^{\bx}_{n,(k)}, n^{-\frac13} \dgr, n^{-\frac23}\mu^{\bx}_{n,k}\big)\Big) \le (3s_{n,k}+5)n^{-\frac13}. 
\end{equation}

\noindent
{\it Proof of Theorem~\ref{thm: thm1} and Proposition~\ref{prop: rnk} in the double finite third moment case. }\\
Recall from~\eqref{def: cC1} that the limit graph $(\cC^{(1)}_k, d^{(1)}_k, \mu^{(1)}_k)$ is defined as $\mathfrak G(\cH^{1, k}, \varpi_k, 0)$. Comparing this with~\eqref{eq: ghp-bd4} and noting~\eqref{eq: ghp}, we see that the convergence in~\eqref{eqcv: thm1} will follow once we have the convergence of $\hat{\mathrm H}^{n,k}$ to $\cH^{1, k}$, and of $\mathrm S^{n,k}$ to $\varpi_k$. The former is a combined consequence of the convergence of $H^{n,m}$ in~\eqref{cv: ZH1}, the convergence of $\Sigma^{n,m}$ in~\eqref{eqcv: Sigma-br}, Lemma~\ref{lem: cv-comp} and the convergence of $(g_{n, (k)}, d_{n, (k)})$ in~\eqref{cv: ex-int}. The latter is a consequence of Proposition~\ref{prop: cv-surplus} combined with~\eqref{cv: ex-int}. This proves Theorem~\ref{thm: thm1}. The first statement in Proposition~\ref{prop: rnk} has been shown in Lemma~\ref{lem: rkcons}. We then follow a similar argument as previously to prove the scaling limit of the connected components, replacing~\eqref{eq: ghp-bd4} with~\eqref{eq: ghp-bd3} along the way.

\begin{proof}[Proof of Theorems~\ref{thm: thm2} and~\ref{thm: thm3} and rest of Proposition~\ref{prop: rnk}]
This is very similar to the previous proof; we therefore omit the detail. 
\end{proof}

\begin{appendix}

\section{Convergence of the Laplace exponents}

Let $F$ be the cumulative distribution function of a probability measure supported on $(0, \infty)$. Define
\[
\varphi(\lambda) = \int_{(0, \infty)}(e^{-\lambda x}-1+\lambda x) x\, dF(x), \quad \lambda \ge 0. 
\]

\begin{lem}
\label{lem-cv: psi-br}
Assume that $\sigma_3(F):=\int_{(0,\infty)} x^3 dF(x)<\infty$. Then for each $\lambda\ge 0$, 
\[
n^{\frac23}\varphi(n^{-\frac13}\lambda) \xrightarrow{n\to\infty} \tfrac12 \sigma_3(F)\lambda^2.
\]
Moreover, for each $\alpha\in (1, 2)$ and $\lambda\ge 0$, we have
\[
n^{\frac{\alpha}{\alpha+1}}\varphi(n^{-\frac{1}{\alpha+1}}\lambda) \xrightarrow{n\to\infty} 0.
\] 
\end{lem}

\begin{proof}
Let us denote $\phi(x) = e^{-x}-1+x$. We note that $\phi(x) = \frac12 x^2 + r(x) x^3$, where $r$ is a continuous and therefore bounded function on $[0, x_0]$ for any fixed $x_0\ge 0$. 
We can write 
\begin{align*}
&n^{\frac23}\varphi(n^{-\frac13}\lambda)  = n^{\frac23}\int_{(0, \infty)}\phi(n^{-\frac13}\lambda x)x\, dF(x)\\
& = \frac12 \lambda^2\int_{(0, n^{\frac13}]}\!\! x^3 dF(x) + n^{-\frac13}\lambda^3\int_{(0, n^{\frac13}]} \!\!r(n^{-\frac13}\lambda x)x^4 dF(x) + n^{\frac23}\int_{(n^{\frac13}, \infty)}\!\!\phi(n^{-\frac13}\lambda x)x dF(x)\\
& =:I_1+I_2+I_3.
\end{align*}
Since $\sigma_3(F)<\infty$, we have $\int_{[M, \infty)}x^3 dF(x)=o(1)$ as $M\to\infty$. As a result, we deduce that for each $\lambda\ge 0$, 
\begin{equation}
\label{limitI1}
\frac12 \sigma_3(F)\lambda^2-I_1 = \frac12 \lambda^2\int_{(n^{\frac13},\infty)} x^3 dF(x) \to 0, \quad \text{as } n\to\infty.
\end{equation}
Using $0\le \phi(x)\le \frac12 x^2$ for all $x\ge 0$, we find that 
\begin{equation}
\label{limitI3}
0\le I_3\le \frac12\lambda^2 \int_{(n^{\frac13},\infty)}x^3 dF(x) \to 0, \quad n\to\infty. 
\end{equation}
For $I_2$, we first note that $C_r:=\sup_{x\in[0, n^{1/3}]}|r(n^{-\frac13}\lambda x)|<\infty$ for each fixed $\lambda$. 
Take some $\epsilon>0$ and let us write
\begin{align*}
\int_{(0, n^{\frac13}]} \!\! x^4 dF(x) &= \int_{(0, \epsilon n^{\frac13}]} \!\! x^4 dF(x) + \int_{(\epsilon n^{\frac13}, n^{\frac13}]} \!\! x^4 dF(x)\le \epsilon n^{\frac13} \sigma_3(F) + n^{\frac13}\int_{(\epsilon n^{\frac13}, \infty)}\!\! \!\! x^3 dF(x)\\
& = \epsilon n^{\frac13} \sigma_3(F) + o(n^{\frac13}), \quad n\to\infty.
\end{align*}
Since $\epsilon$ can be taken arbitrarily small, this shows that $\int_{(0, n^{1/3}]} x^4 dF(x) = o(n^{\frac13})$. Combined with the bound for $r(x)$, this implies that 
\[
|I_2|\le \lambda^3 C_r \cdot n^{-\frac13}\int_{(0, n^{\frac13}]} x^4 dF(x) \to 0, \quad n\to\infty. 
\]
Together with~\eqref{limitI1} and~\eqref{limitI3}, this proves the first convergence in the statement. For the second one, we take $n_k= \lfloor k^{3/(\alpha+1)}\rfloor$ to find that 
\[
(n_k)^{-\frac13}\sim k^{-\frac{1}{\alpha+1}} \quad\text{and}\quad 0\le \varphi(k^{-\frac{1}{\alpha+1}}\lambda)\le \varphi(n_k^{-\frac13}\lambda) = O(n_k^{-\frac23})=O(k^{-\frac{2}{\alpha+1}}). 
\]
The conclusion follows as $\alpha<2$. 
\end{proof}

\begin{lem}
\label{lem-cv: psi-power}
Assume that $1-F(x)\sim C_F x^{-1-\gamma}$ as $x\to\infty$ for some $C_F\in (0,\infty)$ and $\gamma\in (1, 2)$. Then for each $\lambda\ge 0$, 
\[
n^{\frac{\gamma}{\gamma+1}}\varphi(n^{-\frac{1}{\gamma+1}}\lambda) \xrightarrow{n\to\infty} C_F\tfrac{(\gamma+1)\Gamma(2-\gamma)}{\gamma(\gamma-1)} \lambda^{\gamma} .
\]
Moreover, for each $\alpha\in (1, \gamma)$ and $\lambda\ge 0$, we have 
\[
n^{\frac{\alpha}{\alpha+1}}\varphi(n^{-\frac{1}{\alpha+1}}\lambda) \xrightarrow{n\to\infty} 0.
\]
\end{lem}

\begin{proof}
Let us denote $\psi(x) = x(e^{- x}-1+x)$. Then we have
\[
\varphi(\lambda) = \lambda^{-1}\int_{(0,\infty)}\!\!\psi(\lambda x) dF(x) = \int_{(0,\infty)}\int_0^x \psi'(\lambda u)du dF(x) = \int_0^{\infty}\!\! \psi'(\lambda u) (1-F(u))du,
\]
where we have used Fubini's Theorem in the last identity. 
A change of variable then yields 
\[
n^{\frac{\gamma}{\gamma+1}}\varphi(n^{-\frac{1}{\gamma+1}}\lambda) = n\int_0^{\infty} \psi'(\lambda y) \Big(1-F\big(n^{\frac{1}{\gamma+1}}y\big)\Big)dy.
\]
Let $\epsilon>0$. The assumption of $F$ implies that for all $y\ge \epsilon$, we have 
\[
1-F(n^{\frac{1}{\gamma+1}}y) = C_F n^{-1}y^{-\gamma-1}(1+o(1)), \quad n\to\infty,
\]
where the $o(1)$-term is uniform for all $y\ge \epsilon$. It follows that
\begin{equation}
\label{epi-limit}
n\int_{\epsilon}^{\infty} \psi'(\lambda y) \Big(1-F\big(n^{\frac{1}{\gamma+1}}y\big)\Big) dy\xrightarrow{n\to\infty} C_F\int_{\epsilon}^{\infty} \psi'(\lambda y) y^{-\gamma-1}dy.
\end{equation}
An integration by parts yields that
\begin{align*}
C_F\int_{\epsilon}^{\infty}\!\! \psi'(\lambda y) y^{-\gamma-1}dy  = &\, C_F(\gamma+1)\int_{\epsilon}^{\infty}\!\! (e^{-\lambda y}-1+\lambda y)y^{-\gamma-1}dy -C_F\lambda^{-1}\epsilon^{-\gamma-1}\psi(\lambda \epsilon)\\
 \xrightarrow{\epsilon\to 0+} &\, C_F(\gamma+1)\int_{0}^{\infty}\!\! (e^{-\lambda y}-1+\lambda y)y^{-\gamma-1}dy = C_F\tfrac{(\gamma+1)\Gamma(2-\gamma)}{\gamma(\gamma-1)}\lambda^{\gamma},
\end{align*}
where we have used the fact that $\psi(x)\sim \tfrac12 x^3$ as $x\to 0+$. Noting that $\psi'(x) = (e^{-x}-1+x)+x(1-e^{-x})\ge 0$ for all $x\ge 0$, we can apply the monotone convergence theorem and take $\epsilon\to 0$ 
in~\eqref{epi-limit}.  This proves the first convergence in the statement. For the second one, let $n_k= \lfloor k^{(\gamma+1)/(\alpha+1)}\rfloor$, so that 
\[
(n_k)^{-\frac{1}{\gamma+1}}\sim k^{-\frac{1}{\alpha+1}} \quad\text{and}\quad 0\le \varphi(k^{-\frac{1}{\alpha+1}}\lambda)\le \varphi(n_k^{-\frac{1}{\gamma+1}}\lambda) = O(n_k^{-\frac{\gamma}{\gamma+1}})=O(k^{-\frac{\gamma}{\alpha+1}}). 
\]
The conclusion follows as $\alpha<\gamma$. 
\end{proof}

\begin{lem}
\label{lem: psi-bound}
Assume that $1-F(x)\sim C_F x^{-1-\gamma}$ as $x\to\infty$ for some $C_F\in (0,\infty)$ and $\gamma\in (1, 2)$. Then for each $\lambda_0\ge 0$,  we can find some $C=C(\lambda_0)\in (0,\infty)$ so that  
\[
\forall\, \lambda\in [0, \lambda_0]: \quad \varphi(\lambda) \ge C\lambda^{\gamma}. 
\]
\end{lem}

\begin{proof}
As in the beginning of the previous proof, we apply Fubini's Theorem to obtain
\[
\varphi(\lambda) = \int_0^{\infty} \psi'(\lambda u) (1-F(u))du, 
\]
where $\psi(x) =x(e^{-x}-1+x)$. We note that $\psi'(x)= e^{-x}-1+x + x(1-e^{-x})\ge x(1-e^{-x})\ge 0$ for all $x\ge 0$. Meanwhile, the assumption on $F$ implies that we can find some $C'\in (0, \infty)$ so that $1-F(x) \ge C'x^{-1-\gamma}$ for all $x\ge 1$. It follows that 
\begin{align*}
\varphi(\lambda)&\ge \int_1^{\infty} \!\!\psi'(\lambda u)(1-F(u))du \ge  C'\lambda\int_1^{\infty} \!\!(1-e^{-\lambda u}) u^{-\gamma}du
 = C'\lambda^{\gamma} \int_{\lambda}^{\infty} \!\!(1-e^{-y}) y^{-\gamma}dy\\
&\ge C'\lambda^{\gamma} \int_{\lambda_0}^{\infty} (1-e^{-y}) y^{-\gamma}dy.
\end{align*}
This implies the desired inequality, since the last integral is finite. 
\end{proof}

\section{A Girsanov-type theorem for L\'evy processes}
\label{sec: Z-prop}

We follow Appendix A of Conchon-Kerjan and Goldschmidt~\cite{CG23}.  
Let $(L_t)_{t\ge 0}$ be a spectrally positive L\'evy process with Laplace exponent $\Psi$, namely,
\[
\mathbb E[e^{-\lambda L_t}] = \exp\big(t\Psi(\lambda)\big), \quad \lambda\ge 0, t\ge 0,
\] 
In~\cite{CG23}, further conditions were assumed on $\Psi$. However, a careful reading reveals that the arguments there can be readily adapted to all spectrally positive L\'evy processes. 
For each $q>0$, define
\[
\mathcal E_t = \exp\Big(-\int_0^t qs\,dL_s-\int_0^t \Psi(qs)ds\Big), \quad t\ge 0.
\]
Then according to Lemma A.1 in~\cite{CG23}, $(\mathcal E_t)_{t\ge 0}$ is a unit-mean positive martingale with respect to the natural filtration $(\mathfrak F_t)_{t\ge 0}$ of $(L_t)_{t\ge 0}$. This allows us to introduce a new probability measure  $\mathbf Q$ on the canonical space by setting  
\[
\frac{d\mathbf Q}{d\mathbb P}\Big|_{\mathfrak F_t} = \mathcal E_t, \quad t\ge 0. 
\]
Then Proposition A.2 in~\cite{CG23} says that
\begin{equation}
\label{id: lapQ}
\mathbb E_{\mathbf Q}[e^{-\lambda L_t}] = \exp\Big(\int_0^t \big(\Psi(\lambda + qs)-\Psi(qs)\big)ds\Big). 
\end{equation}
Let us note that under $\mathbf Q$, $(L_t)_{t\ge 0}$ still has independent increments, although the distributions of increments are no longer stationary. As a result of this independence, the law of $(L_t)_{t\ge 0}$ can still be determined from its marginal laws. Put another way,~\eqref{id: lapQ} characterises the law of $(L_t)_{t\ge 0}$ under $\mathbf Q$. 
We can also identify the Doob--Meyer decomposition of $(L_t)_{t\ge 0}$ under $\bQ$. 
To that end, let us assume that $\Psi$ takes the following form: 
\begin{equation}\label{def: Psi}
\Psi(\lambda) = \alpha_0\lambda+\tfrac12 \beta \lambda^2 +\int_{(0, \infty)} (e^{-\lambda x}-1+\lambda x) \pi(dx),
\end{equation}
where $\alpha_0\in \R_+$, $\beta\in\R_+$ and $\pi$ is a $\sigma$-finite measure on $(0,\infty)$ satisfying $\int_{(0,\infty)} (x\wedge x^2)\pi(dx)<\infty$. In particular, this covers the case $\Psi(\lambda)=\lambda^{\alpha}$ for $\alpha\in (1, 2)$, as it suffices to take $\alpha_0=\beta=0$ and $\pi(dx) = \alpha(\alpha-1)/\Gamma(2-\alpha) x^{-\alpha-1}dx$. 
Let $(\cM_t)_{t\ge 0}$ be a martingale with independent increments whose marginal laws are characterised by
\[
\mathbb E[e^{-\lambda \cM_t}] = \exp\Big(\tfrac{1}{2}\beta \lambda^2 t+ \int_0^t\int_{(0,\infty)} (e^{-\lambda x}-1+\lambda x) e^{-qxs} \pi(dx)ds\Big), \quad \lambda\ge 0, t\ge 0. 
\]
Then the Laplace transform of $\cM_t - \Psi(qt)/q$ coincides with the right-hand side of~\eqref{id: lapQ}. In other words, 
\[
\{L_t: t\ge 0\} \text{ under $\mathbf Q$ is distributed as } \big\{\cM_t - \tfrac1q \Psi(qt): t\ge 0\big\}. 
\]

\section{A continuous-time version of Duquesne--Le Gall's Theorem}
\label{sec: DuLGthm}

For each $n\in \N$, let $\pi_n$ be a finite measure supported on $(0, \infty)$ and
denote by $q_n=\pi_n(\R_+)$ its total mass. 
Suppose that $\{X^{(n)}_t+t : t\ge 0\}$ is a compound Poisson process with L\'evy measure $\pi_n$. 
We define its associated height process as follows:
\[
\mathrm H^{(n)}_t =\#\Big\{0\le s\le t: X^{(n)}_{s-}< \inf_{u\in [s, t]} X^{(n)}_{u}\Big\}, \quad t\ge 0.
\]
Let $\Psi$ be as in~\eqref{def: Psi} and satisfy the Grey condition: $\int^{\infty} (\Psi(\lambda))^{-1}d\lambda<\infty$. 
Let $(X_t)_{t\ge 0}$ be a spectrally positive L\'evy process with Laplace exponent $\Psi$ and let $(\mathrm H_t)_{t\ge 0}$ be its height process, defined as follows:
\[
\mathrm H_t =\lim_{\epsilon\to 0+} \frac{1}{\epsilon} \int_0^t \mathbf 1_{\{X_s\le \inf_{u\in [s, t]} X_u+\epsilon\}}ds,
\]
where the limit exists in probability, according to Lemma 1.2.1 and Theorem 1.4.3 in~\cite{DuLG02}; see also Section 3.1 in~\cite{DuLeG05}.

Let $W_n$ be a random variable with distribution $\mathbb P(W_n\in dx) = \pi_n(dx)/q_n$ and denote by $g_n$ the generating function for the Poisson$(q_nW_n)$ distribution. Namely, for $s\in [0, 1]$, 
\begin{equation}
\label{def: gn}
g_n(s) = \frac{1}{q_n} \int_{(0,\infty)} \sum_{k\ge 0}s^k\,\mathbb P\big(\text{Poisson}(q_nx)=k\big)\pi_n(dx)=\frac{1}{q_n}\int_{(0,\infty)}e^{(s-1)q_nx}\pi_n(dx),
\end{equation}
For $m\ge 1$, write $g_n^{\circ m}=g_n\circ g_n\circ \cdots \circ g_n$ for the $m$-th iterated composition of $g_n$.  

\begin{thm}
\label{thm: DuLG}
Let $a_n, b_n$ be two sequences of positive numbers satisfying $a_n\to\infty$, $b_n/a_n\to \infty$ and $b_n/a_n^2\to \beta_0\in[0,\infty)$. Assume that $(\pi_n)_{n\in \N}$ is uniformly integrable and as $n\to\infty$, $q_n\to q\in (0,\infty)$ and 
\[
\tfrac{1}{a_n} X^{(n)}_{b_n} \Longrightarrow X_1.
\]
Assume further that there exists some $\delta>0$ so that
\begin{equation}
\label{condH}
\limsup_{n\to\infty} a_n \big(1-g_n^{\circ\lfloor \delta b_n/a_n\rfloor} (0) \big)<\infty.
\end{equation}
Then we have the joint convergence: 
\begin{equation}
\label{cv: DuLG}
\Big\{\tfrac{1}{a_n} X^{(n)}_{b_n t}, \tfrac{a_n}{b_n} \mathrm H^{(n)}_{b_n t}: t\ge 0\Big\} \Longrightarrow \big\{X_t, \mathrm H_t: t\ge 0\big\} \quad\text{in } \mathbb D(\R_+, \R^2).
\end{equation}
Denote 
\begin{equation}
\label{def: Psin}
\Psi_n(u) = b_n \Big(g_n\big(1-\tfrac{u}{a_n}\big)-1 +\tfrac{u}{a_n}\Big), \quad 0\le u\le a_n.
\end{equation}
If we have
\begin{equation}
\label{cond: H}
\lim_{y\to\infty}\limsup_{n\to\infty} \int_y^{a_n} \frac{du}{\Psi_n(u)} = 0,
\end{equation}
then~\eqref{condH} holds for some $\delta>0$. 
\end{thm}

\begin{proof}
We follow closely the arguments in~\cite{BrDuWa21}. In the current set-up, we allow the L\'evy measure $\pi_n$ of $X^{(n)}$ to be a general finite measure on $(0,\infty)$, while in~\cite{BrDuWa21} it takes a specific form: $\pi_n(dx) = (\sum_i w_i)^{-1}\sum_{i\ge 1}w_i\delta_{w_i}(dx)$, where $(w_i)_{i\ge 1}$ is a finite sequence of positive numbers. The specific form of $\pi_n$ actually plays little role in the proof in~\cite{BrDuWa21}. On the other hand, the fact that we allow $\pi_n(\R_+)$ to be different from 1 requires some slight modifications in the proof. Overall, the arguments laid out here are mostly straightforward adaptations from~\cite{BrDuWa21}. We only outline the main steps and highlight the differences. 

\smallskip
\noindent
{\bf Step 1: Coding processes for a Bienaym\'e forest.}   Denote by $\tau_1<\tau_2<\tau_3<\cdots$ the successive jump times of $X^{(n)}$. We also write $\Delta_i = \Delta X^{(n)}_{\tau_i}$. Consider the LIFO-queue with (an infinite number of) clients arriving at $(\tau_i)_{i\ge 1}$ and requesting $\Delta_i$ service time, so that $X^{(n)}_t-\inf_{u\le t}X^{(n)}_u$ corresponds to the service load at time $t$.  Let $\cF=(\cT_i)_{i\ge 1}$ be the forest associated to the queue, where $\cT_i$ is the $i$-th tree component ranked in the arrival times of the customers. 
Thanks to the Markovian nature of $X^{(n)}$, each $\cT_i$ is an independent copy of a Bienaym\'e tree with a Poisson$(q_nW_n)$ offspring distribution. Let $(V^{(n)}_k)_{k\in\Z_+}$ and $(\Ht^{(n)}_k)_{k\in\Z_+}$ be the respective \L{}ukasiewicz path and height process of this forest, which is obtained by running a depth-first traversal of the forest and setting the increment $\Delta V^{(n)}_k+1$ (resp.~$\Ht^{(n)}_k$) to be the number of offspring (resp. height) of the $k$-th vertex in this traversal. We refer to Section 3.1 in~\cite{BrDuWa21} for a precise definition. Let also $(C^{(n)}_t)_{t\ge 0}$ be the contour process of the forest, obtained by tracking the height of an imaginary particle travelling at unit speed in a depth-first fashion in the forest. We again point to Section 3.1 in~\cite{BrDuWa21} for a definition of this classic notion. For $t\ge 0$, denote
\[
N(t) = \sum_{i\ge 1}\mathbf 1_{\{\tau_i\le t\}},\quad M(t) = \sum_{s\le t}\mathbf 1_{\{\Delta\rH^{(n)}_{s}\ne 0\}}.
\]
Note that $N(t)$ counts the number of jumps of $X^{(n)}$ up to time $t$, while $M(t)$ counts the same quantity for $\rH^{(n)}$. 
We have the following analogue of Lemma 3.1 in~\cite{BrDuWa21}. 

\begin{lem}
For each $t\ge 0$, conditional on $X^{(n)}_t-\inf_{u\le t}X^{(n)}_u$, $V^{(n)}_{N(t)}-\inf_{k\le N(t)}V^{(n)}_k$ has a Poisson distribution of mean $q_n(X^{(n)}_t-\inf_{u\le t}X^{(n)}_u)$; conditional on $-\inf_{u\le t}X^{(n)}_u$, $-\inf_{k\le N(t)}V^{(n)}_k$ has a Poisson distribution of mean $-q_n\inf_{u\le t}X^{(n)}_u$. In consequence, for all $a, x>0$, we have
\begin{equation}
\label{bd: XR}
\mathbb P\Big(\Big|V^{(n)}_{N(t)}-q_n X^{(n)}_t\Big|>2a\Big) \le 1\wedge \frac{4q_n x}{a^2} + \mathbb P\Big(-\inf_{u\le t}X^{(n)}_u\ge x\Big) + \mathbb E\Big[1\wedge \frac{q_n}{a^2}\big(X^{(n)}_t-\inf_{u\le t}X^{(n)}_u\big)\Big].
\end{equation}
\end{lem}

The proof of the lemma proceeds exactly as in the proof of Lemma 3.1 in~\cite{BrDuWa21}. Informally, to see why the first statement is true, we note that $V^{(n)}_{N(t)}-\inf_{k\le N(t)}V^{(n)}_k$ counts the number of those clients that will arrive after $t$ and appear as the offspring of some client currently in the queue. 
Meanwhile, the same arguments leading to Eq.~(95) and (96) in~\cite{BrDuWa21} yield in the current case that $C^{(n)}_{M(t)} = \rH^{(n)}_t$ for all $t\ge 0$ and for each $t, a\in (0,\infty)$,
\begin{equation}
\label{bd: M}
\mathbb P\Big(\sup_{s\le t}\big|M(s)-2q_ns\big| > 2a) \le 1\wedge\frac{16q_n t}{a^2} + \mathbb P\Big(1+\sup_{s\le t}\Ht^{(n)}_{N(s)}>a\Big).
\end{equation}

\smallskip
\noindent
{\bf Step 2: Convergence of $V^{(n)}$.} Let us show that under the conditions of the theorem, we have
\begin{equation}
\label{cv: Rn}
\Big\{\tfrac{1}{a_n} V^{(n)}_{\lfloor b_n t\rfloor} : t\ge 0\Big\} \Longrightarrow \big\{qX_{t/q}: t\ge 0\}.
\end{equation}
We follow the arguments in Section 7.1 of~\cite{BrDuWa21}. In the first instance, let us introduce a random walk $(S^{(n)}_k)_{k\ge 0}$, where $S^{(n)}_0=0$ and for $k\ge 1$, 
\[
\Delta S^{(n)}_k = q_n\Delta_k - 1.
\]
We further assume that $(\cN_k)_{k\ge 1}$ are i.i.d.~$\cN(0, 1)$ variables, independent of $(\Delta_k)_{k\ge 1}$. 
Combining the convergence of $X^{(n)}_{b_n}/a_n$ to $X_1$ with $q_n\to q$, we have
\[
\tfrac{q_n}{a_n} X^{(n)}_{b_n/q_n}\quad \Longrightarrow\quad qX_{1/q}.
\]
This is equivalent to the convergence of their Laplace exponents.  As $(\pi_n)_{n\in \N}$ is uniformly integrable, this is further equivalent to the following (Theorem 2.9, Chapter VII,~\cite{JaSh-book}):
\begin{align}\label{cvX1}
& \frac{b_n}{a_n} \Big(1- \int_{\R_+}x \pi_n(dx) \Big)\to \alpha_0, \\ \label{cvX2}
& \frac{b_n}{q_n} \int_{\R_+} \Big(1\wedge \frac{q_n^2 x^2}{a_n^2}\Big) \pi_n(dx) \to \beta q + \frac{1}{q}\int_{\R_+}(1\wedge q^2x^2) \pi(dx), \\ \label{cvX3}
& \forall\, h\in \cC_0: \quad \frac{b_n}{q_n}\int_{\R_+}h\Big(\frac{q_n x}{a_n}\Big) \pi_n(dx) \to \frac{1}{q}\int_{\R_+}h(qx)\pi(dx),
\end{align}
where $\cC_0$ denotes the set of all bounded continuous functions from $\R_+$ to $\R$ vanishing on a neighbourhood of $0$. 
Recall that $b_n/a_n^2\to\beta_0\in \R_+$. Combined with the central limit theorem, this implies the following distributional convergence as $n\to\infty$: 
\begin{equation}
\label{cv: bc}
\frac{1}{a_n} \sum_{k=1}^{\lfloor b_n\rfloor}\cN_k \; \Longrightarrow\; \sqrt{\beta_0}\,\cN_1. 
\end{equation}
Our aim is to show that
\begin{equation}
\label{pfcv: Sn}
\frac{1}{a_n}S^{(n)}_{\lfloor b_n\rfloor}+\frac{1}{a_n} \sum_{k=1}^{\lfloor b_n\rfloor}\cN_k\; \Longrightarrow\; q X_{1/q}\,,
\end{equation}
Applying Theorem 2.35 in Ch.~VII of~\cite{JaSh-book} and taking into account~\eqref{cv: bc}, we see that~\eqref{pfcv: Sn} will follow once we show that 
\begin{align}\label{cvS1}
& -\frac{b_n}{a_n} \mathbb E\big[S^{(n)}_1\big]  \to \alpha_0, \\ \label{cvS2}
& \frac{b_n}{a_n^2} \mathbb E\big[\big(S^{(n)}_1\wedge a_n)^2\big] \to \beta q-\beta_0 + \frac{1}{q}\int_{(0,\infty)}(1\wedge q^2x^2) \pi(dx), \\ \label{cvS3}
& \forall\, h\in \cC_0: \quad b_n\mathbb E\big[h(S^{(n)}_1/a_n)\big]\to \frac{1}{q}\int_{(0,\infty)}h(qx)\pi(dx).
\end{align}
We point out that~\eqref{cvS1} follows directly from~\eqref{cvX1}, as $\mathbb E[S^{(n)}_1]=\int x\pi_n(dx)-1$. For the convergence in~\eqref{cvS3}, from~\eqref{cvX3} and the fact that $h$ is bounded and continuous, we deduce 
\[
b_n\mathbb E\Big[h\Big(\frac{S^{(n)}_1}{a_n}\Big)\Big] 
 = \frac{b_n}{q_n}\int_{\R_+}h\Big(\frac{q_n x-1}{a_n}\Big)\pi_n(dx) \to \frac{1}{q}\int_{\R_+}h(qx)\pi(dx).
\]
It remains to check~\eqref{cvS2}. To that end, we note that if $Z$ is a real-valued random variable and $c\in \R_+$, we have $\mathbb E[Z\wedge c] = \mathbb E[Z]-\mathbb E[Z; Z\ge c]+c\mathbb P(Z\ge c)$. As $\Delta_1$ only takes positive values, $(q_n\Delta_1-1)^2\ge a_n^2$ implies $q_n\Delta_1\ge a_n+1$ for $n$ large enough. Then for those $n$, 
\begin{align*}
&\mathbb E\big[\big(S^{(n)}_1\wedge a_n)^2\big] - \mathbb E\big[(q_n\Delta_1)^2\wedge a_n^2\big] \\
=& -2\mathbb E[q_n\Delta_1]+1 -\mathbb E\big[(q_n\Delta_1-1)^2; q_n\Delta_1\ge a_n+1\big] + \mathbb E\big[(q_n\Delta_1)^2; q_n\Delta_1\ge a_n\big] \\
&-a_n^2 \mathbb P(a_n\le q_n\Delta_1<a_n+1)\\
= & -2\mathbb E[q_n\Delta_1]+1 +\mathbb E\big[(q_n\Delta_1)^2; a_n\le q_n\Delta_1< a_n+1\big] + 2\mathbb E[q_n\Delta_1; q_n\Delta_1\ge a_n+1] \\
&- \mathbb P(q_n\Delta_1\ge a_n+1) -a_n^2 \mathbb P(a_n\le q_n\Delta_1<a_n+1)\
\end{align*}
We deduce from~\eqref{cvX1} that $\mathbb E[q_n\Delta_1] = 1+O(a_n/b_n) = 1+o(1)$. The convergence in~\eqref{cvX3} implies that $\mathbb P(q_n\Delta_1\ge a_n) \sim \mathbb P(q_n\Delta_1\ge a_n+1) = O(1/b_n)$. It follows that 
\[
0\le \mathbb E\big[(q_n\Delta_1)^2; a_n\le q_n\Delta_1< a_n+1\big]\le (a_n+1)^2 \mathbb P(a_n\le q_n\Delta_1<a_n+1)=o(a_n^2/b_n), 
\]
as $n\to\infty$. Similarly, $a_n^2\mathbb P(a_n\le q_n\Delta_1<a_n+1)=o(a_n^2/b_n)$. 
Since $(\pi_n)_{n\ge 1}$ is uniformly integrable, $\mathbb E[q_n\Delta_1; q_n\Delta_1\ge a_n+1]=o(1)$. Finally, $\mathbb P(q_n\Delta_1\ge a_n+1) = o(1)$. All put together, we find that 
\[
\mathbb E\big[\big(S^{(n)}_1\wedge a_n)^2\big] - \mathbb E\big[(q_n\Delta_1)^2\wedge a_n^2\big]  = -1 +o(a_n^2/b_n), \quad n\to\infty,
\]
where we have used the fact that $o(1) = o(a_n^2/b_n)$, since $b_n/a_n^2\to\beta_0$ and is therefore bounded. Comparing this with~\eqref{cvX2}, we conclude that 
\[
\frac{b_n}{a_n^2} \mathbb E\big[\big(S^{(n)}_1\wedge a_n)^2\big]  = \beta q+\frac{1}{q}\int_{\R_+}(1\wedge q^2x^2)\pi(dx) - \frac{b_n}{a_n^2} +o(1),
\]
which yields~\eqref{cvS2}. 

Having proved~\eqref{pfcv: Sn}, let us write $L(\lambda) = \log\mathbb E[\exp(-\lambda S^{(n)}_1)]$; we obtain (Lemma A.3 in~\cite{BrDuWa21}) that for each $\lambda\ge 0$, 
\begin{equation}
\label{cvSla}
b_n L(\lambda/a_n)+ \frac12\frac{b_n\lambda^2}{a_n^2} \to \Psi(q\lambda)/q.
\end{equation}
Now let $\cL(\lambda)=\log\mathbb E[\exp(-\lambda V^{(n)}_1)]$. Noting that $V^{(n)}_1+1$ is a Poisson random variable of mean $q_n\Delta_1=S^{(n)}_1+1$, we find that
\[
\cL(\lambda) = \lambda +\log\mathbb E\big[\exp\big(-q_n\Delta_1(1-e^{-\lambda})\big)\big] = \lambda + L(1-e^{-\lambda})-(1-e^{-\lambda}). 
\] 
It then follows from~\eqref{cvSla} that $b_n\cL(\lambda/a_n)\to \Psi(q\lambda)/q$ for each $\lambda\ge 0$. The convergence in~\eqref{cv: Rn} now follows (Corollary 3.6, Ch.~VII,~\cite{JaSh-book}). 

\smallskip
\noindent
{\bf Step 3: Joint convergence of $X^{(n)}$ and $\mathrm H^{(n)}$.} Thanks to Theorem 2.5.2 and Corollary 2.5.1 in~\cite{DuLG02}, the convergence in~\eqref{cv: Rn} together with the condition~\eqref{condH} yields the following:
\begin{equation}
\label{cvjt}
\Big\{\tfrac{1}{a_n}V^{(n)}_{\lfloor b_n t\rfloor}, \tfrac{a_n}{b_n} \Ht^{(n)}_{\lfloor b_n t\rfloor}, \tfrac{a_n}{b_n} C^{(n)}_{b_n t}: t\ge 0\Big\} \Longrightarrow \big\{qX_{t/q}, \tilde\rH_t, \tilde\rH_{t/2}: t\ge 0\big\} 
\end{equation}
in $\mathbb D(\R_+, \R^3)$, where $(\tilde\rH_t)_{t\ge 0}$ is the height process for $(qX_{t/q})_{t\ge 0}$ and satisfies
\[
\tilde\rH_t = \lim_{\epsilon\to 0} \frac{1}{\epsilon} \int_0^t \mathbf 1_{\{qX_{s/q}\le \inf_{u\in [s, t]} qX_{u/q}+\epsilon\}}ds = \mathrm H_{t/q} \quad\text{ almost surely}. 
\]
We note that $N(t), t\ge 0$ is a Poisson process of rate $q_n$. It is then not difficult to show that for all $t\ge 0$,
\[
\sup_{s\le t}\big| \tfrac{1}{b_n} N(b_n t) - qt\big|  \xrightarrow{n\to\infty} 0 \quad\text{in probability}. 
\]
Combining this with~\eqref{bd: XR} and the convergence of $X^{(n)}$, we deduce that jointly with the convergence in~\eqref{cvjt}, we also have $(X^{(n)}_{t})_{t\ge 0}$ converging in distribution to $(X_{t})_{t\ge 0}$. Similarly, with~\eqref{bd: M} and the identity $C^{(n)}_{M(t)}=\rH^{(n)}_t$, we can further incorporate $\rH^{(n)}$ into the convergence to conclude that 
\begin{align*}
& \Big\{\tfrac{1}{a_n}V^{(n)}_{\lfloor b_n t\rfloor}, \tfrac{a_n}{b_n} \Ht^{(n)}_{\lfloor b_n t\rfloor}, \tfrac{a_n}{b_n} C^{(n)}_{b_n t}, \tfrac{a_n}{b_n}X^{(n)}_{b_nt}, \tfrac{a_n}{b_n}\rH^{(n)}_{b_nt}: t\ge 0\Big\} \\
& \qquad\qquad\qquad\qquad\qquad\qquad\qquad\qquad\qquad\qquad\Longrightarrow \big\{qX_{t/q}, \tilde\rH_t, \tilde\rH_{t/2}, X_{t}, \tilde\rH_{qt}: t\ge 0\big\} 
\end{align*}
in $\mathbb D(\R_+, \R^3)\times \mathbb D(\R_+, \R^2)$. As $(\tilde\rH_{qt})_{t\ge 0}=(\rH_t)_{t\ge 0}$ almost surely, 
the joint convergence in~\eqref{cv: DuLG} follows. Finally, to see why~\eqref{cond: H} is a sufficient condition for~\eqref{condH}, one can follow the same arguments given in the proof of Proposition 7.3 in~\cite{BrDuWa21}. 
\end{proof}

\end{appendix}
%
%
%
%

\noindent
{\bf Acknowledgments.}
We are very grateful to the anonymous referee for their comments and suggestions. 

\end{document}